\documentstyle[german,12pt,epsfig]{report}

\def\H{{\bf H}}
\def\C{{\bf C}}
\def\R{{\bf R}}
\def\F{{\bf F}}
\def\Z{{\bf Z}}
\def\N{{\bf N}}
\def\Q{{\bf Q}}

\def\pslz{{\rm PSL}_2(\Z)}
\def\slz{{\rm SL}_2(\Z)}

\def\qed{\hfill\framebox[4mm][t1]{\phantom{x}}}

\def\me#1{\ifx -#1\else{\ifx .#1\else{\ifx ,#1\else{\ifx :#1\else\ 
\fi}\fi}\fi}\fi #1}

\def\VM{V^{\natural}}
\def\VBE{{\natural}}
\def\VB{V\!B^{\,\VBE}}
\def\VF{V_{\rm Fermi}}
\def\V#1{V_{(#1)}}
\def\VSO#1{V_{{\bf SO}(#1)}}
\def\so#1{{\bf so}(#1)}
\def\SO#1{{\bf SO}(#1)}

\def\VY{{}^V\! Y}
\def\UY{{}^U\! Y}
\def\WY{{}^W\! Y}
\def\LY{{}^L\! Y}
\def\BY{{}^{\VB}\! Y}
\def\LVirY{{}^{L_{1/2}}\! Y}
\def\MY{{}^{\VM}\! Y}

\def\ou{\overline{u}}
\def\ov{\overline{v}}
\def\ow{\overline{w}}

\def\COM{\mathop{\rm Com}} 
\def\rang{{\rm Rang\ }}

\def\da{\downarrow}

\def\VOA#1{VOA\me{#1}}
\def\VOAs#1{VOAs\me{#1}}
\def\SVOA#1{SVOA\me{#1}}
\def\SVOAs#1{SVOAs\me{#1}}
\def\OVOA#1{(S)VOA\me{#1}}
\def\OVOAs#1{(S)VOAs\me{#1}}

\def\CFT#1{CFT\me{#1}}

\def\nett#1{"`sch"on"'\me{#1}}
\def\nette#1{"`sch"o\-ne"'\me{#1}}
\def\netten#1{"`sch"o\-nen"'\me{#1}}
\def\sehrnett#1{"`sehr sch"on"'\me{#1}}
\def\sehrnette#1{"`sehr sch"o\-ne"'\me{#1}}
\def\sehrnetten#1{"`sehr sch"o\-nen"'\me{#1}}
\def\sehrnetter#1{"`sehr sch"o\-ner"'\me{#1}}

\def\sde#1{selbst\-duale\me{#1}}
\def\sdn#1{selbst\-dualen\me{#1}}

\def\GT{\Gamma_{\theta}}
\def\GTN{\Gamma_{\theta,N}}
\def\jt{j_{\theta}}

\def\Lacht{L_{1/2}^{(48)}(0)}
\def\Lt{L_{1/2}^{\otimes 48}(0)}
\def\La{L_{1/2}(0)}
\def\Lb{L_{1/2}(\frac{1}{2})}
\def\Lc{L_{1/2}(\frac{1}{16})}
\def\L#1{L_{(#1)}}

\def\Lcliff#1{{\it L-Clif\-ford\/\me{#1}}}
\def\Llie#1{{\it L-Lie\/\me{#1}}}
\def\Lvir#1{{\it L-Vir\/\me{#1}}}

\def\badset{\{10,11,\frac{25}{2},13,\frac{27}{2},\frac{29}{2}\} }

\def\tx#1{\bf #1}

\title{Selbstduale Vertexoperatorsuperalgebren
\newline \hfill \phantom{xxxxxxxxx} und \newline
\hfill \phantom{}das Babymonster
}
\author{Gerald H"ohn}
\date{28.~Juni 1995
\vskip3cm
Dissertation Bonn
}

\begin{document}

\bibliographystyle{amsalpha}


\newtheorem{satz}{Satz}[section]
\newtheorem{lemma}[satz]{Lemma}
\newtheorem{korollar}[satz]{Korollar}
\newtheorem{definition}[satz]{Definition}
\newtheorem{vermutung}[satz]{Vermutung}

\renewcommand{\baselinestretch}{1.2}

\pagenumbering{roman}

\phantom{x}
\vspace{3cm}
\LARGE
\centerline{Selbstduale Vertexoperatorsuperalgebren}
\centerline{und}
\centerline{das Babymonster}

\vspace{5cm}

\normalsize
\centerline{Inaugural-Dissertation zur Erlangung des Doktorgrades}
\centerline{der Mathematisch-Naturwissenschaftlichen Fakult"at}
\centerline{der Rheinischen Friedrich-Wilhelms-Universit"at zu Bonn}

\vspace{2.5cm}
\centerline{vorgelegt von}

\medskip
\centerline{\sc Gerald H"ohn}

\vspace{2cm}

\centerline{Bonn, Juli 1995}

\newpage

\setcounter{bottomnumber}{2}

\phantom{x}
\vspace{18cm}
\begin{tabular}{ll}
Referent: & Prof.~Dr.~{\sc Yuri~Manin} \\
Koreferent: & Prof.~Dr.~{\sc Friedrich~Hirzebruch}
\end{tabular}

\newpage

\renewcommand{\baselinestretch}{1.0}\small\normalsize

\phantom{xxxx}
\vspace{16cm}

\begin{picture}(0,0)%
\epsfig{file=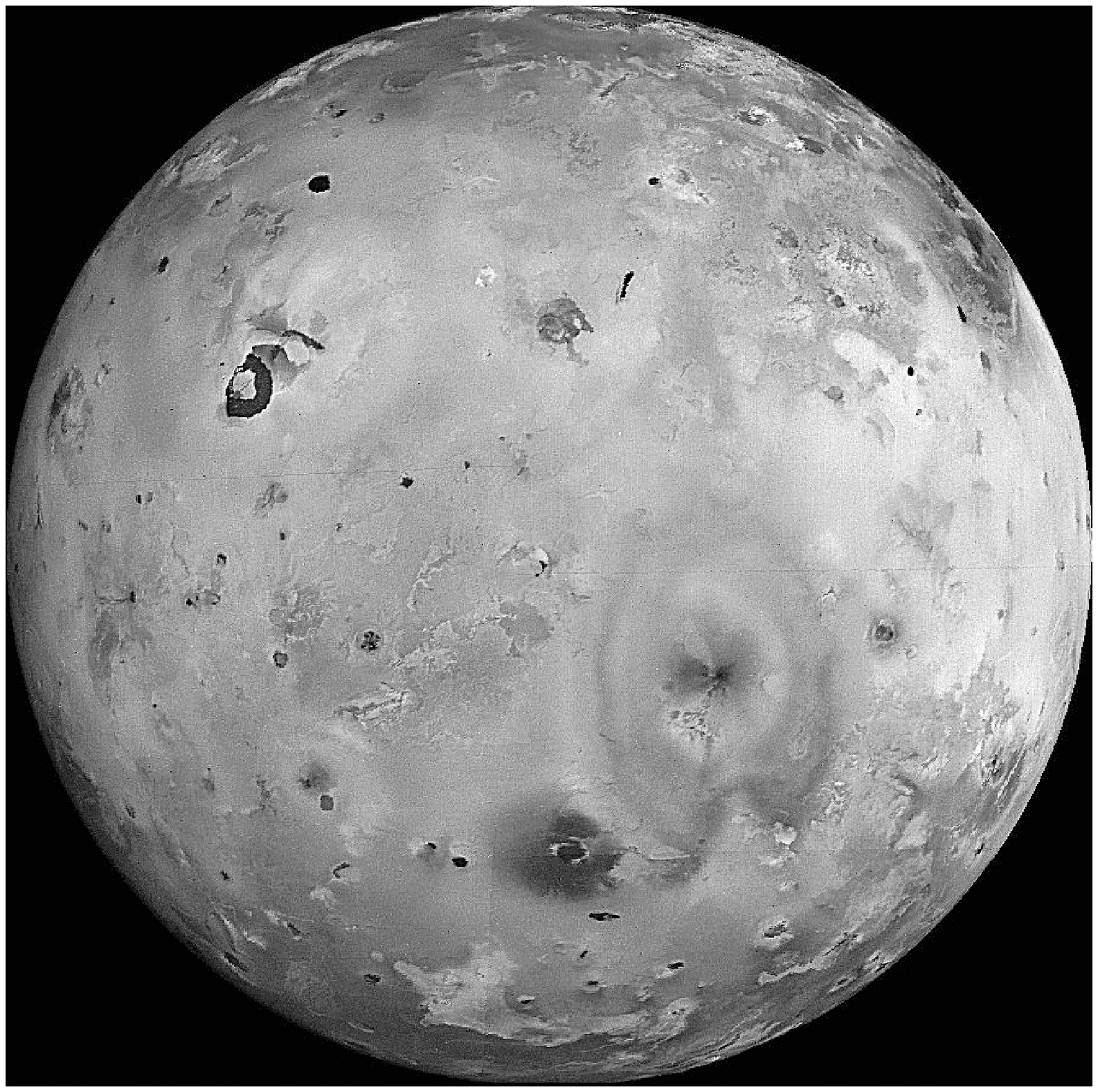}%
\end{picture}%

\vspace{2cm}
\centerline{\bf Jupiter Mond Io}

{\sl
Die Daten dieses von der Raumsonde {\em Voyager II}
aufgenommenen Bildes erreichten nach etwa~$45$
Minuten die Erde. Die "Ubertragung 
erfolgte mit Hilfe des bin"aren Golaycodes~$g_{24}$.
}

\newpage\phantom{x}

{\bf Danksagung }

\bigskip
Es ergaben sich aufgrund der vielen Querverbindungen der Theorie der 
Vertexoperator-Algebren
zu anderen Gebieten anregende Diskussionen mit verschiedenen Mathematikern.

Herrn Prof.~{\sc D.~Zagier} danke ich f"ur einige kl"arende Diskussionen 
"uber Modul\-funktionen, die zu einer Ver\-einfachung zweier Abschnitte f"uhrten. 
Viele Diskussionen mit {\sc G.~Mersmann}
halfen mir oft, meinen Ideen eine klare
Gestalt zu geben. Hierf"ur und die Hilfe beim Beweis eines Lemmas sei ihm gedankt.
Herrn {\sc W.~Neutsch} verdanke ich mehrere anregende Diskussionen "uber
Idempotente in der Griess Algebra. 
In der Endphase der Arbeit luden mich {\sc J.~Lepowsky} und 
{\sc Y.~Huang} nach Rutgers ein, was mir erstmals die Gelegenheit gab
"uber alle Aspekte der Arbeit zu diskutieren. Schlie"slich halfen mir 
{\sc G.~Mason}, {\sc C.~Dong} und {\sc R.~Griess} einige Unklarheiten zu
beseitigen und Argumente zu vervollst"andigen.  
Das {\sc GAP-Team} in Aachen half bei einer gruppen\-theoretischen
Charakter\-rechnung
und {\sc M.~Schr"oder} unterzog Teile der Arbeit einer kritischen Durchsicht.
Ich m"ochte mich bei allen Genannten f"ur ihre Hilfe sehr bedanken.
Herrn Prof.~{\sc P.~Slodowy} danke ich f"ur sein Interesse an der Arbeit und 
Herrn Prof.~{\sc Y.~Manin} f"ur Vortr"age "uber 
Vertexoperator-Algebren am MIT und am MPI f"ur Mathematik, die mein 
Interesse an diesem Gebiet wachgehalten haben.   

Mein besonderer Dank gilt Prof.~{\sc F.~Hirzebruch} f"ur seine
Unterst"utzung w"ahrend der Promotion,
f"ur die M"oglichkeit, am Max-Planck-Institut f"ur Mathematik 
arbeiten zu k"onnen
und vor allem f"ur die vielen anregenden Vorlesungen, die 
meine mathematische Entwicklung wesentlich beeinflu"st haben und das Fundament
zu dieser Arbeit legten.
 
Meinen {\sc Eltern} sei schlie"slich f"ur die stets begleitende Unterst"utzung
w"ahrend meines Studiums gedankt.

\newpage
\phantom{x}
\vspace{20cm}
{\footnotesize  
Fassung vom 25.3.1996}
\newpage

\pagenumbering{arabic}

\setcounter{bottomnumber}{1}
\addcontentsline{toc}{chapter}{Inhalt}
\vbox{\tableofcontents}

\newpage\phantom{leer}\newpage
\renewcommand{\baselinestretch}{1.2}\small\normalsize

\chapter*{
Einleitung}
\addcontentsline{toc}{chapter}{Einleitung}

\vspace{-5mm}
Ende der 70-er Jahre wurde ein \"uberraschender und zun\"achst
unerkl\"arlich erscheinender Zusammenhang zwischen
Fourier\-koeffizienten von Modul\-funktionen und Dimensionen irreduzibler
Darstellungen der gr\"o\ss{}ten sporadischen einfachen Gruppe, dem Monster,
gefunden.
Eine teilweise Erkl\"arung fand dieses "`moonshine"' genannte Ph\"anomen
Mitte der 80-er Jahre mit Hilfe der Vertexoperator-Algebren. 
In der vorliegenden Arbeit werden 
{\it selbstduale Vertex\-operator-Super\-algebren} unter\-sucht, 
dies sind Vertex\-operator Super\-algebren mit genau
einer irre\-duziblen Dar\-stellung. 
Die Vertex\-operator Super\-algebren haben neben der Theorie der 
Modul\-funktionen und endlichen Gruppen Verbindungen zu vielen weiteren 
Strukturen in Mathematik und Physik, wie affinen Kac-Moody Algebren, 
elliptischen Geschlechtern, Operaden
und konformen Quantenfeld\-theorien. Ein Leit\-gedanke der Arbeit 
ist die Analogie  zwischen {\it Codes,  Gittern und Vertex\-operator 
Algebren}, die auch den 
folgenden Satz motivierte, eines der Haupt\-resultate der Arbeit:

{\bf Satz:} {\sl Es existiert eine 
Vertexoperator-Superalgebra $\VB$ vom Rang $23\frac{1}{2}$ auf der das 
Babymonster $B$ --- die zweitgr\"o\ss{}te  spora\-dische end\-liche
ein\-fache Grup\-pe ---
in nat\"urlicher Weise operiert. Sie besitzt den Charakter 
$$\chi_{\VB}=\chi_{1/2}^{47}-47\,\chi_{1/2}^{23}=
q^{-\frac{47}{48}}\,(1+4371\,q^{3/2}+96256\, q^2+ 1143745\,q^{5/2}+ \cdots),$$ 
wobei $\chi_{1/2}=
\sqrt{\sum_{n \in {\bf Z}}q^{\frac{1}{2} n^2}/
(q^{1/24}\,\prod_{n=1}^{\infty}(1-q^n))}$.}

\medskip

Diese Babymonster Vertexoperator-Superalgebra ist das Analogon zum k"urzeren Golay Code bzw.~zum k"urzeren Leechgitter. 
Das Hauptanliegen der Arbeit ist aber allgemeiner und besteht in der Entwicklung
einer neuen Methode zur systematischen Klassifikation der selbstdualen 
Vertexoperator-Superalgebren.
Das oben angegebene Resultat ordnet sich so in das 
in diesem Zusammenhang erzielte Klassifikationsresultat
f"ur die spezielle Klasse der {\it extremalen\/} selbstdualen 
Vertexoperator-Superalgebren ein.


\medskip
Wir geben in der Einleitung 
zuerst eine kurze Einf"uhrung in die verschiedenen
Aspekte der Theorie und erl"autern den von 
uns gew"ahlten Zugang. 
In einem erg"anzenden Abschnitt werden einige
der Zusammenh"ange zu anderen Gebieten beschrieben,
die diese Arbeit motiviert haben. Dort werden auch
einige der zugeh"origen Begiffe erl"autert.
Schlie"slich geben wir eine Zusammenfassung "uber den Aufbau der Arbeit und
der erzielten Resultate. \nopagebreak[2]
Der vor allem hieran interessierte Leser sollte gleich 
mit diesem Abschnitt beginnen.

\pagebreak[3]
\medskip
{\bf Einf"uhrung in die Theorie der VOAs und die vorliegende Arbeit}

Die Theorie der Vertexoperator-Algebren (VOAs) wurde vor etwa zehn Jahren 
von R.~Bor\-cherds \cite{Bo-ur} sowie I.~Frenkel, 
J.~Lepowsky und A.~Meurmann \cite{FLM} 
eingef"uhrt, um einen Mond\-schein\-modul $\VM$ f"ur das Monster~\cite{Gr}
 --- der gr"o"sten spo\-radischen end\-lichen ein\-fachen Gruppe --- 
zu konstruieren und damit
Teile des als {\it Monstrous Moon\-shine} bekannten Ph"anomens~\cite{CoNo}
zu verstehen. Gleich\-zeitig sind \VOAs die mathe\-matische
Pr"a\-zisierung des physika\-lischen Begriffs der Chiralen Algebra
einer konformen Quantenfeldtheorie~\cite{BPZ}.

Die wesentlichen Daten einer Vertexoperator-Superalgebra (\SVOA{}) sind ein 
graduierter komplexer Vektorraum $V=\bigoplus_{n\in \frac{1}{2}
\Z_{\geq 0}} V_n$
mit endlichdimensionalen $V_n$, zusammen mit einer linearen
Abbildung  $Y(.,z):V\rightarrow {\rm End} (V)[[z,z^{-1}]]$ von $V$ in den 
Vektorraum der formalen Laurentreihen mit Koeffizienten in den Endomorphismen
von $V$ und zwei ausgezeichneten Elementen ${\bf 1}\in V_0$ und $\omega
\in V_2$. Die Reihen $Y(v,z)$ f"ur Elemente $v$ aus $V$ werden als 
Vertexoperatoren bezeichnet.

Zus"atzlich zu einer Reihe von Endlichkeitsbedingungen gelten zwei 
Hauptaxiome. Das eine ist die {\it Kommutativit"at\/} der Vertexoperatoren:
$$ Y(u,z_1) Y(v,z_2)\sim\pm\, Y(v,z_2)Y(u,z_1), $$
wobei die Tilde $\sim$ so zu verstehen ist, da"s es eine nat"urliche 
Zahl $n$ gibt,
so da"s nach Multiplikation beider Seiten mit $(z_1-z_2)^n$ eine Gleichheit
zwischen formalen Laurentreihen besteht.
Die technische Schwierigkeit hier ist, 
da"s Reihenentwicklungen von  $(z_1-z_2)^{-n}$ in verschieden Gebieten von 
$\C^2$ verglichen werden m"ussen.
Das andere Hauptaxiom betrifft das Element $\omega$. Die Koeffizienten 
des Vertexoperators
 $Y(\omega,z)=\sum_{n \in \Z} L_n z^{-n-2}$ sollen eine Darstellung
der {\it Virasoroalgebra\/} bilden:
$$[L_n,L_m]=(m-n)L_{m+n}+\frac{1}{12}(m^3-m)\delta_{m+n}c\cdot {\rm id}_V.$$
Die reelle Zahl $c$ wird als der {\it Rang\/} der VOA bezeichnet.
Die Operation der Virasoro\-algebra mu"s mit der Vertexoperator Abbildung 
$Y(.,z)$ vertr"aglich sein, was durch die Beziehung
$$\frac{d}{dz}Y(v,z)=Y(L_{-1} v,z), \quad\hbox{f"ur $v\in V$}$$
gefordert wird.

\medskip
Der Schwerpunkt der bisherigen mathemati\-schen Forschung auf diesem Gebiet lag 
einer\-seits in der Kon\-struktion spezieller Bei\-spiele, wie \VOAs
zu H"ochst\-gewichts\-dar\-stellungen von affinen Kac\--Moody
Al\-geb\-ren~\cite{FreZhu} oder der 
Virasoro\-algebra~\cite{Wan},
und andererseits in der
Entwicklung der allgemeinen Theorie, insbesondere einer Darstellungs\-theorie 
f"ur \VOAs sowie m"oglichen Verall\-gemeinerungen der Struktur wie 
Vertex\-operator-Super\-algebren (SVOAs) oder, noch allgemeiner, abelschen Inter\-twiner\-algebren~\cite{DoLe}.

Die allgemeine Theorie umfa"st die Entwicklung der Dar\-stellungs\-theo\-rie 
mit  den Begrif\-fen Modul, Inter\-twiner\-operator \cite{FHL} und 
Tensor\-produkt \cite{HuLe-tensor,Li-dr} sowie
 ei\-ne geo\-metri\-sche Inter\-pretation
mit Hilfe des Modul\-raumes punktierter Riemann\-scher Fl"achen~\cite{Hu-geo}
und durch Operaden~\cite{HuLe-operaden}.
Da sich Fl"achen h"oheren Geschlechtes
durch Zusammen\-kleben $3$-fach punktierter
Sph"aren erhalten lassen, sollten Eigen\-schaften
f"ur h"oheres Geschlecht, wie z.B.~in den Arbeiten \cite{Zhu-dr,Zhu-riem}, sich als Folgerungen ergeben. 

Wir werden in dieser Arbeit auf den geometrischen Standpunkt und den
Zusammenhang
zu Operaden nicht weiter eingehen. Das einzige wichtige Resultat, das 
wir aus der Theorie f"ur h"oheres Geschlecht ben"otigen, ist 
das Modultransformationsverhalten des {\it Charakters}
$$\chi_V=q^{-\frac{c}{24}}\sum_{n\in \Z} {\rm dim}\,V_n\cdot q^n$$
einer \VOA $V$ sowie der Charaktere ihrer Moduln~\cite{Zhu-dr} und die 
Verallgemeinerung f"ur \SVOAs.

Die nat"urliche Frage nach der Klassifikation von \SVOAs ist von
der mathematischen Seite bislang noch kaum untersucht worden.
F"ur die Vektoren in
$V_{1/2}$, $V_1$ oder $V_2$  einer \SVOA $V$ ist bekannt,
da"s sie unter gewissen Bedingungen \OVOAs erzeugen,
die zur unendlichdimensionalen Clifford, zu affinen Kac-Moody bzw.~Virasoro
Algebren assoziert sind.
Systematische Klassifkations\-ans"atze sind bisher allerdings
noch nicht entwickelt worden, zumindest nicht von der
mathematischen Seite.
Eine Klassifikation aller \VOAs wird auch nicht m"oglich sein: Jedes
gerade Gitter definiert eine VOA und schon die selbstdualen geraden Gitter
ab Rang $32$ sind nicht mehr explizit zu klassifizieren, da ihre
Anzahl f"ur h"ohere R"ange sehr stark anw"achst.

Andererseits wurde schon in~\cite{FLM} f"ur die Monster-VOA $\VM$ vom Rang $24$ vermutet, da"s sie die einzige selbstduale VOA vom Rang $24$ mit $V_1=0$ ist.
Diese --- bisher noch nicht bewiesene --- 
Eigenschaft w"are das Analogon entsprechender Aussagen f"ur
den Golay Code und das Leechgitter in Dimension $24$.
Ganz allgemein scheint eine tieferliegende, nicht restlos verstandene
Analogie zwischen {\it Codes}, {\it Gittern} und {\it VOAs} zu bestehen
(vgl.~\cite{Go-mero}), wenn auch
Codes und Gitter und ihre
Beziehung untereinander recht 
gut verstanden sind. 
Es existiert eine injektive Abbildung von der Menge
der Codes in die der Gitter sowie --- in einem n"achsten Schritt ---
eine injektive Abbildung von der Menge
der Gitter in die der \VOAs. 
Trotzdem sollte jede der drei Objektklassen auch f"ur sich untersucht werden, 
da jede der drei ihre eigene Hierarchie an Resultaten besitzt.

Schellekens~\cite{schellekens1} hat eine Liste von 
$71$ VOAs angegeben, die vermutlich alle selbstdualen \VOAs vom Rang $24$ 
beschreiben, allerdings ist nur f"ur einen Teil der Liste die volle
VOA-Struktur konstruiert, und die Klassifikationsmethode bedarf noch einer
mathematischen Rechtfertigung.
Dieses Resultat ist das Analogon der Klassifikation der selbstdualen geraden 
Gitter vom Rang $24$ durch Niemeier~\cite{niemeier,CoSl} und eines entsprechenden Resultates f"ur Codes ($9$ Codes,~\cite{MacSl}). 
Selbstduale ungerade Codes und selbstduale ungerade Gitter sind ebenfalls
bis zum Rang $24$ klassifiziert (s.~\cite{CoPless} f"ur Codes und 
\cite{CoSl-ul23},
Kapitel 16, f"ur Gitter). Ihre Klassifikation kann auf die Klassifikation der
geraden Codes und Gitter in Dimension $24$ zur"uckgef"uhrt 
werden.
Das Analogon zu {\it ungeraden\/} Codes und Gittern bilden 
Vertex\-operator-{\it Super\/}algebren. 
Ihre Klassifikation kann ganz analog auf die der 
selbstdualen \VOAs zur"uckgef"uhrt werden.

Als Folgerung des Modultrans\-formations\-verhaltens des Charakters 
l"a"st sich zeigen, da"s  SVOAs stets einen halbganzen Rang besitzen.
Bei Rang $23\frac{1}{2}$ findet sich als ein interessantes Beispiel die schon
zu Beginn erw"ahnte Babymonster-\SVOA $\VB$, 
auf der das von B.~Fischer gefundene und in~\cite{LeSi} konstruierte
Babymonster $B$ in nat"urlicher Weise operiert.
Diese Babymonster-\SVOA besitzt als Charakter die Modulfunktion
$$\chi_{\VB}=q^{-\frac{47}{48}}(1+4371\, q^{\frac{3}{2}}+96256\, q^2
 + 1143745\,{q^{{5\over 2}}} + 
  9646891\,{q^3} + 64680601\,{q^{{7\over 2}}}
  + \cdots \,),$$
deren Berechnung den Ausgangspunkt zu dieser Arbeit bildete.
Die \SVOA $\VB$ ist das Analogon des k"urzeren Golay Codes
$g_{22}$ in Dimension $22$ und des k"urzeren Leechgitters $O_{23}$ in
Dimension~$23$. Wir vermuten, da"s $\VB$ selbstdual ist.

\medskip

Die in dieser Arbeit entwickelte Klassi\-fikations\-methode f"ur selbst\-duale 
\OVOAs ver\-wendet Struk\-tur\-aussagen f"ur von Vek\-toren kleinen Gewichts
erzeugte \SVOAs, 
das Modul\-trans\-for\-mations\-verhalten der Charaktere sowie die Beziehung 
zwischen selbstdualen \SVOAs und \VOAs.
Einige der in den Resul\-taten der Arbeit gemachten technischen An\-nahmen
werden sich wohl durch weitere Unter\-suchungen als unn"otig erweisen.
Eine Hoffnung ist, da"s das Ver\-st"andnis von \SVOAs bis zum Rang $26$
n"utzlich f"ur Anwendungen sein wird, "ahnlich
wie z.B.~Cartans Klassi\-fikation der halb\-einfachen Lie\-algebren oder die Klassi\-fikation der selbst\-dualen Codes und Gitter in Dimen\-sionen bis $26$. 
Zumindest wird eine f"ur sich interessante Struktur sichtbar.

\bigskip

{\bf Zusammenh"ange zu anderen Gebieten}

Wir erl"autern in diesem Abschnitt kurz einige der Zusammenh"ange zu anderen Gebieten, die f"ur diese Arbeit von Bedeutung waren. Ansonsten sei auf die 
angegebene Literatur verwiesen, insbesondere auf die ausf"uhrliche
Einleitung in~\cite{FLM}.

\medskip
{\it Kombinatorik: Codes, Gitter und \VOAs{} }

Eine Einf"uhrung und Standardreferenz zur Codierungstheorie ist das 
Buch~\cite{MacSl}. F"ur Gitter und Kugelpackungen bietet~\cite{CoSl} eine
ziemlich vollst"andige "Ubersicht. Codes und Gitter, die ersten beiden 
Objektklassen der dreifachen Analogie zwischen Codes, Gittern und \VOAs, werden auch in~\cite{Eb} ausf"uhrlich untersucht.
Dieses Buch basiert auf einer
Vorlesung von Prof.~Hirzebruch, die ich im Wintersemester 86/87 in Bonn geh"ort habe.

Wir fassen hier kurz die wichtigsten Definitionen f"ur Codes und Gitter 
zusammen.

Ein (bin"arer linearer) Code $C$ der L"ange $n$ ist ein linearer Unterraum des 
Vektorraumes $\F_2^n$ "uber dem K"orper mit $2$ Elementen.
Auf $\F_2^n$ hat man das Standard\-skalar\-produkt 
$(c,d)=\sum_{i=1}^n c_i\cdot d_i$ f"ur Vektoren $c,d\in \F_2^n$.
Mit $C^{\perp}$ wird der zu $C$ bez"uglich $(.,.)$ orthogonale Code bezeichnet.
Ein Code $C$ hei"st {\it selbstdual\/}, falls $C=C^{\perp}$ gilt und 
{\it (doppelt) gerade\/}
falls das Gewicht $w(c)=\sum_{i=1}^n c_i\in \Z$ f"ur alle $c\in C$ durch
$4$ teilbar ist.
Das {\it Gewichts\-z"ahlerpolynom\/} von $C$ ist die erzeugende Funktion
$$W_C(x)=\sum_{c\in C}x^{n-w(c)}$$
f"ur die Anzahl der Codew"orter vom festem Gewicht.

Ein Gitter $L$ vom Rang $n$ ist eine diskrete Untergruppe von maximalen
Rang in dem Vektorraum $\R^n$, der die standard euklidischen Struktur
besitzt. 
Mit $L^{*}=\{x\in \R^n\mid (x,y)\in \Z\ \hbox{f"ur alle $y\in L$}\}$ 
wird das zu $L$ duale Gitter bezeichnet.
Das Gitter $L$ hei"st ganz, falls $L\subset L^*$ ist,
d.h.~das Skalarprodukt zwischen zwei beliebigen Gittervektoren eine ganze 
Zahl ist.
Ein Gitter $L$ hei"st {\it selbstdual\/} oder unimodular, falls $L=L^*$ gilt 
und {\it gerade\/},
falls die Quadratl"ange $(x,x)$ f"ur alle $x\in L$ durch
$2$ teilbar ist.
Die {\it Thetareihe\/} von $L$ ist die erzeugende Funktion
$$\Theta_L(q)=\sum_{x\in L}q^{\frac{1}{2}(x,x)}$$
f"ur die Anzahl der Gittervektoren von fester Quadratl"ange.
\label{KA}

Jedem  Code $C$ kann ein Gitter $L_C$
zugeordnet werden:
$$L_C:=\frac{1}{\sqrt{2}}\rho^{-1}(C),$$
wobei $\rho:\Z^n\longrightarrow \F_2^n$ die Reduktion modulo $2$ darstellt.
Ist $C$ selbstorthogonal, selbstdual bzw.~gerade, so ist $L_C$ ganz,
selbstdual bzw.~gerade. 

F"ur \SVOAs hat man analoge Definitionen: Eine \SVOA hei"st {\it selbstdual\/},
falls sie genau einen irreduziblen Modul besitzt. Eine gerade \SVOA ist eine
{\it \VOA{}\/} und der {\it Charakter\/} ist die erzeugende Funktion f"ur die 
Dimensionen der homogenen Komponenten.
Jedem ganzen Gitter $L$ kann eine \SVOA $V_L$ zugeordnet werden.

Einen weiteren vierten Schritt in der Analogie zwischen Codes, Gittern und 
\VOAs bilden Codes "uber der Kleinschen Vierergruppe $\Z_2\times\Z_2$
(s.~\cite{Ho-klein}).

\medskip

{\it Gruppentheorie: Endliche einfache Gruppen und Liealgebren}

Die kompakten einfachen zusammenh"angenden Lieschen Gruppen sind vor 
hundert Jahren unabh"angig von W.~Killing und E.~Cartan klassifiziert worden. 
Die Klassifikation der 
endlichen einfachen Gruppen kam vor etwa 15 Jahren zu einem Abschlu"s:
Die Liste der endlichen einfachen Gruppen besteht aus den zyklischen Gruppen 
von Primzahlordnung, den alternierenden Gruppen, den Gruppen von Lieschem Typ
sowie einer Liste von $26$ Ausnahmegruppen, den sogenannten
sporadischen\nopagebreak[4]
einfachen Gruppen~\cite{atlas,Asch}. Die gr"o"ste sporadische einfache Gruppe 
ist das Monster, die zweitgr"o"ste das Babymonster.

"Ahnlich wie die Darstellungstheorie von kompakten Lieschen Gruppen bzw.~ihren Liealgebren 
kann die Darstellungs\-theorie von affinen 
Kac-Moody Algebren, den Lie\-algebren von Schleifen\-gruppen kompakter Liescher 
Gruppen, entwickelt werden~\cite{Kac}. Es stellt sich heraus, da"s die
integrablen 
H"ochst\-gewichts\-darstellungen von affinen Kac-Moody Algebren die gr"o"sere
Struktur einer \VOA bzw.~von Moduln hier"uber besitzen~\cite{FreZhu}.
Die Automorphismen\-gruppe dieser \VOAs ist die zu der affinen Kac-Moody Algebra
geh"orige Liesche Gruppe.

Andererseits finden sich auch \VOAs mit endlichen Gruppen als 
Auto\-morphismen\-gruppen, wie z.B.~der
Mondschein\-modul $\VM$, der
das Monster als Auto\-morphismen\-gruppe besitzt~\cite{FLM} oder die
in dieser Arbeit konstruierte \SVOA 
mit dem Baby\-monster als Auto\-morphismen\-gruppe. 
Vom Standpunkt der \VOA-Theorie stehen also endliche und Liesche Gruppen
auf einer Stufe.  

Die f"unf Mathieugruppen lassen sich am einfachsten mit Hilfe der 
Golay Codes verstehen, sieben weitere sporadische Gruppen mit Hilfe des Leechgitters. 
Die restlichen durch das Monster gegebenen sporadischen Gruppen
sollten sich alle in nat"urlicher Weise als Automorphismen\-gruppen von \VOAs beschreiben lassen.
Ob die "ubrigen sechs als "`Parias"' bezeichneten  sporadischen Gruppen
eine angemessene Beschreibung durch \VOAs besitzen ist offen. Das
Studium von \VOAs "uber endlichen K"orpern k"onnte vielleicht Resultate
in dieser Richtung liefern (vgl.~\cite{Ry,BoRy}).

\medskip

{\it Mathematische Physik: Konforme Quantenfeldtheorie und Stringtheorie}

Die der Quantenmechanik zugrunde liegende Mathematik ist im wesentlichen die lineare Algebra von hermiteschen Vektorr"aumen sowie
die Darstellungstheorie von endlichen Gruppen und Liealgebren.

F"ur Quanten\-feld\-theorien (QFTs) wie die Quanten\-elektro\-dynamik 
ist ge\-gen\-w"ar\-tig keine voll\-st"andig be\-frie\-digende 
mathe\-matische Be\-schreibung
ver\-f"ugbar, obwohl dies f"ur einzelne Bau\-steine wie 
die Differential\-geo\-metrie von Prinzipal\-b"undeln zutrifft.
Einen recht all\-gemeinen Rahmen zur Be\-schreibung von QFTs liefern die 
{\it Wightman Axiome\/}~\cite{StWi-PCT,Ha}. Lange Zeit waren allerdings keine 
vollst"andig kon\-struierten nicht\-trivialen Bei\-spie\-le bekannt.

Diese "anderte sich mit dem Auffinden von konformen 
Quanten\-feld\-theorien~\cite{BPZ,Go-mero}. 
Sie sind gerade die QFTs, die den Wightman Axiomen
und den (geeignet er\-weiter\-ten) Atiyah'schen Axiomen einer topo\-logischen 
QFT~\cite{At-knot} f"ur zwei\-dimensionale Raum\-zeiten mit konformer Struktur gen"ugen. \VOAs bilden die einer konformen QFT unter\-liegendende mathematische
Struktur, sie entsprechen der Chiralen Algebra. Die Einbeziehung von Fermionen
als Feldern zu halbganzen Spindarstellungen entspricht der Erweiterung zu 
\SVOAs.

Schlie"slich sei auf Verbindungen zur String\-theorie~\cite{LuTh,Ma-string,Bo-ecm} hingewiesen. 
Nur in Dimension $26$ l"a"st sich die bosonische String\-theorie "`anomaliefrei"' formulieren. Dies ist gerade der Rang der von
Borcherds in~\cite{Bo-lie} aus dem Mond\-scheinmodul kon\-struierten Monster 
Lie\-algebra, die dem Raum der "`physikalischen Zust"ande des Strings"' 
entspricht. In der heterotischen String\-theorie kommt den beiden selbst\-dualen
\VOAs \hbox{$V_{E_8}\otimes V_{E_8}$} und $V_{D_{16}^+}$ vom Rang $16$ 
eine wichtige Rolle zu.

\medskip

{\it Differential-Topologie: Elliptische Geschlechter} 

Der analytische Index ${\rm Ind}(D\otimes E)$ eines mit einem reellen 
Vektorb"undel $E$ getwisteten reellen Diracoperators "uber einer 
$n$-dimensionalen Spin Mannigfaltigkeit $X$ ist ein Element in 
$\widetilde{{\cal M}_n}/i^*\widetilde{{\cal M}}_{n+1}$, den graduierten  Cliffordalgebren $Cl_{n}(\R)$-Moduln modulo denen, die schon
$Cl_{n+1}(\R)$-Moduln sind. Der topologische Index kann durch eine Abbildung
$$\alpha\ :\ \Omega_*^{\rm Spin}({\rm BO})\longrightarrow KO^{-\,*}({\rm pt})$$
in die $KO$-Theorie beschrieben werden.
Der Atiyah-Singer Indexsatz besagt, da"s folgende Gleichheit gilt:
$$ {\rm Ind}(D\otimes E)=\alpha(X,E).$$
Hierbei wird nach Atiyah-Bott-Shapiro  $KO^{-n}({\rm pt})$ mit
$\widetilde{{\cal M}_n}/i^*\widetilde{{\cal M}}_{n+1}$ identifiziert.

Elliptische Geschlechter~\cite{HBJ}
sind die {\it formale\/} string\-theoretisch motivierte
Er\-wei\-terung des topologischen Indexes auf den Schleifen\-raum ${\cal L}X$.
Die \VOAs bilden hier in gewisser Weise das string\-theoretische Analogon zur
Struktur des Vektor\-raumes im klassischen Fall. F"ur eine fixierte \OVOA $V$
betrachtet man die Abbildung
$$ {\cal L}\alpha\ :\ \Omega_*^{\rm Spin}({\rm B} {\rm Aut}(V))\longrightarrow 
KO^{-n}({\rm pt})[[q]],\quad  (X,\widetilde{V})\mapsto 
\alpha(X,\bigotimes_{n=1}^{\infty}S_{q^n}TX_{\bf C}\otimes\widetilde{V}), $$
wobei $\widetilde{V}=\bigoplus_{n\geq 0}\widetilde{V_n}\,q^n$ das zu dem 
${\rm Aut}(V)$-Prinzipalb"undel assoziierte Element in $KO(X)[[q]]$ bezeichnet.
Ist z.B.~$V$ die zur fundamentalen Stufe $1$ Darstellung der affinen Liealgebra 
$\widetilde{{\bf so}}(n)$ 
assoziierte \SVOA $\VF^{\otimes n}$, so erh"alt man f"ur das zum Tangentialb"undel geh"orige 
${\bf SO}(n)$-Prinzipalb"undel das elliptische Geschlecht der Stufe $2$, 
welches auch als die formale $S^1$-"aquivariante Signatur des Schleifenraumes definiert
ist.

Gegenw"artig sind noch kaum Definitionen oder S"atze bekannt, die den
Schleifenraum verwenden oder die volle \VOA-Struktur beinhalten.
Insbesondere fehlt eine zur \hbox{$K$-Theorie} analoge geometrische 
Definition von elliptischer Kohomologie  (vgl.~\cite{Se-ellcoh}).
Ein Resultat in dieser Richtung ist von H.~Tamanoi 
erzielt worden~\cite{Ta-ellvoa}. Er zeigt u.a.~da"s das elliptische Geschlecht
ein \OVOA-Modul "uber der \OVOA der parallelen Schnitte in $\widetilde{V}$ ist.
Das vom Verfasser in seiner Diplomarbeit~\cite{Ho-diplom} f"ur komplexe 
Mannigfaltigkeiten \nopagebreak[4] mit erster Chernklassse $c_1=0$ 
eingef"uhrte und untersuchte universelle komplexe elliptische Geschlecht 
$\varphi_{\rm ell}$ steht in Verbindung zu {$N=2$} supersymmetrischen 
\SVOAs{}~\cite{KaYaYa}.

\bigskip

{\bf Aufbau der Arbeit und Zusammenfassung der Resultate}

In {\it Kapitel 1\/} 
finden sich die Definitionen von Vertex\-operator-Super\-algebren, 
deren Moduln und den Intertwine\-rr"aumen. Es werden die wichtigsten bekannten 
Beispiele und Klassifikations\-resultate zusammen\-gestellt.
Nach einer Einf"uhrung der wichtigsten Begriffe aus der Theorie der 
Modul\-funktionen beschreiben wir die Resultate von Zhu "uber die 
Korrelations\-funktionen auf Riemanschen
Fl"achen vom Geschlecht $1$ und ihre Verallgemeinerungen f"ur 
\SVOAs:

{\bf Satz:} {\sl Die $n$-Punkt Korrelationsfunktionen 
auf dem Torus einer vollst"andigen 
Liste $M_1$, $M_2$, $\ldots$, $M_m$ von irreduziblen Moduln einer 
\netten rationalen \SVOA $V$ formen einen $\GT$-Modul.} 

Die Thetagruppe $\GT$ --- die Untergruppe von $\slz$, die eine der 
Spinstrukturen auf dem Torus $\C/(\Z\tau+\Z)$ fixiert --- operiert dabei in 
der "ublichen Weise auf $\C^n\times\H$.

\medskip
In {\it Kapitel 2\/} werden selbstduale \OVOAs betrachtet.
Au"ser deren Definition und der Diskussion 
der bisher bekannten Beispiele werden 
die folgenden beiden Aussagen "uber den Charakter einer selbstdualen
\OVOA bewiesen:

{\bf Satz:} {\sl (a) Der Charakter einer selbstdualen \netten rationalen \VOA ist ein homogenes Polynom in den Charakteren der \VOAs $V_{E_8}$ und $\VM$.

(b) Der Charakter einer selbstdualen \sehrnetten unit"aren rationalen \SVOA ist
ein homogenes Polynom in den Charakteren der \SVOAs $\VF$ und $V_{E_8}$.}

Insbesondere erhalten wir als Korollar, da"s der Rang einer wie im Satz
betrachteten \SVOA eine halbganze Zahl ist.

\medskip
{\it Kapitel 3\/}
untersucht die Beziehung zwischen selbstdualen \SVOAs und \VOAs. Es
wird der folgende Zusammenhang gefunden:

{\sl Seien $c\in\frac{1}{2}\Z$ und $d>c$, $d\in 8\Z$. Dann besteht eine
$1:1$-Korrespondenz zwischen

 (1) Isomorphieklassen von selbstdualen \sehrnetten unit"aren
\SVOAs vom Rang $c$.

(2) Isomorphieklassen von Paaren $(V,V_{{\bf SO}(k)})$ von 
selbstdualen \netten unit"aren \VOAs $V$ vom Rang $d$ zusammen mit einer
Unter-\VOA $V_{{\bf SO}(k)}$ vom Rang $\frac{k}{2}=d-c$.
(Falls $k=8$, ist wegen der {\it Triality} von ${\bf so}(8)$ zus"atzlich ein 
$V_{{\bf SO}(8)}$-Modul bis auf Isomorphie zu fixieren.)}

Bewiesen werden in dieser Arbeit allerdings nur Teile der Richtung (2) nach (1)
unter einer weiteren schwachen Rationalit"atsforderung. Der Zusammenhang 
erlaubt es uns, die vermutlich vollst"andige Liste aller selbstdualen 
\SVOAs bis zum Rang~$16$ anzugeben.

\medskip
In {\it Kapitel 4\/}
wird dieser Zusammenhang auf den Mondscheinmodul $\VM$ angewendet,
und man erh"alt die Babymonster-\SVOA:

{\bf Satz:} {\sl Es existiert eine \SVOA $\VB$ vom Rang $23\frac{1}{2}$ mit
Charakter
$$\chi_{\VB}= -\frac{31}{16}\,\chi_{\VF}^{47}+\frac{47}{16}\,\chi_{\VF}^{31}
\chi_{E_8}=
q^{-\frac{47}{48}}(1+4371\, q^{\frac{3}{2}}+96256\, q^2
 + 1143745\,{q^{{5\over 2}}} 
  + \cdots \,), $$ 
auf der die Gruppe $2\times B$, wobei $B$ das
Babymonster ist, in nat"urlicher Weise operiert.}

Wir vermuten, da"s $\VB$ selbstdual und die einzige solche \SVOA $V$ vom Rang
$23\frac{1}{2}$ mit $V_{1/2}=V_1=0$ ist. Die \SVOA $\VB$ sollte
als das nat"urliche Objekt angesehen werden, welches das Babymonster $B$,
die zweitgr"o"ste sporadische einfache endliche Gruppe, definiert.

\medskip
In {\it Kapitel 5\/}
schlie"slich werden allgemein {\it extremale selbst\-duale \OVOAs{}\/} 
betrachtet. 
Dies sind selbst\-duale \OVOAs, f"ur die die ersten Koeffizien\-ten des
Charakters so klein wie m"oglich sind, bei Be\-r"uck\-sich\-tigung der 
in Kapitel 2 beschriebenen Struktur der Charaktere. Man erh"alt
in Analogie zu Codes und Gittern das folgende Klassi\-fikations\-resultat:

{\bf Satz:} {\sl Extremale selbstduale \sehrnette unit"are rationale \SVOAs
existieren nur f"ur die R"ange $\frac{1}{2}$, $1$, $\ldots$, $7\frac{1}{2}$, 
$8$, $12$, $14$, $15$, $15\frac{1}{2}$, $23\frac{1}{2}$ und $24$.
F"ur jeden Rang ist genau ein Beispiel bekannt: $\VF$, $\VF^{\otimes 2}$,
$\ldots$, $\VF^{\otimes 15}$, $V_{E_8}$, $V_{D_{12}^+}$, 
$V_{(E_7+E_7)^+}$, $V_{A_{15}^+}$, $V_{E_{8,2}^+}$, $\VB$ und $\VM$.}

F"ur die Beispiele $V_{E_{8,2}^+}$ und  $\VB$ sind allerdings nicht alle 
Eigenschaften bewiesen. Es wird vermutet, da"s die angegebene Liste von
extremalen \SVOAs vollst"andig ist.

\chapter{Grundlagen der Theorie der Vertexoperator-Superalgebren}

In diesem ersten Kapitel finden sich die wichtigsten Begriffe und Resultate
aus der Theorie der \VOAs, die in dieser Arbeit ben"otigen werden.
Eine Einf"uhrung in dieses Gebiet geben die bei\-den B"u\-cher
\cite{FLM} und \cite{FHL}. F"ur eine Zu\-sammen\-fassung des gegen\-w"artigen
Standes der Theorie sei auf die drei "Ubersichtsartikel 
\cite{dong-ueb,li-ueb,huang-ueb} und die dortigen Literaturangaben verwiesen.
\medskip 

Abschnitt 1 enth"alt die Definition von \SVOAs, ihren Moduln und den 
Intertwinerr"aumen. Wir geben die Umformulierung der Jacobi-Identit"at in 
Termen von Korrelations\-funk\-tionen auf der Sph"a\-re, welche uns sp"a\-ter
in Kapitel 3 erlauben wird, die Axiome der dort konstruierten \SVOAs zu 
verifizieren. Zus"atzlich werden symmet\-rische und hermite\-sche 
Bilinear\-formen
betrachtet und die Konstruktion des ("au"seren) Tensorproduktes besprochen.
In Abschnitt 2 beschreiben wir die wichtigsten Beispiele von \VOAs und 
beweisen einige Strukturaussagen "uber \OVOAs, die von Vektoren vom Gewicht
$\frac{1}{2}$, $1$ oder $2$ erzeugt werden.
Abschnitt 3 stellt die in dieser Arbeit ben"otigten Grundlagen aus der Theorie
der Modulfunktionen zusammen. Schlie"slich finden sich in Abschnitt~4 die
Resultate "uber Korrelationsfunktionen auf Riemanschen Fl"achen vom 
Geschlecht~$1$, 
die sich gegenw"artig schon aus dem axiomatischen algebraischen Zugang zur
konformen Quanten\-feld\-theorie mittels \VOAs her\-leiten lassen. Wir formulieren
eine \SVOA-Erwei\-terung des Zhu'chen Resultates~\cite{Zhu-dr} "uber das 
Modultransformationsverhalten der Korrelations\-funktionen der Moduln einer
\VOA.

\section{Definitionen und grundlegende S"at\-ze 
}

Die Axiome von 
Vertexoperator-Algebren sind aus der mathematischen Pr"azisierung des Begriffs 
der chiralen Algebra einer konformen Quantenfeldtheorie auf der
Riemannschen Zahlenkugel entstanden.
\medskip

Die $n$-Punkt Korrelationsfunktionen von Vertexoperatoren sind rationale 
Funktionen, deren Potenzreihenentwicklung im allgemeinen nicht in
ganz $\C^n$ konvergieren werden. Um trotzdem Entwicklungen in verschiedenen
Gebieten vergleichen zu k"onnen, verwenden wir in der Definition von
\OVOAs die formale Entwicklung der $\delta$-Distribution bei $0$:
$$\delta(z)=\sum_{n\in \Z}z^n\in\C[[z,z^{-1}]].$$
Die Reihe $\delta(z_1+z_2)$ wird definiert durch
$$\delta(z_1+z_2)=\sum_{n\in \Z}\sum_{m\in \N}{n \choose m} z_1^{n-m}z_2^m,$$
d.h.~alle binomische Ausdr"ucke werden in nichtnegative ganze Potenzen
in der zweiten Variable $z_2$ entwickelt. Man beachte, da"s die Definition 
in $z_1$ und $z_2$ nicht symmetrisch ist.
\begin{definition}\label{SVOA}
Eine Vertexoperator-Superalgebra (\SVOA{}) ist ein Tupel 
$(V,Y,{\bf 1},\omega)$,
be\-ste\-hend aus ei\-nem $\frac{1}{2}\Z$-gra\-duier\-ten $\C$-Vektor\-raum $V=
\bigoplus_{n\in \frac{1}{2}\Z} V_n$, einer linea\-ren Ab\-bild\-ung
$Y(\,.\,,z):V\longrightarrow {\em End}(V)[[z,z^{-1}]]$,
$a\mapsto Y(a,z)=\sum_{n\in \Z} a_n\, z^{-n-1}$ (dem Vertex\-operator) und 
zwei 
Elementen ${\bf 1}\in V_0$ (dem Vakuum) und $\omega\in V_2$ (dem 
Virasoro\-element), das den fol\-genden Axiomen gen"ugt:

(A1) (Regularit"at)\hfill\newline
--- $\dim V_n<\infty$ f"ur alle $n$ und  $\dim V_n=0$ f"ur $n$
hinreichend klein,\hfill\break
 --- f"ur alle $a$, $b\in V$ ist $a_n b=0$, wenn $n$ gen"ugend 
gro"s,\hfill\break
 --- $Y(a,z)=0$ genau dann, wenn $a=0$.

(A2) $Y({\bf 1},z)={\rm id}_V$.

(A3) Die Koeffizienten von $$Y(\omega,z)= \sum_{n\in \Z} 
{\omega_n}\,z^{-n-1}=\sum_{n\in \Z} {L_n}\,z^{-n-2}$$
erf"ullen die Bedingungen\hfill\break
--- $L_0|_{V_n}=n\cdot {\rm id}_{V_n}$,\hfill\break
--- $Y(L_{-1}a,z)=\frac{d}{dz}Y(a,z)$ f"ur alle $a\in V$ (Translationseigenschaft),
\hfill\break
--- $[L_m,L_n]=(m-n)L_{m+n}+\frac{m^3-m}{12}\, c\, \delta_{m+n,0}\cdot
\rm{id}_V$ (Virasoro Algebra),\hfill\break
 wobei $c\in \C$ als der {\it Rang} von $V$ bezeichnet wird.

(A4) F"ur alle homogenen $a$, $b\in V$ gilt (Jakobiidentit"at):
$$t^{-1}\delta\left(\frac{z-w}{t}\right)Y(a,z)Y(b,w)
-(-1)^{|a||b|}t^{-1}\delta\left(\frac{w-z}{-t}\right)Y(b,w)Y(a,z)$$
$$\qquad\qquad\qquad=w^{-1}\delta\left(\frac{z-t}{w}\right)(Y(Y(a,t)b),w)$$
mit
$$|u| =\cases{0, & f"ur $u\in \bigoplus_{n\in\Z} V_n$ \cr
 1, & f"ur $u\in \bigoplus_{n\in\Z+\frac{1}{2}} V_n$.}$$
\end{definition}
Wir setzen $V_{(0)}=\bigoplus_{n\in\Z} V_n$ bzw.~$V_{(1)}= 
\bigoplus_{n\in\Z+\frac{1}{2}} V_n$ und bezeichnen $V_{(0)}$ als den {\it geraden} und  $V_{(1)}$ als den {\it ungeraden} Teil von $V$. Elemente in
 $V_{(0)}$ (bzw.~$V_{(1)}$) hei"sen {\it gerade} (bzw.~{\it ungerade}).
F"ur das Tupel $(V,Y,{\bf 1},\omega)$ schreiben wir meist nur $V$, falls 
klar ist, welche Vertexoperator Struktur auf $V$ betrachtet wird.
Eine {\it Vertexoperator-Algebra} (kurz \VOA{}) ist eine \SVOA $V=V_{(0)}$,
die nur aus geraden Elementen besteht.

Ein Element $a$ hat das {\it (konforme) Gewicht} $k$, falls $a\in V_k$.
Die Koeffizienten $a_n\in {\rm End}(V)$ des
Vertexoperators $Y(a,z)$ sind dann homogen vom Grad $k-n-1$:
$$ a_n(V_m)\subset V_{m+k-n-1}. $$
Zwei wichtige Folgerungen aus der Jakobiidentit"at sind die 
{\it Kommutatorrelationen: }
\begin{equation}\label{kommutatorrelationen}
[a_n,b_m]=\sum_{i\geq 0}{ m \choose i} (a_ib)_{m+n-i}, 
\end{equation}
und die {\it Assoziativit"atsrelationen: }
\begin{equation}\label{assorelationen}
(a_lb)_n=\sum_{i\geq 0}(-1)^i{l \choose i}\left(a_{l-i}b_{n+i}
- (-1)^l a_{l+n-i}b_i\right)
\end{equation}
Die drei Abbildungen $L_{-1}$, $L_0$ und $L_1$ erzeugen eine Liealgebra
${\bf sl}_2(\C)$, unter denen die Vertexoperatoren zu Elementen
$a\in V$ die folgenden Transformationseigenschaft besitzen:
\begin{eqnarray}
[L_{-1},Y(a,z)] & = & Y(L_{-1}a,z),\label{l-1bed}\\ 
{[}{L_0},Y(a,z)] 
 & = & Y(L_{0}a,z)+z\,Y(L_{-1}a,z),\label{l0bed}\\
{[}{L_{1}},Y(a,z)]  
& = & Y(L_{1}a,z)+2z\,Y(L_{0}a,z)+z^2\,Y(L_{-1}a,z).\label{l1bed}
\end{eqnarray}

Den wichtigen Begriff der Darstellung einer \OVOA beschreibt die
\begin{definition}\label{SVOA-Modul}
Ein Modul einer \SVOA $V$ ist ein Paar $(M,Y_M)$, bestehend aus einem
$\Q$-graduierten $\C$-Vektorraum $M=\bigoplus_{n\in \Q} M_n$
und einer linearen Abbildung \hbox{$Y_M(\,.\,,z):V\longrightarrow {\em End}(M)[[z,z^{-1}]]$},
$a\mapsto Y_M(a,z)=\sum_{n\in \Z} a_n\, z^{-n-1}$, das den folgenden
Axiomen gen"ugt:

(B1) (Regularit"at)\hfill\newline
--- $\dim M_n<\infty$ f"ur alle $n$ und  $\dim M_n=0$ f"ur $n$
hinreichend klein,\hfill\break
 --- f"ur alle $a\in V$, $b\in M$ ist $a_n b=0$, wenn $n$ gen"ugend gro"s,\hfill\break
 --- $Y_M(a,z)=0$ genau dann, wenn $a=0$.

(B2) $Y_M({\bf 1},z)={\rm id}_M$.

(B3) Die Koeffizienten von $$Y_M(\omega,z)= \sum_{n\in \Z} 
{\omega_n}\,z^{-n-1}=\sum_{n\in \Z} {L_n}\,z^{-n-2}$$
erf"ullen die Be\-dingungen\hfill\break
--- $L_0|_{M_n}=n\cdot {\rm id}_{M_n}$,\hfill\break
--- $Y_M(L_{-1}a,z)=\frac{d}{dz}Y_M(a,z)$ f"ur alle $a\in V$  (Translationseigenschaft),
\hfill\break
--- $[L_m,L_n]=(m-n)L_{m+n}+\frac{m^3-m}{12}\, {\rm rank}(V)\, \delta_{m+n,0}\cdot\rm{id}_M$ (Virasoro Algebra).

(B4) F"ur alle homogene $a$, $b\in V$ gilt 
(Jakobiidentit"at):
$$t^{-1}\delta\left(\frac{z-w}{t}\right)Y_M(a,z)Y_M(b,w)
- (-1)^{|a||b|}t^{-1}\delta\left(\frac{w-z}{-t}\right)Y_M(b,w)Y_M(a,z)$$
$$\qquad\qquad\qquad =w^{-1}\delta\left(\frac{z-t}{w}\right)(Y_M(Y(a,t)b),w).$$
\end{definition}
Wir bezeichen den $V$-Modul $(M,Y_M)$ h"aufig mit $M$ und den Vertexoperator
$Y_M$ einfach mit $Y$. Man hat die "ublichen Modulnotationen 
wie Homomorphismus,
direkte Summe und Irreduzibilit"at. Das {\it (konforme) Gewicht eines
irreduziblen Moduls\/} $M$ ist das kleinste $h\in\Q$ mit ${\rm dim}\,M_h>0$.

Offensichtlich ist eine \OVOA stets ein 
Modul "uber sich selbst, den wir als den {\it adjungierten Modul\/} bezeichnen. 
Eine \OVOA hei"st {\it einfach\/}, wenn sie als adjungierter Modul 
irreduzibel ist.

Zu jedem Modul $M$ einer \VOA $V$ gibt es einen dualen Modul $(M',Y_{M'})$
(s.~\cite{FHL}, Abschnitt 5.2):
Der unterliegende Vektorraum ist der ein\-geschr"ank\-te Dual\-raum $M':=
\bigoplus_{n\in \Q} M_n^*$, und der Vertex\-opera\-tor $Y_{M'}$ ist de\-finiert
durch 
$$\langle Y_{M'}(a,z)w',w\rangle=\langle w', Y_M(e^{z L_1}(-z^{-2})^{L_0}
a,z^{-1})w\rangle,$$
wobei $w\in M$, $w'\in M'$ und  $a\in V$.

Wir sagen da"s ein Modul einer \SVOA die Bedingung ${\bf P}_{Vir}$ erf"ullt,
wenn er eine direkte Summe von Virasoroh"ochstgewichtsdarstellungen ist.

In~\cite{Zhu-dr} ist 
von Zhu die folgende Endlichkeitsbedingung ${\bf C}_2$ eingef"uhrt
worden: Sei $V_{-2}V$ der von Elementen des Typs $a_{-2}b$ aufgespannte 
Untervektorraum von $V$. Eine \SVOA $V$ erf"ullt die Bedingung ${\bf C}_2$, 
falls der Quotient $V/(V_{-2}V)$ endlichdimensional ist. 
 
Es ist n"utzlich, zus"atzlich zu den Axiomen in Definition~\ref{SVOA} und
\ref{SVOA-Modul} noch die folgenden drei Bedingungen zu fordern:
\begin{definition}[\nett{}]\label{defnett}
Eine \OVOA $V$ hei"st \nett, wenn die folgenden Eigen\-schaf\-ten 
gelten:\hfill\break
1) $V_m=0$ f"ur $m<0$ und $V_0=\C\cdot{\bf 1}$,\hfill\break
2) alle Moduln von $V$ und $V_{(0)}$ erf"ullen die Bedingung 
${\bf P}_{Vir}$,\hfill\break
3) $V$ erf"ullt die Bedingung ${\bf C}_2$ von Zhu.
\end{definition}
Die Bedingungen 1) und 2) sind auch aus physikalischer Sicht sinnvoll. 
Vie\-le S"at\-ze gelten nur unter diesen Voraussetzungen
(vgl.~z.B.~\cite{Lian}, Kap.~3),
ansonsten lassen sich Gegen\-beispiele finden. 
Ve\-rschiedene Autoren nehmen daher unter\-schied\-liche Tei\-le 
von~\ref{defnett}
in die De\-fini\-tion einer \OVOA mit auf.
Der fol\-gende Satz be\-n"o\-tigt z.B.~die Eigen\-schaf\-ten 1) und 2)
(s.~\cite{Li-bili}).
\begin{satz}[Bilinearform]\label{bilinearform}
Jede \nette \VOA $V$ besitzt eine ein\-deu\-tige nicht aus\-gearte\-te
symmet\-rische in\-varian\-te Bi\-linear\-form $(\,.\,,\,.\,)$, die durch
$$(a,b)\cdot {\bf 1}={\rm Res}_z z^{-1}(Y(e^{zL_1}(-z^{-2})^{L_0}a,z^{-1})b) 
\quad\hbox{f"ur }a,b\in V,$$
gegeben ist.
\end{satz}
Invarianz ist hierbei durch 
$$(Y(u,z)v,w)=(v,Y(e^{zL_1}(-z^{-2})^{L_0}u,z^{-1})w)
\quad\hbox{f"ur alle $u$, $v$, $w\in V$}$$
definiert. 

\medskip
Eine Beschreibung von Vertexoperatoren zwischen drei Moduln 
gibt die fol\-gen\-de 
De\-fini\-tion, in der wir uns auf den sp"ater be\-n"otig\-ten Fall von \VOAs
beschr"anken. Dazu sei f"ur einen Vektorraum $V$ mit $V\{z\}=\{\sum_{n\in\Q}v_nz^n\mid v_n\in V \}$ 
der Raum der $V$-wertigen formalen Reihen in rationalen $z$-Potenzen bezeichnet.
\begin{definition}[Intertwiner Operatoren]\label{SVOA-Intertwiner}
Sei $V$ eine \VOA und seien $(M_i,Y_i)$, $(M_j,Y_j)$ und $(M_k,Y_k)$
drei $V$-Moduln. Der Raum ${M_k \choose M_i\ M_j}_V $ der Intertwineroperatoren
vom Typ ${k \choose i\ j}$ ist definiert als der $\C$-Vektorraum von linearen 
Abbildungen ${\cal Y}(\,.\,,z):W_i\longrightarrow {\rm Hom}(W_j,W_k)\{z\}$,
${\cal Y}(w,,z)=\sum_{n\in \Q}w_n\, z^{-n-1}$, die den folgenden Axiomen 
gen"ugen:

(C1) (Regularit"at)\hfill\break
 --- f"ur alle $m^{(i)}\in M_i$, $m^{(j)}\in M_j$ ist ${m^{(i)}}_n m^{(j)}=0$, 
wenn $n$ gen"ugend gro"s,\hfill\break
 --- ${\cal Y}(m^{(i)},z)=0$ genau dann, wenn $m^{(i)}=0$.\hfill\break

(C2) (Translationseigenschaft)\hfill\break
$${\cal Y}(L_{-1}m^{(i)},z)=\frac{d}{dz}{\cal Y}(m^{(i)},z)\quad\hbox{f"ur alle 
$m^{(i)}\in M_i$}.$$
(C3) F"ur alle $a\in V$, $m^{(i)}\in M_i$ gilt (Jakobiidentit"at):
$$t^{-1}\delta\left(\frac{z-w}{t}\right)Y_k(a,z){\cal Y}(m^{(i)},w)
-t^{-1}\delta\left(\frac{w-z}{-t}\right)
{\cal Y}(m^{(i)},w)Y_j(a,z)$$
$$\qquad\qquad\qquad =w^{-1}\delta\left(\frac{z-t}{w}\right)
({\cal Y}(Y_i(a,t)m^{(i)}),w).$$
\end{definition}
Wir setzen $N_{ij}^k=\dim {M_k \choose M_i\ M_j}_V$ und bezeichnen 
diese Zahlen als die {\it Fusionsregeln\/}.
Sie beschreiben 
in gewisser Weise die 
Zerlegung des (inneren) Tensorproduktes zweier irreduzibler Moduln in 
irreduzible Komponenten.
Die Tensorprodukttheorie war von Huang und Lepowsky in 
den Arbeiten~\cite{HuLe-tensor,HuLe-te1,HuLe-te2,HuLe-te3,HuLe-te4} 
eingef"uhrt worden.
Ein anderer Ansatz hierf"ur findet sich in~\cite{Li-dr}.
\begin{definition}[Fusionsalgebra]\label{deffusion}
Die Fusions\-algebra ${\cal F}(V)$ einer \SVOA $V$ ist der von den 
ir\-re\-du\-zib\-len Moduln von $V$ er\-zeug\-te freie $\Z$-Modul 
zu\-sammen mit der auf den Er\-zeugern 
durch $M_i\times M_j = \sum_{k} N_{ij}^k\,M_k$ de\-finier\-ten
bi\-linearen Ab\-bildung
\hbox{$\times:{\cal F}(V)\otimes {\cal F}(V)\longrightarrow {\cal F}(V)$.} 
\end{definition}
H"aufig wird ${\cal F}(V)$ auch als $\C$-Vektorraum angesehen. Es wird vermutet, da"s ${\cal F}(V)$ eine kommutative und assoziative Algebra ist. Dies ist der 
Fall, wenn die Tensorprodukttheorie von Huang und Lepowsky anwendbar ist.

\medskip
Die Jacobiidentit"at f"ur \SVOAs und ihre Moduln l"a"st sich mit Hilfe der "`Dualit"at"' umformulieren.
Sei $V$ eine \SVOA und $M$ ein Modul. F"ur Elemente $a_1$, $\dots$, $a_n\in V$
und $v\in M$, $v'\in M'$ konvergiert
in dem Gebiet $|z_1|>\cdots>|z_n|>0$ die Reihe
$$\langle v',Y(a_1,z_1)\ldots Y(a_n,z_n)v\rangle\in \C[[z_1,z_1^{-1},\ldots,
z_n,z_n^{-1}]]$$ 
zu einer rationalen Funktion
$R_{v',v}((a_1,z_1),\ldots,(a_n,z_n))\in \C[z_i,z_i^{-1},\frac{1}{z_i-z_j}]$,
die m"ogliche Pole bei $0$, $\infty$ und $z_j$ als Funktion in $z_i$ besitzt.
Wir bezeichnen diese Funktion als die {\it $n$-Punkt Korrelations\-funk\-tion
auf der Sph"are.}
\begin{satz}[Kommutativit"at]
Seien $a_1$, $\ldots$, $a_n\in V$ homogene Elemente und $\sigma\in S_n$ eine
Permutation. Dann gilt f"ur die $n$-Punktkorellationsfunktionen 
(Kommutativit"at):
\begin{equation}\label{kommutativitaet}
R_{v',v}((a_1,z_1),\ldots,(a_n,z_n))=(-1)^w
R_{v',v}((a_{\sigma(1)},z_{\sigma(1)}),\ldots,(a_{\sigma(n)},z_{\sigma(n)})),
\end{equation}
wobei $w$ die Anzahl der Vertauschungen von ungeraden Elementen 
$a_i\in V_{(1)}$ in $\sigma$ ist.
\end{satz}\nopagebreak[2]
Physikalisch entspricht die (Anti-)Kommutativit"at der Kausalit"at --- eine
der funda\-menta\-len Forderungen an eine Quantenfeldtheorie.
\pagebreak[2]

In~\cite{DoLe}, Kapitel 7, wurde gezeigt:
\begin{satz}\label{Kommutativitaet}
In Definition~\ref{SVOA} bzw.~\ref{SVOA-Modul} kann die Jacobi\-iden\-ti\-t"at 
(Axiom (A4) bzw. (B4)) durch die Ra\-tio\-nali\-t"at und Kommu\-ta\-tivit"at
(\ref{kommutativitaet}) der $2$-Punkt Korre\-lations\-funktionen auf der 
Sph"a\-re und den beiden Be\-dingungen (\ref{l-1bed}) und (\ref{l0bed}) f"ur 
$L_{-1}$ und $L_0$ ersetzt werden.
\end{satz}
Die Rationalit"at und Kommutativit"at ist in den Anwendungen oft leichter
nachzupr"ufen als die Jacobiidentit"at.

Eine wichtige Klasse von \OVOAs besitzt zus"atzlich die folgenden Eigenschaften:
\begin{definition}[rational]\label{rational}
Eine \OVOA hei"st rational, falls der Rang rational ist, sie nur endlich viele
irreduzible Moduln besitzt und jeder endlich erzeugte Modul eine direkte Summe
von irreduziblen Moduln ist. Zus"atzlich sollen f"ur \SVOAs diese Eigenschaften
auch f"ur die gerade Unter-\VOA gelten.
\end{definition}
Es ist zu vermuten, da"s die Rationalit"at des Ranges aus den anderen 
Bedingungen folgt (vgl.~\cite{AnMo}).

\medskip

Wir sagen: Eine \VOA $V$ ist von {\it reellem Typ}, falls sie die 
Komplexifizierung einer "uber $\R$ definierten \VOA $V_{\R}$ ist:
$V=V_{\R}\otimes_{\R}\C$.
\begin{definition}[unit"ar]\label{defunitaer}
Eine \VOA $V$ hei"st unit"ar, wenn sie von reellem Typ ist
und jeder $V$-Modul $M$ eine positiv definite
in\-varian\-te hermi\-tesche Bilinear\-form $\langle\,.\,,\,.\,\rangle_M$
besitzt. Invarianz bedeutet, da"s f"ur Elemente $u$, $v \in M$ und einen 
homogenen ${\tx sl}_2(\C)={\rm Span}_{\C}(L_{-1},L_0,L_1)$ 
H"ochstgewichtsvektor
$a\in (V_{\R})_h$ f"ur alle $n\in\Z$ die Beziehung 
$\langle a_{-n-2(h-1)}u,v\rangle_M=\langle u,a_nv\rangle_M$
gilt. 
\end{definition}
Eine Unter-\VOA einer unit"aren \VOA ist auch wieder unit"ar.
Eine \SVOA hei"st unit"ar, wenn ihre gerade Unter-\VOA unit"ar ist.
\begin{lemma}\label{lemmaunitaer}
F"ur den Rang $c$ einer \netten unit"aren \VOA $V$
und das konforme Gewicht $h$ eines $V$-Moduls $M$ gelten $c\geq 0$, $h\geq 0$. 
\end{lemma}
{\bf Beweis: }
Sei $v\in M$ ein normierter Virasoroh"ochstgewichtsvektor vom Gewicht $(c,h)$.
Aus der Beziehung $L_nL_{-n}v=L_{-n}L_nv+2nh\,v+\frac{n(n^2-1)}{12}c\,v$
folgt $\langle L_{-n} v, L_{-n} v \rangle_M=2nh+\frac{n(n^2-1)}{12}c$ f"ur 
$n > 0$.
Daraus ergibt sich f"ur $n=1$ die Ungleichung $h\geq 0$, und f"ur
$n$ sehr gro"s die Ungleichung $c\geq 0$. \qed
\begin{lemma}\label{lemmareduktiv}
Der Gewicht $1$ Anteil $V_{1}$ einer \netten  unit"aren \SVOA $V$ ist
eine reduktive Liealgebra, d.h.~die direkte Summe einer abelschen und einer
halbeinfachen Liealgebra.
\end{lemma}
{\bf Beweis:} 1) ${\tx g}:=V_{1}$ ist eine Lie Algebra mit Lieklammer $[a,b]=a_0b$.
\hfill\newline
2) Ist ${\tx u}$ ein Ideal von ${\tx g}$, so ist es wegen der Invarianz
der hermiteschen Form auch das orthogonale Komplement ${\tx u^{\bot}}$.
\hfill\newline
3) Da die Form positiv definit ist, gilt Argument 2) 
auch f"ur alle Lie Unteralgebren von ${\tx g}$.\hfill\newline
Es folgt, da"s ${\tx g}$ eine direkte Summe von einfachen bzw.~abelschen
Summanden ist.\qed

\medskip

Seien $(V_1,Y_1,{\bf 1}_1,\omega_1)$, $\ldots$, 
$(V_n,Y_n,{\bf 1}_n,\omega_n)$ \OVOAs.
Das Tensorprodukt
$$ V=V_1\otimes \cdots \otimes V_n$$
kann in nat"urlicher Weise mit einer \OVOA-Struktur $(V,Y,{\bf 1},\omega)$
versehen werden:
\begin{eqnarray}
Y(v_1\otimes\cdots\otimes v_n,z) & = &  Y_1(v_1,z)\otimes\cdots\otimes
 Y_n(v_n,z) \quad (v_i\in V_i),\label{deftensorprodukt}\\
{\bf 1} & = & {\bf 1}_1\otimes\cdots\otimes{\bf 1}_n,\nonumber\\
\omega & = & \omega_1\otimes {\bf 1}_2 \otimes\cdots\otimes {\bf 1}_n +\cdots +
{\bf 1}_1\otimes \cdots\otimes {\bf 1}_{n-1}\otimes\omega_n .\nonumber
\end{eqnarray}
Das Tensorprodukt in (\ref{deftensorprodukt}) ist dabei in einem
$\Z_2$-graduierten Sinne zu verstehen.
Analog definiert man das Tensorprodukt f"ur Moduln und Intertwineroperatoren.
Wir fassen die wichtigsten Eigenschaften des Tensorproduktes zusammen in
\begin{satz}\label{tensorprodukt}
Das Tupel $(V,Y,{\bf 1},\omega)$ ist eine \OVOA, deren Rang die Summe
der R"ange der ein\-zelnen Fak\-toren ist. Ihre irre\-duzib\-len Moduln 
(wobei alle Eigen\-werte der Vira\-soro\-opera\-toren $(L_i)_0$ als rational
vor\-aus\-gesetzt seien)
sind
gerade die Tensorprodukte der irreduziblen
Moduln der einzelnen $V_i$, $i=1$, $\ldots$, $n$.
Besitzen alle Faktoren die Eigenschaften \nett, unit"ar oder rational, so
besitzt diese auch das Tensorprodukt.
\end{satz}
Zum Beweis der \OVOA-Eigenschaft des Tensorproduktes und der Aussage "uber die
irreduziblen Moduln siehe~\cite{FHL} und~\cite{DoLe}, Kap.~10.
Die Rationalit"at wurde in~\cite{DoMaZhu} bewiesen, die Eigenschaften
\nett und unit"ar folgen leicht aus der Definition.

Das ("au"sere) Tensorprodukt "`$\otimes$"' der \OVOAs $V_1$, $\ldots$, $V_n$
ist von dem (inneren) Tensorpodukt 
von Moduln $M_1$, $\ldots$, $M_n$ einer festen \OVOA zu unterscheiden,
das in Anschlu"s an Definition~\ref{SVOA-Intertwiner} erw"ahnt wurde.

\section{Beispiele von \OVOAs und Struktur\-s"atze}

Die Angabe einer Basis ist die einfachste M"oglichkeit, einen Code oder
ein Gitter zu definieren, und umgekehrt kann jeder Code oder jedes Gitter auf
diese Weise beschrieben werden. 
Die gleiche Rolle "ubernimmt die Cartanmatrix f"ur eine
endlichdimensionale halbeinfache komplexe Liealgebra. Wesentlich schwieriger
ist die Situation bei \OVOAs{}: Es ist keine einfache explizite Beschreibung
einer \OVOA bekannt, aus der sich die G"ultigkeit aller Axiome unmittelbar 
erg"abe (siehe aber~\cite{Li-local}). 
Umgekehrt ist auch nicht bekannt, ob {\it endlich viele\/}
komplexe Zahlen gen"ugen, um alle \SVOAs von z.B.~festem Rang in einfacher
nat"urlicher Weise zu beschreiben.~\footnote{Mengentheoretisch ist dies
trivial, aber dies ist hier nicht das Problem.} Alle bisher bekannten Beispiele
von \netten rationalen \VOAs und ihren Moduln gehen auf vier Beispielklassen
zur"uck, in denen die \VOA-Struktur aus einer einfacheren bekannten Struktur 
konstruiert wird, und dann mit deren Hilfe die Axiome verifiziert werden.

Im einzelnen sind dies die \SVOAs, die zu den H"ochst\-gewichts\-darstellungen
einer unendlich
dimensionalen Clifford\-algebra, einer affinen Kac-Moody Algebra
oder der Virasoro\-algebra assoziert sind, 
sowie die \SVOAs, die mit Hilfe ganzer
Gitter konstruiert werden. Die ersten drei Beispielklassen
lassen sich dadurch charakterisieren, da"s sie von den Vertexoperatoren
zu Vektoren vom Gewicht $\frac{1}{2}$, $1$ oder $2$ erzeugt werden. Dies
wird f"ur die Klassifikations\-"uberlegungen in dieser Arbeit wichtig sein.
Alle weiteren Beispiele sind hieraus mittels des 
Tensorproduktes~\cite{FHL,DoMaZhu},
der Kommutantenkonstruktion~\cite{FreZhu} oder mittels 
Orbifoldkonstruktionen~\cite{FLM,DoMa-moonshine} erhalten worden.

\medskip
Nach Definition ist jede \nette \VOA eine direkte Summe von Moduln "uber
der vom Virasoroelement erzeugten Virasoroalgebra. Eine besonders einfache 
Struktur haben daher die \VOAs, die nur aus einem solchen Modul bestehen,
n"amlich dem vom Vakuum erzeugten. Sei $L_{c_{p,q}}(0)$ die irreduzible 
H"ochstgewichtsdarstellung der Virasoralgebra vom H"ochstgewicht $(c_{p,q},0)$, 
mit $c_{p,q}=1-6(p-q)^2/pq$ und  $p$, $q\in \{2,3,4,\ldots\}$, $(p,q)=1$;
seien $L_{c_{p,q}}(h_{m,n})$ die irreduziblen Darstellungen vom
H"ochstgewicht $(c_{p,q},h_{n,m})$ mit $h_{n,m}=\frac{(np-mq)^2-(p-q)^2}{4pq}$,
$0<m<p$, $0<n<q$. Der Vektorraum $L_{c_{p,q}}(0)$ besitzt eine nat"urliche
\VOA-Struktur (s.~\cite{FreZhu}, Abschnitt 4) und 
die $L_{c_{p,q}}(h_{m,n})$ sind \VOA-Moduln hier"uber.
Der folgende Satz von Wang (s.~\cite{Wan}, Satz 4.2) sagt aus, da"s die in der 
Physik als 
"`Minimale Modelle"' bezeichneten \VOAs $L_{c_{p,q}}(0)$, rational sind.
\begin{satz}[Virasoro \VOAs{}]\label{satzminimalemodelle}
Die \VOA $L_{c_{p,q}}(0)$ ist rational, und die $L_{c_{p,q}}(h_{m,n})$,
$0<m<p$, $0<n<q$ sind gerade alle irreduziblen Darstellungen
von $L_{c_{p,q}}(0)$.
\end{satz}
Die \VOAs $L_{c_{p,q}}(0)$ sind f"ur $q=p+1$, d.h.~$c_{p,q}=1-\frac{6}{m(m+1)}$,
$m\in\{3,4,\dots\}$,  unit"ar.

Unter gewissen Bedingungen ist der von den Koeffizienten des Vertex\-operators
eines Elementes vom Gewicht $2$ vom Vakuum erzeugte Unter\-vektor\-raum eine
zu einer Vira\-soro\-algebra asso\-ziierte Un\-ter-\VOA:
\begin{satz}[Eigenschaft \Lvir{}]\label{L2-vir}
Sei $V$ eine \nette $\VOA$ und $a\in V_2$ ein Element vom Gewicht $2$, so
da"s die Koef\-fizien\-ten $a_n$ des Vertex\-operators 
$Y(a,z)=\sum_{n\in\Z} a_n\, z^{-n-1}$
die folgenden drei Bedingungen erf"ullen:
\begin{equation}\label{Vir-bedingung}
a_1a=2\cdot a,\qquad a_2a=0 \quad\hbox{und}\qquad  a_3a=\frac{d}{2}.
\end{equation}
Dann gilt:
\begin{list}{}{}
\item
1) Das System $\{a_n\mid n\in \Z\}\cup\{{\rm id}_V\}$ erzeugt eine 
Virasoroalgebra ${\rm Vir}_a$ unter der Lieklammer.
\item
2) Der von ${\rm Vir}_a$ erzeugte Untervektorraum $W_a:={\cal U}({\rm Vir}_a)
{\bf 1}={\cal U}({\rm Vir}_a^-){\bf 1}$ 
ist eine Unter-VOA von $V$ vom Rang $d$ mit Virasoroelement $a$.
\item
3) Die Unter-\VOA $W_a$ ist als \VOA isomorph zu einem Quotienten von
$M_{d}:={\cal U}({\rm Vir}^-)/\langle L_{-1}{\bf 1} \rangle$.
\end{list}
\end{satz}
Wir bezeichnen die Bedingungen (\ref{Vir-bedingung}) an ein Element $a$
als {\it Eigenschaft \Lvir{}}.

{\bf Beweis:}\newline
Zu 1) Aus den Kommutatorrelationen~(\ref{kommutatorrelationen})
und den Beziehungen (\ref{Vir-bedingung}),
sowie $a(n)a=0$ f"ur $n\geq 3$, erh"alt man bei Beachtung der Indexverschiebung
$x(n):=x_{n+1}$ die Gleichungen
\begin{eqnarray*}
[a(n),a(m)] & = &\sum_{i \geq -1}{n+1 \choose i+1}(a(i)a)(m+n-i) \\
&= &(a(-1)a)(m+n+1) + 2 (n+1) \cdot a(m+n) + {n+1 \choose 3}\frac{d}{2}
\cdot {\bf 1}(m+n-2).
\end{eqnarray*}
Dies vereinfacht sich zu 
$$[a(n),a(m)]=\frac{1}{2}\left([a(n),a(m)]-[a(m),a(n)]\right)
 =(n-m)\cdot a(m+n) + \frac{n^3-n}{12}\delta_{m+n,0}\cdot {\rm id}_V, $$
also gerade zu den Kommutatorrelationen der Virasoroalgebra.

Zu 2) Ein Element $w\in W_a$ ist wegen der Virasororelation und $a_0{\bf 1}=0$
eine Linearkombination von Elementen des Typs 
$v=a_{-i_1}a_{-i_2}\dots a_{-i_k}{\bf 1}$, wobei $i_{\nu}\geq 0$, $\nu=1$,
$\ldots$, $k$. Wir zeigen durch vollst"andige Induktion "uber $k$, da"s der
Koeffizient $v_n$ des Vertexoperators $Y(v,z)$ eine Linarkombination
von Abbildungen $a_{-j_1}a_{-j_2}\dots a_{-j_l}$, $j_{\mu}\in \Z$ ist.
F"ur $k=1$ gilt $(a_{-i_1}{\bf 1})_n={-i_1 \choose -n+1} a_{-i_1+n+1}$;
allgemein gilt wegen der Assoziativit"atsrelation~(\ref{assorelationen})
\begin{eqnarray}\label{rekursion}
v_n & = & (a_{-i_1}a_{-i_2}\dots a_{-i_k}{\bf 1})_n  \\
 & = & \sum_{\mu\geq 0}(-1)^{\mu}{- i_1 \choose \mu}\left(
a_{-i_1-\mu}(a_{-i_2}\dots a_{-i_k}{\bf 1})_{n+\mu}
- (-1)^{-i_1}(a_{-i_2}\dots a_{-i_k}{\bf 1})_{-i_1+n-\mu}a_{\mu} \right),
\nonumber
\end{eqnarray}
d.h.~aufgrund der Induktionsannahme hat $v_n$ die verlangte Gestalt.
Wir haben damit gezeigt, da"s $W_a$ in $V$ abgeschlossen ist, d.h.~eine 
Unter-\VOA bildet.

Zu 3)\, Gleichung (\ref{rekursion}) zeigt auch, da"s die \VOA-Struktur von
$W_a$ schon vollst"andig durch die Abbildungen $a_n$ bestimmt ist.
Andererseits is $W_a$ nach 1) und 2) ein H"ochst\-gewichts\-modul "uber
der Virasoroalgebra vom H"ochstgewicht $(d,0)$ mit $a_0 {\bf 1}=0$. Er ist
daher ein \VOA-Quotient des Vermamodulquotienten ${\cal U}({\rm Vir}^-){\bf 1}/
\langle a_0 {\bf 1} \rangle=M_d$, der nach \cite{FreZhu}, Abschnitt 4, eine
\VOA-Struktur besitzt. \qed

\medskip
Eine wichtige Klasse von \SVOAs wird von unendlich\-dimensionalen 
Cliffordalgebren erzeugt.

Sei $A$ ein $l$-dimensionaler Vektorraum zusammen mit einer nichtentarteten
Bilinearform $\langle .,. \rangle$. Wir setzen
$$ A(\Z+\frac{1}{2})=\bigoplus_{n \in \Z +\frac{1}{2}} A(n),$$
wobei $A(n)$ eine Kopie von $A$ vom Gewicht $n\in \Z+\frac{1}{2}$ ist.
Auf $ A(\Z+\frac{1}{2})$ sei f"ur 
Elemente $a(n)\in A(n)$, $b(m)\in A(m)$ und $m$, $n\in\Z+\frac{1}{2}$
die Form $\langle .,. \rangle$ durch
$\langle a(n),b(m)\rangle=\langle a,b \rangle\delta_{n,-m}$ fortgesetzt.
Die Untervektorr"aume $A^{\pm}(\Z+\frac{1}{2})=\bigoplus_{0<n\in \Z+\frac{1}{2}}
A(\pm n)$ sind dann maximale isotrope Teilr"aume.
Bez"uglich dieser Bilinearform erzeugt $A(\Z+\frac{1}{2})$ eine 
unendlich dimensionale Cliffordalgebra ${\bf Cliff}(\Z+\frac{1}{2})$.
Sei ${\cal J}$ das von $A^+(\Z+\frac{1}{2})$ erzeugte Ideal
in ${\bf Cliff}(\Z+\frac{1}{2})$. Der Vektorraum
$$V_{{\rm Fermi},l}:={\bf Cliff}(\Z+\frac{1}{2})/ {\cal J} =
\Lambda^*(A^-(\Z+\frac{1}{2}))\,{\bf 1}=\bigotimes_{0<n\in\Z+\frac{1}{2}}
\Lambda^*(A(-n))\, {\bf 1}$$
ist dann ein irreduzibler ${\bf Cliff}(\Z+\frac{1}{2})$-Modul. Tats"achlich
besitzt $V$ sehr viel mehr Struktur. Wir fassen dies zusammen in
\begin{satz}[Clifford \SVOAs{}]\label{fermion}
Der Vektorraum $V_{{\rm Fermi},l}$ besitzt die Struktur einer 
\netten unit"aren \SVOA vom Rang
$\frac{l}{2}$. Die Koeffizienten $a_n$ des zu $a\in V_{1/2}\cong A$ geh"origen
Vertexoperators $Y(a,z)$ operieren auf $V$ gerade wie $a(n)$ durch die
Cliffordmultiplikation. Die \SVOA $V_{{\rm Fermi},l}$ ist rational und 
besitzt genau eine irreduzible Darstellung n"amlich $V$ selbst.
\end{satz}
F"ur diese Konstruktion vgl.~\cite{tsukada,ffr} und \cite{KacWang}, Abschnitt 4.
In den ersten beiden Referenzen wird $l$ als gerade vorausgesetzt, 
dies ist aber eine nicht notwendige Einschr"ankung. In der Physik werden diese
\SVOAs als "`freie Fermionen"' bezeichnet.

Die Koeffizienten der Vertexoperatoren zu Elementen aus $V_{1/2}$ einer
\SVOA $V$ bilden ganz allgemein eine unendlich dimensionale Cliffordalgebra.
F"ur die von ihnen erzeugte Unter-\SVOA gilt (vgl.~Prop.~6.2~\cite{Ta},
der Beweis dort bleibt auch f"ur ungerade $l$ g"ultig):
\begin{satz}[Eigenschaft \Lcliff{}\/]\label{Leinhalb-cliff}
Sei $V$ eine \nette unit"are \SVOA mit ${\rm dim}\, V_{1/2}=l$.\linebreak[2]
Dann ist die von $V_{1/2}$ erzeugte Unter-\SVOA isomorph zu
$V_{{\rm Fermi},l}$.
\end{satz}
Die \SVOA $V_{{\rm Fermi}}:=V_{{\rm Fermi},1}$ hat den Rang $\frac{1}{2}$
und die gerade Unter-\VOA ist $L_{1/2}(0)$. Als $L_{1/2}(0)$-Modul besteht die
Zerlegung $\VF=L_{1/2}(0)\oplus L_{1/2}(\frac{1}{2})$.

\medskip
Eine dritte Klasse sind die zu affinen Kac-Moody Algebren
assozierten und in der Physik als "`WZW-Modelle"' bekannten \VOAs.
Sei ${\tx g}$ eine einfache komplexe Liealgebra und 
$$\hat{{\tx g}}={\tx g}\otimes \C[t,t^{-1}]\oplus \C\cdot c=
\hat{{\tx g}}^+\oplus {\tx g}\oplus\hat{{\tx g}}^-\oplus \C\cdot c$$
die zugeh"orige affine Liealgebra, sei $\{\Lambda_0,\ldots,\Lambda_n\}$ ein
System von funda\-mentalen Ge\-wichten von $\hat{{\tx g}}$, so da"s
$\{\Lambda_1,\ldots,\Lambda_n\}$ fundamentale Gewichte f"ur ${\tx g}$ sind.
F"ur eine H"ochstgewichtsdarstellung $V_{\lambda}$ von ${\tx g}$ mit 
H"ochstgewicht $\lambda\in \Z_{\geq 0}\,\Lambda_1\oplus\ldots\oplus
\Z_{\geq 0}\,\Lambda_n$ und ein $k\in \Z_{\geq 0}$ bezeichnen wir mit
$\hat{V}_{k\Lambda_0+\lambda}={\cal U}(\hat{\tx g})\otimes_{
{\cal U}({\tx g}\oplus\hat{{\tx g}}^-\oplus \C\cdot c})V_{\lambda}$ den 
Vermamodul vom H"ochstgewicht $k\Lambda_0+\lambda$. Er enth"alt einen 
maximalen eindeutigen Untermodul $I(k,\lambda)$. Sei $L_{k\Lambda_0+\lambda}$ 
der Quotientenmodul $\hat{V}_{k\Lambda_0+\lambda}/I(k,\lambda)$. Bezeichne
schlie"slich mit $\check{g}$ die duale Coxeterzahl, mit $\theta$ die l"angste
Wurzel von ${\tx g}$ und schreibe $(\,.\, , \, .\,)$ f"ur die invariante
Bilinearform auf~${\tx g}$, die so normiert ist, da"s $2$ die Quadratl"ange 
einer langen Wurzel ist. 
\begin{satz}[Kac-Moody \VOAs{}]\label{satzlievoas}
Der $\hat{{\tx g}}$-Modul $L_{k\Lambda_0}$ besitzt eine nat"urliche 
\VOA-Struk\-tur vom Rang $\frac{k\cdot\dim {\tx g}}{k+\check{g}}$.
Die so er\-haltene VOA ist \nett, uni\-t"ar und ra\-tional und besitzt als
irre\-duzible Moduln gerade die irre\-duziblen inte\-grablen 
H"ochst\-gewichts\-darstel\-lungen $\{L_{k\Lambda_0+\lambda}\mid 
(\lambda,\theta)
\leq k\}$ von $\hat{{\tx g}}$ der Stufe $k$.
\end{satz}
\medskip

Zu jedem ganzen positiv definiten Gitter $L\subset \R^n$ kann eine 
\SVOA $V_L$ asso\-ziiert werden. Dem geraden Untergitter $L_0\subset L$
entspricht dabei die gerade Unter-\VOA $V_{(0)}=V_{L_0}$, d.h.~$V_L$ ist
eine \VOA, falls $L$ gerade. Bezeichne mit 
$L^*=\{x\in \R^n\mid \langle x,y\rangle \in \Z \ \hbox{f"ur alle\ }\hbox{$y\in 
L$}\}$ das zu $L$ duale Gitter.
\begin{satz}[Gitter-\SVOAs{}]\label{satzgittersvoa}
Die Gitter-\SVOA $V_L$ ist eine \nette unit"are rationale
\SVOA vom Rang $n$, deren irreduzible Moduln den Nebenklassen 
$L^*/L$ entsprechen.
Die Fusionsalgebra ist isomorph zum Gruppenring $\Z[L^*/L]$. 
\end{satz}
Die \VOA-Struktur war in \cite{FLM} (gerade Gitter) und \cite{DoLe}
(ungerade Gitter) konstruiert worden. Die irreduziblen Moduln waren im
geraden Fall in~\cite{Do-gitter} bestimmt worden. Der etwas allgemeinere Fall
von beliebigen ganzen Gittern und deren assozierten \SVOAs
ist zumindest in der physikalischen Literatur ebenfalls untersucht worden
(vgl.~z.B.~\cite{Go-mero}). Dort werden diese \SVOAs durch die Quantisierung
eines "`Strings auf dem Torus $\R^n/L$"' konstruiert.

Es besteht der folgende Zusammenhang zwischen Gitter-\VOAs und den zu affinen
Kac-Moody Algebren geh"origen \VOAs (s.~\cite{FLM}, Kapitel 7):
\begin{satz}\label{satzgitterlieiso}
Die zu den Wurzelgittern der Liealgebren vom Typ $A_n$, $D_n$ ($n\geq 1$) und
$E_n$ ($n=6$, $7$, $8$) geh"origen \VOAs sind isomorph zu den \VOAs, die zu den 
fundamentalen Stufe $1$ Darstellungen der entsprechenden affinen Kac-Moody 
Algebren assoziiert sind.
\end{satz}
Die zum geraden Gitter $L=\sqrt{2k}\,\Z$; $k\in \N$, $k>1$ assoziierte \VOA
kann als die zu der abelschen Liealgebra ${\rm Lie}({\bf U}(1))=:{\tx t}$ 
assozierte \VOA bei Stufe $k$ aufgefa"st werden.

Die Koeffizienten von Vertexoperatoren zu Elementen vom
Gewicht $1$ bilden eine affine Liealgebra. Zusammen mit 
Lemma~\ref{lemmareduktiv} ist daher zu vermuten, da"s {\it unit"are\/} \nette
rationale \VOAs die folgende Eigenschaft besitzen: 
\begin{definition}[Eigenschaft \Llie{}]\label{L1-lie}
Sei $V$ eine \VOA, $\widetilde{V}_1$ die von $V_1$ erzeugte Unter-\VOA und
$W={\rm Com}_V(\widetilde{V}_1)$ die Kommutante von $\widetilde{V}_1$ in $V$.
Wir sagen: $V$ besitzt die Eigenschaft \Llie, wenn gilt:
\begin{list}{}{}
\item a) $W$ ist eine rationale Unter-\VOA von $V$.
\item b) $\widetilde{V}_1$ ist ein Tensorprodukt von Kac-Moody \VOAs
$L_{k_i\Lambda_0^{(i)}}^{{\tx g}_i}$, 
wobei die ${\tx g}_i$ die einfachen
bzw.~abelschen Summanden der reduktiven Liealgebra $V_1={\tx g}=
\bigoplus_{i=1}^n{\tx g}_i$ sind.
\end{list}
\end{definition}
Hat eine \VOA $V$ die Eigenschaft \Llie{}, so ist sie eine direkte Summe
von Moduln der rationalen Unter-\VOA $\widetilde{V}_1\otimes W$.
Eine Schwierigkeit \Llie{} f"ur beliebige unit"are rationale \VOAs zu beweisen,
ist u.a.~die geforderte Rationalit"at von $W$.
\medskip

Schlie"slich sei noch der Mondscheinmodul $\VM$ aufgef"uhrt. Er wurde
in \cite{FLM} kon\-stru\-iert, um Teile des "`Mondschein des 
Monsters"'~\cite{CoNo}
zu erkl"aren, und war gleichzeitig das Beispiel, welches in der Mathematik
die Definition der \VOA-Struktur motivierte~\cite{Bo-ur,FLM-ur}.
Zus"atzlich ist er das erste Beispiel einer Orbifold-\VOA. Wir notieren das Hauptresultat aus~\cite{FLM} --- mit den Bezeichnungen wie dort --- als
\begin{satz}[Mondscheinmodul $\VM$]\label{monstermodul}
Der als $\Z_2$-Orbifold der Leechgitter-\VOA{} kon\-struier\-te 
Mond\-schein\-modul $\VM=V_{\Lambda_{24}}^+\oplus (V_{\Lambda_{24}}^T)^+$
ist eine uni\-t"are \VOA vom Rang $24$ mit dem Monster als 
Auto\-morphismen\-gruppe.
\end{satz}

\section{Modulgruppen und Modulfunktionen}

Wir f"uhren die in dieser Arbeit verwendeten Notationen und elementaren
Resultate im Zusammenhang mit Modulfunktionen ein. Dies, sowie die meisten der
ben"otigten expliziten Reihenentwicklungen finden sich in gr"o"serer 
Ausf"uhrlichkeit in Anhang I von \cite{HBJ}.

Sei $\H=\{ \tau\in\C\mid {\rm Im}(\tau)>0\}$ die obere Halbebene der komplexen
Zahlen. Die Gruppe $\slz$ operiert auf $\H$ verm"oge
$$ (A,\tau)\mapsto A\tau:=\frac{a \tau +b}{c \tau +d},\qquad \hbox{f"ur 
$A=\left({a\, b\atop c \, d}\right)\in \slz$ und $\tau\in\H$.} $$

Eine {\it Kongruenz\-unter\-gruppe} $\Gamma$ ist eine Unter\-gruppe von $\slz$, 
die eine Haupt\-kongruenz\-unter\-gruppe $\Gamma(N):=\{A \in \Gamma \mid 
A \equiv \left({1\, 0\atop 0\, 1}\right) \pmod{N}\}$ f"ur ein $N\in \N$ 
enth"alt.
Die Spitzen von $\Gamma$ sind die Bahnen der Operation von $\Gamma$ auf
${\bf P}_1(\Q)\subset {\bf P}_1(\R)$. Der Bahnenraum $\H/\Gamma$ kann durch
die Hinzunahme der Spitzen von $\Gamma$ zu einer kompakten Riemannschen
Fl"ache $\overline{\H/\Gamma}$ erweitert werden.

Die Operation von $\slz$ auf $\H$ induziert eine Operation von 
$\slz$ auf den Funktionen auf $\H$:
$$ f(\tau)\mapsto f\mid_A(\tau):=f(\frac{a\tau+b}{c\tau+d})\qquad
\hbox{f"ur ein $f:\H\longrightarrow\C$ und $A=\left({a\, b\atop c \,d}
\right)\in \slz$.} $$
Eine Modulfunktion zur Gruppe $\Gamma$ ist eine unter der $\mid$-Operation
invariante meromorphe Funktion auf $\H$, die in den Spitzen eine
Fourierreihenentwicklung besitzt und dort h"ochstens einen Pol hat. Dies sind
gerade die meromorphen Funktionen ${\cal M}(\overline{\H/\Gamma})$ auf der
Fl"ache $\overline{\H/\Gamma}$.

Im Zusammenhang mit selbstdualen \SVOAs sind wir insbesondere an $\slz$ selbst 
und der Thetagruppe $\GT:=\{A \in \slz \mid A \equiv 
\left({1\, 0\atop 0\, 1}\right) \hbox{oder} \left({0\, 1\atop 1\, 0}\right)
\pmod{2} \}$
interessiert.

Die Gruppe $\slz$ wird erzeugt von den beiden Elementen
$S=\left({0\, -1\atop 1\ \, 0}\right)$ und $T=\left({1\, 1\atop 0\, 1}\right)$.
Sie ist genauer das freie Produkt von $S$ und $T$ modulo der Relationen
$S^4=(ST)^6=1$ und $S^2=(ST)^3$.
Die inhomogene Gruppe $\pslz=\slz/\langle \pm { 1} \rangle$ ist das freie 
Produkt der ebenfalls wieder mit $S$ und $T$ bezeichneten Bildern von $S$ 
und $T$ in $\pslz$ zusammen mit den beiden Relationen $S^2=1$ und $(ST)^3=1$.

Die Fl"ache $\overline{\H/\slz}$ hat eine mit $\infty$ bezeichnete
Spitze und ist vom Geschlecht~$0$. Der K"orper der meromorphen Funktionen
von $\overline{\H/\slz}$ ist daher ein rationaler Funktionenk"orper.
Die Modulfunktion $j=\frac{E_4^3}{\Delta}$ zu $\slz$ ist holomorph in
$\H$ und besitzt in der Spitze $\infty$  (lokale Koordinate $q=e^{2\pi i\tau}$) 
die Fourierreihenentwicklung
\begin{equation}\label{modulfktj}
j=q^{-1} + 744 + 196884\, q + 21493760 \, q^2 + \cdots ,
\end{equation}
d.h.~es gilt ${\cal M}(\overline{\H/ \Gamma})\cong\C(j)$.
Hierbei ist $E_4=1 + 240 \sum_{n=1}^{\infty} \sigma_3(n)\, q^n$ eine 
Eisensteinreihe vom Gewicht $4$ zu $\slz$, und die Diskriminante $\Delta$
ist die $24$-te Potenz der $\eta$-Funktion, die die Produktentwicklung
$$\eta=q^{\frac{1}{24}}\prod_{n=1}^{\infty} (1-q^n)$$
besitzt.

Die Thetagruppe $\GT$ (vgl.~\cite{rankin}) wird erzeugt von $S$ und $T^2$.
Sie ist das freie Produkt von $S$ und $T^2$ modulo der 
Relationen $S^4=1$ und $S^2T^2=T^2S^2$. 
Die inhomogene Thetagruppe $\overline{\Gamma}_{\theta}=
\GT/\langle \pm {1} \rangle\subset\pslz $ ist das freie Produkt von
$S$ und $T^2$ zusammen mit der Relation $S^2=1$. Die Kongruenzuntergruppe
$\GT$ hat Index $3$ in $\slz$, und die drei Nebenklassen $\slz/\GT$ werden
repr"asentiert von $1$, $T$ und $ST$.
Die Fl"ache $\overline{\H/\GT}$ hat ebenfalls das Geschlecht $0$ und besitzt
die beiden Spitzen $1$ und $\infty$.
Die Funktion $\jt=\left(\frac{\Theta_{\Z}}{\eta}\right)^{12}$, wobei 
$\Theta_{\Z}=\sum_{n\in \Z} q^{\frac{1}{2} n^2}$ 
die Thetareihe zum Gitter $\Z$ ist, ist
eine Modulfunktion zu $\GT$. Sie hat in der Spitze $\infty$ 
(lokale Koordinate $q^{\frac{1}{2}}$) die Produktentwicklung
\begin{equation}\label{modulfktjt}
\jt=
\left(\frac{\eta^2(\tau)}{\eta(2\tau)\eta(\tau/2)}\right)^{24}=
q^{-\frac{1}{2}}\prod_{n=0}^{\infty}(1+q^{n+\frac{1}{2}})^{24},
\end{equation}
hat dort also einen einfachen Pol. In der Spitze $1$ hat sie die Entwicklung
\begin{equation}\label{spitze-1}
\jt=
-2^{12}\left(\frac{\eta(2\tau)}{\eta(\tau)}\right)^{24}=
-2^{12}\cdot q
\prod_{n=1}^{\infty}(1+q^n)^{24},
\end{equation}
d.h.~dort hat sie eine einfache Nullstelle.
Aufgrund der Produktentwicklung ist sie eine in $\H$ holomorphe Funktion,
und es gilt ${\cal M}(\overline{\H/\GT})\cong\C(\jt)$.

\medskip
Mit $\Theta_L=\sum_{x\in L} q^{\frac{1}{2}\langle x,x\rangle}$ bezeichnen 
wir die Thetareihe eines ganzen Gitters $L\subset {\R}^n$. Sie ist eine
Modulform vom Gewicht $\frac{n}{2}$ zur Thetagruppe $\GT$.

Wir werden eine Modulfunktion oft mit ihrer Reihen\-entwick\-lung in einer 
Spit\-ze iden\-tifizie\-ren.

\section{\OVOAs bei Geschlecht Eins}

Die Beschreibung einer konformen Quantenfeldtheorie f"ur Riemannsche Fl"achen von h"oherem Geschlecht sollte sich aus der Beschreibung bei
Geschlecht Null ergeben, denn Fl"achen gr"o"seren Geschlechtes lassen sich 
durch Zusammenkleben der dreifach punktierten Zahlenkugel erhalten 
(vgl.~\cite{MoSei}).

In diesem Abschnitt stellen wir die Resultate zusammen, die sich aus der 
Theorie der Vertexoperator-Algebren (mathematische Beschreibung von \CFT 
bei Geschlecht $0$) f"ur Riemansche Fl"achen bei Geschlecht $1$ (Tori) ergeben. Die wesentlichen S"atze finden sich in der Arbeit~\cite{Zhu-dr} von Zhu.
Wir formulieren weiter die entsprechenden Resultate f"ur \SVOAs.

\medskip
Sei $V$ eine \nette rationale \VOA vom Rang $c$ und $M$ ein irreduzibler $V$-Modul vom konformen Gewicht $h$.
F"ur ein $n$-Tupel $(a_1,a_2,\ldots,a_n)$ von homogenen Elementen $a_i\in V$ 
definiert die formale $q$-Spur\pagebreak[2]
\begin{equation}\label{kor-def1}
 F_{M}((a_1,x_1),(a_2,x_2),\ldots,(a_n,x_n),q):=\qquad\qquad\qquad\qquad
\qquad\qquad $$ $$\qquad\qquad\qquad\qquad
{x_1}^{{\rm deg\,}a_1}\cdot \cdots \cdot {x_2}^{{\rm deg\,}a_n}
\,{\rm tr}\mid_{M} Y(a_1,x_1)Y(a_2,x_2)\ldots Y(a_n,x_n)q^{L_0}
\end{equation}\pagebreak[2]
ein Element in $q^h\C[[x_1,x_1^{-1},x_2,x_2^{-1},\ldots,x_n,x_n^{-1},q]]$.
\hfill\break
Die Reihe $F_{M}((a_1,x_1),(a_2,x_2),\ldots,(a_n,x_n),q)$ konvergiert 
nach Satz 4.2.1  aus~\cite{Zhu-dr} zu einer in dem Gebiet
$$\{(x_1,x_2,\ldots,x_n,q)\mid 1>|x_1|>\cdots|x_n|>|q|\}$$
holomorphen Funktion, und der Limes kann zu einer in dem Gebiet
$$\{(x_1,x_2,\ldots,x_n,q)\mid x_i\not=0 , |q|<1\}$$
meromorphen Funktion $\widetilde{F}_{M}$ fortgesetzt werden.
Wir ersetzen die Variable $x_i$ durch $e^{2\pi i z_i}$, $q$ durch
$e^{2\pi i\tau}$ und bezeichnen 
\begin{equation}\label{kor-def2}
T_{M}((a_1,z_1),\ldots,(a_n,z_n),\tau)=
q^{-\frac{c}{24}}\widetilde{F}_{M}((a_1,e^{2\pi i z_1}),
\ldots,(a_n,e^{2\pi i z_n}),e^{2\pi i\tau})
\end{equation}
als die \it $n$-Punkt-Korrelations\-funktion der Vertex\-opera\-toren 
$Y(a_1,x_1)$, $\ldots$, $Y(a_n,x_n)$ auf dem To\-rus mit Para\-meter $\tau$.
\rm Die Funktion $T_{M}((a_1,z_1),(a_2,z_2),\ldots,(a_n,z_n),\tau)$ ist
nach Satz 4.5.1 aus~\cite{Zhu-dr} 
eine in jeder Variablen~$z_i$ doppeltperiodische Funktion mit den Perioden $1$ 
und $\tau$, und die einzigen m"oglichen Singularit"aten sind die Stellen
$z_i=z_j+m+k\tau$ ($i\not=j$; $m$, $k\in \Z$).
Weiter h"angt die Funktion nicht von der Reihenfolge der Vertexoperatoren
$Y(a_i,x_i)$ ab, d.h.~f"ur jede Permutation $\sigma\in S_n$ gilt
$$T_{M}((a_1,z_1),(a_2,z_2),\ldots,(a_n,z_n),\tau)=
T_{M}((a_{\sigma(1)},z_{\sigma(1)}),(a_{\sigma(2)},z_{\sigma(2)}),
\ldots,(a_{\sigma(n)},z_{\sigma(n)}),\tau).$$
\begin{definition}\label{konformer-Block}
Sei $V$ eine rationale \nette \VOA und $M_1$, $\ldots$, $M_m$ eine vollst"andige
Liste von irreduziblen $V$-Moduln. Den $\C$-Vektorraum $B$, der
von den $m$~Funktionen 
\hbox{$T_{M_i}\,:\,\bigcup_{n=1}^{\infty}((V\times \C)^n\times 
\H)\longrightarrow\C$,($i=1$, $\ldots$, $m$)
} aufgespannt wird, bezeichnen wir
als den {\em konformen Block auf dem Torus.}
\end{definition}
{\it Anmerkung:} Unsere Definition des konformen Blockes auf dem Torus
stimmt mit derjenigen in~\cite{Zhu-dr} nach den dort bewiesenen
S"atzen "uberein.
Insbesondere sind die $m$~Funktionen $T_{M_i}$ linear unabh"angig, d.h.~$B$ 
ist ein $m$-dimensionaler $\C$-Vektorraum mit der ausgezeichneten Basis
$\{T_{M_i}\}$. Das Hauptresultat aus~\cite{Zhu-dr} (Satz 5.3.2) ist:
\begin{satz}[Zhu]\label{trans-charvoa}
Es existiert eine Darstellung $\rho:\slz\longrightarrow{\rm End}(B)$ der 
Modulgruppe, so da"s f"ur Virasoroh"ochstgewichtsvektoren $a_1$, $a_2$, 
$\ldots$, $a_n$ und $A=\left({a\, b\atop c \, d}\right)\in\slz$ 
die \hbox{$n$-Punkt-Korrelationsfunktionen} das folgende
Transformationsverhalten besitzen:
$$T_{M_i}((a_1,\frac{z_1}{c\tau +d}),(a_2,\frac{z_2}{c\tau +d}),
\ldots,(a_n,\frac{z_n}{c\tau +d}),\frac{a\tau+b}{c\tau +d})=
\qquad\qquad\qquad\qquad\qquad$$ $$\qquad\qquad\qquad
(c\tau+d)^{\sum_{i=1}^n{\deg a_i}}\sum_{j=1}^m\rho(A)_{ij}
T_{M_j}((a_1,z_1),(a_2,z_2),\ldots,(a_n,z_n),\tau).$$
Hierbei ist $(\rho(A)_{ij})$ die $m\times m$-Matrix der in der Basis 
$\{T_{M_i}\}$ geschriebenen Abbildung $\rho(A)\in{\rm End}(B)$.
\end{satz}
Der Vektor $\{T_{M_i}((a_1,z_1),\ldots,(a_n,z_n),\tau)\}$ kann also als
eine vektorwertige meromorphe Jacobiform vom Gewicht $\sum_{i=1}^n{\deg a_i}$
zur vollen Modulgruppe {\small $\Z^2>\!\!\!\!\lhd\,\slz$} aufgefa"st werden 
(vgl.~\cite{EiZa}).

F"ur $n=1$ hat die in $z_1$ elliptische Funktion $T_{M}((a_1,z_1),\tau)$ 
keine Polstellen, ist also konstant bzgl.~$z_1$.
Wir bezeichnen 
$$\chi_{M}:=T_{M}(({\bf 1},z_1),\tau)=q^{-\frac{c}{24}}
\sum_{n=1}^{\infty}{\rm dim}\, M_n\,q^n={\rm tr}|_M\, q^{L_0-\frac{c}{24}}$$
als den {\it (konformen) Charakter} eines $V$-Moduls $M$. Er ist eine in
$\tau\in\H$ holomorphe Funktion.
\begin{korollar}[Zhu]\label{trans-charvoa-korollar}
Die nat"urliche Projektion $\pi:B\longrightarrow \overline{B}:={\rm Span}
\{\chi_{M_i}\}$ induziert eine Darstellung
$$\overline{\rho}:\pslz\longrightarrow{\rm End}(\overline{B})$$
der inhomogenen Modulgruppe auf dem Vektorraum der konformen Charaktere.
\end{korollar}
Im allgemeinen ist $\pi$ nicht injektiv, z.B.~gilt f"ur ein $M_i$ und seinen
dualen Modul $M_i^*=M_{i'}$ die Beziehung $\chi_{M_i}=\chi_{M_{i'}}$, aber
$M_i$ und $M_{i'}$ sind nicht notwendig isomorph.
Wir schreiben f"ur ein Element $A\in\slz$ kurz $\widetilde{A}=\rho(A)$ und
$\overline{A}=\overline{\rho}(A)$.

Verlinde hat in~\cite{verlinde} den folgenden Zusammenhang zwischen der
Fusionsalgebra und der Darstellung $\rho$ aus Satz~\ref{trans-charvoa}
hergestellt:
\begin{vermutung}[Verlindeformel]\label{verlindeformel}
Die Dimensionen ${N_{ij}^k}$ der Intertwinerr"aume einer einfachen \netten 
rationalen \VOA $V$ mit irreduziblen Moduln $M_1=V$, $M_2$, $\ldots$, $M_m$ 
lassen sich durch die $m\times m$-Matrix 
$\widetilde{S}=\rho(S)$, $S=\left({0\, -1\atop 1\ \, 0}\right)$ ausdr"ucken:
$$ N_{ij}^k=\sum_{n=1}^m\,\frac{\widetilde{S}_{in}\widetilde{S}_{jn}
\widetilde{S}^{-1}_{nk}}{\widetilde{S}_{1n}}.$$
\end{vermutung}
F"ur beliebige \nette rationale \VOAs ist die Verlindeformel bisher noch nicht
bewiesen. Die Definiton des Raumes der Vakua in~\cite{Zhu-riem} f"ur
Riemannsche Fl"achen h"oheren Geschlechtes erlaubt die allgemeine Formulierung 
der Verlindeformel im Rahmen von \VOAs.

\medskip
Analog zu \VOAs lassen sich auch f"ur \SVOAs die Korrelationsfunktionen 
auf dem Torus definieren und entsprechende Resultate erzielen. Wir werden
daher hier nur die Aussagen formulieren und an einigen Stellen etwas zur
Beweismethode sagen. Die nachfolgenden Kapitel sind auch unabh"angig von diesen
"Uberlegungen, falls in der Definiton der \SVOA zus"atzlich das entsprechende 
Modultransformationsverhalten der Charaktere gefordert wird.

Die Definiton der $n$-Punkt Korrelationsfunktion auf dem Torus ist f"ur eine
\SVOA $V=V_{(0)}\oplus V_{(1)}$ die gleiche wie in~(\ref{kor-def1}) 
und~(\ref{kor-def2}), allerdings folgen nur f"ur \hbox{$a_1$, $\ldots$, $a_n 
\in V_{(0)}$} die Konvergenz und die analytische Fortsetzbarkeit direkt aus dem 
Resultat f"ur \VOAs. F"ur beliebige $a_i\in V$  ist 
$T_M((a_1,z_1),(a_2,z_2),\ldots,(a_n,z_n),\tau))$ eine in den Variablen
$z_i$ doppeltperiodische Funktion mit den beiden
Perioden $1$ und $2\tau$. Beim "Ubergang $z_i\mapsto z_i+\tau$
transformiert sich $T_M$ mit dem Vorzeichen $(-1)^{|a_i|}$. F"ur eine beliebige
Permutation $\sigma\in S_n$ gilt
$$T_{M}((a_1,z_1),\ldots,(a_n,z_n),\tau)=(-1)^w\cdot
T_{M}((a_{\sigma(1)},z_{\sigma(1)}),\ldots,
(a_{\sigma(n)},z_{\sigma(n)}),\tau),$$
wobei $w$ die Anzahl der Vertauschungen von ungeraden Elementen $a_i\in V_{(1)}$
in $\sigma$ ist. 
Entsprechend 
dem konformen Block definiert man den {\it superkonformen Block:}

\begin{definition}\label{superkonformer-Block} 
Sei $V$ eine ra\-tio\-nale \nette \SVOA und sei $M_1$, $\ldots$, $M_m$ eine 
voll\-st"an\-dige Lis\-te von irre\-duziblen $V$-Mo\-duln. 
Der {\em superkonforme Block auf dem Torus} ist der von den $m$ Funktionen 
$T_{M_i}\,:\,\bigcup_{n=1}^{\infty}((V\times \C)^n\times \H)
\longrightarrow\C$,($i=1$, $\ldots$, $m$) aufgespannte $\C$-Vektorraum $SB$. 
\end{definition}
Die "`Super-Erweiterung"' von Satz~\ref{trans-charvoa} ist
\begin{satz}\label{trans-charsvoa}
Es existiert eine Darstellung $\rho:\Gamma_{\theta}
\longrightarrow{\rm End}(SB)$ der 
Thetagruppe, so da"s f"ur Virasoroh"ochstgewichtsvektoren $a_1$, $a_2$, 
$\ldots$, $a_n$ und $A=\left({a\, b\atop c \, d}\right)\in\Gamma_{\theta}$ 
die $n$-Punkt-Korrelationsfunktionen das folgende
Transformationsverhalten besitzen:
$$T_{M_i}((a_1,\frac{z_1}{c\tau +d}),(a_2,\frac{z_2}{c\tau +d}),
\ldots,(a_n,\frac{z_n}{c\tau +d}),\frac{a\tau+b}{c\tau +d})=
\qquad\qquad\qquad\qquad\qquad$$ $$\qquad\qquad\qquad
(c\tau+d)^{\sum_{i=1}^n{\deg a_i}}\sum_{j=1}^n\rho(A)_{ij}
T_{M_j}((a_1,z_1),(a_2,z_2),\ldots,(a_n,z_n),\tau).$$
Hierbei ist $(\rho(A)_{ij})$ die $m\times m$-Matrix der in der Basis 
$\{T_{M_i}\}$ geschriebenen Abbildung $\rho(A)\in{\rm End}(SB)$.
\end{satz}
F"ur die Charaktere gilt das 
\begin{korollar}\label{trans-charsvoa-korollar}
Die nat"urliche Projektion $\pi:SB\longrightarrow \overline{SB}:={\rm Span}
\{\chi_{M_i}\}$ induziert eine Darstellung
$$\overline{\rho}:\overline{\Gamma}_{\theta}\longrightarrow
{\rm End}(\overline{B})$$
der inhomogenen Thetagruppe auf dem Vektorraum der konformen Charaktere.
\end{korollar}
{\it Anmerkung:\/} 
Vertexoperatoren $Y(a,x)$ zu Virasoro\-h"ochst\-gewichts\-vektoren $a\in V_k$ 
trans\-formieren sich wie Schnit\-te in der 
$2k$-ten Tensor\-potenz des Spinor\-b"undels 
"uber dem Torus  $\C/(\Z+\Z\tau)$ zu der durch 
$\kappa:(1,\tau)\mapsto (1,-1)$ ge\-gebenen Spin-Struktur. 
Zwei Tori mit dieser Spin-Struktur sind genau dann konform "aquivalent, wenn
die Modulparameter durch $\Gamma_{\theta}$ ineinander "ubergef"uhrt werden
k"onnen. Um zu einer Darstellung von $\slz$ zu gelangen, mu"s man alle
drei nichttrivialen (ungeraden) Spin-Strukturen betrachten und den superkonformen
Block $SB$ um die Korrelationsfunktionen von getwisteten $V$-Moduln
erweitern (vgl.~S.~\pageref{twistsektor}). Korrelationsfunktionen
zur geraden Spin-Struktur liefern ebenfalls eine $\slz$-Darstellung.

{\bf Beweis von Satz~\ref{trans-charsvoa} (grobe Skizze):}
Kac und Wang haben in~\cite{KacWang} das Analogon der Zhu'schen Algebra
$A(V)$ f"ur \SVOAs
definiert und gezeigt, da"s eine $1:1$-Korrespondenz zwischen 
irreduziblen $A(V)$-Moduln und irreduziblen $V$-Moduln besteht.
Dies ist der "`Geschlecht $0$ Teil"'.

F"ur den "`Geschlecht $1$ Teil"' kann man ebenfalls analog zu der Arbeit
von Zhu~\cite{Zhu-dr} vorgehen und zuerst den superkonformen Block
auf dem Torus "`axiomatisch"' definieren und dann zeigen:\hfill\break
(1) Die Korrelationsfunktionen auf dem Torus liegen in dem superkonformen 
Block. \hfill\break
(2) Bei Transformationen mit $\GT$ bleibt man im superkonformen 
Block. \hfill\break
(3) Ein Element im superkonformen 
Block l"a"st sich als Spur auf $A(V)$ ausdr"ucken.

Dazu mu"s man zeigen, da"s man wegen der Endlichkeits\-bedingung ${\bf C}_2$
eine Differential\-gleichung erh"alt, deren Koeffizienten
elliptische Funktionen zu dem Paar Kurve mit festem $2$-Teilungs\-punkt
sind. \qed
\medskip

Man erh"alt leicht die folgenden Eigenschaften des (super)konformen Blockes
bei \VOA-kategoriellen Konstruktionen.

Seien $V_1$ und $V_2$ zwei \nette rationale \OVOAs. F"ur den als 
$\slz$-Modul (bzw.~$\Gamma_{\theta}$-Modul) aufgefa"sten (super)konformen Block
des Tensorproduktes $V_1 \otimes V_2$ gilt
$$  (S)B_{V_1 \otimes V_2} \cong (S)B_{V_1}\otimes (S)B_{V_2}.$$
Dies folgt unmittelbar aus Satz~\ref{tensorprodukt} und 
Definition~\ref{konformer-Block} bzw.~\ref{superkonformer-Block}.

Eine Einbettung $i:U\subset V$ von zwei \netten rationalen \OVOAs $U$ 
und $V$ induziert die nat"urliche Abbildung 
$i^*\,:\,(S)B_V\longrightarrow (S)B_U$ 
zwischen $\slz$- (bzw.~$\Gamma_{\theta}$-) Moduln, die vermutlich injektiv ist.
Schlie"slich liefert die Inklusion $i:W_{(0)}\subset W$ der Unter-\VOA $W_{(0)}$
in einer \netten rationalen \SVOA $W$ das folgende kommutative Diagramm:
$$\begin{array}{ccc}
\Gamma_{\theta}\times SB_{W} & {\textstyle \nu\,\times\, i^*}\atop
{\textstyle \longrightarrow} & \slz\times B_{W_{(0)}} \\
\rho\downarrow &  &\rho\downarrow \\
 SB_{W} & {\textstyle i^*}\atop {\textstyle \longrightarrow} & B_{W_{(0)}}. 
\end{array}$$
Hierbei ist $\nu$ die Inklusion $\Gamma_{\theta}\subset\slz$, und $i^*:
SB_W\longrightarrow B_{W_{(0)}}$ bildet die Funktion
$T_M:\bigcup_{n=1}^{\infty}((W\times \C)^n,\H)\longrightarrow\C$
auf die Einschr"ankung $i^*(T_M)=T_{M\mid i}\circ i:
\bigcup_{n=1}^{\infty}(\hbox{$(W_{(0)}\times\C)^n$},\H)
\longrightarrow\C$ ab. Da $W_{(0)}$ auch rational ist, zerlegt sich
$M|i$ als direkte Summe von irreduziblen $W_{(0)}$-Moduln.

Alle drei Konstruktionen sind mit der Projektion $\pi$ auf die konformen
Charaktere vertr"aglich.

\medskip
Wir notieren noch die konformen Charaktere der wichtigsten Beispiele aus dem 
vor\-letzten Abschnitt. Diese ergeben sich aus der De\-fini\-tion der der \OVOA
und ihren Moduln unter\-liegenden Vektor\-r"aume.
\label{beispielcharaktere}

Die Charaktere der $3$ irreduziblen $L_{1/2}(0)$-Moduln sind:
\begin{eqnarray}\label{char-vir}
\chi_{L_{1/2}(0)}& = & \frac{1}{2}\left(
\sqrt{\frac{\Theta_{\Z}}{\eta}}+\sqrt{\frac{\Theta_{\Z}}{\eta}}\mid_T\right)
=q^{- {1 \over 48}}(
1 + {q^2} +  {q^3} + + 2\,{q^4} + 2\,{q^5} +\cdots \,),
\nonumber \\
\chi_{L_{1/2}(\frac{1}{2})} & = & 
 \frac{1}{2}\left(
\sqrt{\frac{\Theta_{\Z}}{\eta}}-\sqrt{\frac{\Theta_{\Z}}{\eta}}\mid_T\right)
= q^{- {1 \over 48}}(
{q^{{1\over 2}}} + {q^{{3\over 2}}} + {q^{{5\over 2}}} +  
  {q^{{7\over 2}}} + 2\,{q^{{9\over 2}}} + 
  2\,{q^{{{11}\over 2}}} +\cdots \,),
\nonumber \\
\chi_{L_{1/2}(\frac{1}{16})}& = & 
\frac{1}{\sqrt{2}}\sqrt{\frac{\Theta_{\Z+\frac{1}{2}}}{\eta}}
= q^{- {1 \over 48}}(
{q^{{1\over {16}}}} + {q^{{{17}\over {16}}}} + {q^{{{33}\over {16}}}} + 
  2\,{q^{{{49}\over {16}}}} + 2\,{q^{{{65}\over {16}}}} + +\cdots \,).
\end{eqnarray}
F"ur den Charakter der Fermi-\SVOA $\VF\cong L_{1/2}(0)\oplus 
L_{1/2}(\frac{1}{2})$ ergibt sich somit:
\begin{equation}\label{char-fermi}
\chi_{\VF}=\sqrt{\frac{\Theta_{\Z}}{\eta}}=\sqrt[24]{j_{\theta}}.
\end{equation}
Die Charaktere der irreduziblen Moduln $M_{L+[i]}$, $[i]\in L^*/L$
einer Gitter-\OVOA $V_L$ zum ganzen Gitter $L\subset\R^n$ sind
$$\chi_{M_{L+[i]}}=\frac{\Theta_{{L+[i]}}}{\eta^n}.$$

F"ur die Moduln $V_{k\Lambda_0+\lambda}^{\tx g}$ zum H"ochstgewicht $\lambda$
der zu einer Kac-Moody Algebra $\hat{{\tx g}}$ assoziierten \VOA bei 
Stufe $k$
ergibt sich der Charakter aus der Weil-Kac Charakterformel (s.~\cite{Kac}, 
Prop.~10.10, Bezeichnungen wie dort):
$$\chi_{V_{k\Lambda_0+\lambda}^{\tx g}}
=q^{-\frac{c}{24}}\prod_{\alpha\in\Delta_+^{\vee}}
\left(\frac{1-q^{\langle k\Lambda_0+\lambda+\rho,\alpha\rangle}}
{1-q^{\langle\rho,\alpha\rangle}}\right)^{{\rm mult}\,\alpha}.$$
Schlie"slich erh"alt man f"ur den Charakter des Mondscheinmoduls
$\VM$ (s.~\cite{FLM}, (10.5.55)):
$$\chi_{\VM}=j(q)-744=J(q).$$

\chapter{Selbstduale Vertexoperator-Superalgebren}

\vspace{-4mm}
F"ur selbstduale gerade und ungerade Codes erh"alt man aus der Mac\-Williams 
Iden\-ti\-t"at~\cite{MacSl} mit Hilfe von Invarianten\-theorie eine 
Be\-schreibung des Gewichts\-z"ahler\-polynomes des Codes. F"ur selbstduale gerade und ungerade Gitter
ergibt sich aus der Jacobi-Umkehrformel und der Theorie der Modulformen eine 
Beschreibung der Thetareihe des Gitters. In diesem Kapitel geben wir eine
analoge Beschreibung f"ur den Charakter von selbstdualen \VOAs und \SVOAs.
Hierzu benutzen wir einerseits den Satz von Zhu (\ref{trans-charvoa}) 
bzw.~seine Verallgemeinerung \ref{trans-charsvoa} f"ur \SVOAs und andererseits
elementare Eigenschaften von meromorphen Funktionen. F"ur \SVOAs stellt sich
die Bedingung unit"ar als wichtig heraus. Ein Ganzzahligkeitsargument 
erlaubt es uns weiter zu zeigen, da"s der Rang einer selbstdualen \SVOA ein 
Element von $\frac{1}{2}\Z$ ist, was schon in \cite{verlinde} vermutet 
worden war.

Selbstduale Codes und Gitter besitzen viele interessante Eigenschaften,
so da"s auch die Klassifkation von selbstdualen \VOAs nat"urlich
erscheint. Die Klassifikation der $24$-dimensionalen selbstdualen geraden
Gitter steht in Verbindung mit den verschiedenen Typen von tiefen L"ochern im 
Leechgitter~\cite{CoSl}, Kapitel 23 und 24,
und den Bahnen von lichtartigen Vektoren im $26$-dimensionalen
selbstdualen Lorentzgitter unter der Operation der
Automorphismengruppe~\cite{CoSl}, Kapitel 26 und 27.
Wir geben hier f"ur R"ange $c\leq 24$ eine Liste \cite{schellekens1} an, von der
vermutet wird, da"s sie alle selbstdualen unit"aren \VOAs beschreibt. 
Die analoge Klassifikation f"ur selbstduale Gitter war von 
Kneser~\cite{kneser} 
und Niemeier~\cite{niemeier} erhalten worden. Weiter
geben wir in diesem Kapitel eine Beschreibung der selbstdualen unit"aren
\SVOAs vom Rang kleiner $8$: f"ur jeden Rang gibt es genau eine \SVOA.  Mit den Methoden des n"achsten  Kapitels wird es m"oglich sein, die selbstdualen 
unit"aren \SVOAs vom Rang kleiner $24$
zu klassifizieren~\cite{hoehn2}. Viel weiter wird man in der Klassifikation
nicht kommen k"onnen, da die {\it Minkowski-Siegelsche 
Ma"sformel\/}~\cite{dir,min,sie,weil} zeigt, da"s die Anzahl der selbstdualen
Gitter und damit auch die Anzahl der \OVOAs jenseits von $24$ sehr stark 
anw"achst.
\pagebreak[2]

Eine interessante Fragestellung ist es, eine analoge Ma"s\-for\-mel auch f"ur 
\VOAs zu finden. Pro\-bleme hier\-bei sind, da"s man gleich\-zeitig end\-liche 
und unend\-liche 
Auto\-morphismen\-gruppen ber"ucksichtigen mu"s, da"s man keinen
Raum kennt (${\bf F}_2^n$ bei Codes, $\R^n$ bei Gittern), in den sich 
alle \VOAs festen Ranges einbetten lassen, und da"s man keine offensichtliche
Gruppe hat, die operiert ($S_n$ bei Codes, 
${\bf SO}_n({\bf A_Q})$ die orthogonale Gruppe "uber den Adelen bei
Gittern $\cong$ quadratischen Formen). Insbesondere w"are eine "Uber\-tragung
des Begriffes {\it Geschlechter\/}~\cite{cassels,are} von quadratischen Formen
auf \VOAs interessant. Zu erwarten ist dann zumindest, da"s die 
{\it Klassen\-zahl\/}, d.h.~die Anzahl der Isomorphie\-klassen von \VOAs
\hbox{von festem Geschlecht (z.B.~dem selbstdualen) endlich ist.}

\vspace{-6mm}

\section{Selbstduale \VOAs{ }}\vspace{-4mm}

Die Struktur der in \ref{deffusion} definierten Fusionsalgebra
kann recht kompliziert sein, es sind keine allgemeinen 
Klassifikationss"atze bekannt. Es erscheint daher sinnvoll, insbesondere die 
\OVOAs zu untersuchen, f"ur die die Fusionalgebra trivial ist. 
Wie das Beispiel des Mondscheinmoduls zeigt, k"onnen diese \OVOAs trotzdem
eine interessante Struktur besitzen.
Wir machen daher die \vspace{-2mm}
\begin{definition}\label{selbstdual}
Eine \OVOA $V$ hei"st selbstdual, wenn  sie bis auf Isomorphie genau einen
irreduziblen Modul besitzt, n"amlich $V$ als adjungierten Modul selbst.
\end{definition}
\vspace{-2mm}
Anmerkungen: Au"ser der Be\-zeichnung {\it selbstdual\/} werden auch die
Be\-zeich\-nun\-gen {\it holo\-morph\/} \cite{DM-M24} oder {\it mero\-morph\/}
\cite{schellekens1} ver\-wendet. Wegen der Analogie zu Codes und Gitter 
verwenden wir in dieser Arbeit die Bezeichnung selbstdual.

Eine selbstduale \OVOA ist insbesondere einfach, da sie als ad\-jungier\-ter
Mo\-dul irre\-duzi\-bel ist.

Nach Satz~\ref{satzgittersvoa} ist die Fusionsalgebra einer \SVOA $V_L$ zu
einem ganzen Gitter $L$ isomorph zum Gruppenring $\Z[L^*/L]$, d.h.~$V_L$ ist
genau dann selbstdual, wenn $L$ selbstdual ist.

\medskip

Im Rest dieses Abschnittes betrachten wir selbstduale \VOAs und im
n"achsten selbstduale \SVOAs. Alle Aussagen dieses Abschnittes finden 
sich auch in der physikalischen Literatur, sind dort teilweise aber nicht 
bewiesen oder nur ohne genaue Angabe der Voraussetzung.

In Verallgemeinerung der entsprechenden S"atze f"ur gerade Codes (Gleason,
vgl.~\cite{MacSl}, Kap.~19)
und Gitter (Hecke, vgl.~\cite{CoSl}, Kap.~7.6) gilt der 

\vspace{-2mm}
\pagebreak[2]\vbox{
\begin{satz}[vgl.~z.B.~\cite{Go-mero}]\label{charsdvoa}
Sei $V$ eine selbstduale \nette rationale \VOA vom Rang $c$. Dann ist der
Charakter $\chi_V=q^{-c/24}\sum_{n=0}^{\infty}\dim V_n\, q^n$ ein homogenes 
Polynom $P(x,y)$ vom Gewicht $c$ in den Funktionen $x=\sqrt[3]{j}$
(Gewicht $8$) und $y=1$ (Gewicht~$24$), oder --- "aquivalent dazu ---
ein homogenes Polynom in $x=\chi_{V_{E_8}}$ und $y=\chi_{\VM}$, den
Charakteren der selbstdualen \VOA zum Gitter $E_8$ und der selbstdualen
Monster-\VOA $\VM$.
\end{satz}}\pagebreak[3]
\begin{korollar}
Der Rang $c$ einer selbst\-dualen \netten rationalen \VOA ist ein 
ganz\-zahliges Viel\-faches von 8.
\end{korollar}
Falls $V$ unit"ar ist, gilt nach Satz~\ref{lemmaunitaer} nat"urlich $c\geq 0$.

{\bf Beweis:}
Da $V$ nur einen irreduziblen Modul besitzt, liefert nach 
Korollar~\ref{trans-charvoa-korollar}
zum Satz von Zhu der von 
$\chi_V=q^{-c/24}\sum_{n=0}^{\infty}\dim V_n q^n$ aufgespannte 
$\C$-Vektorraum eine eindimensionale Darstellung von $\pslz$.
F"ur die Erzeuger 
$S=\left({0\, -1\atop 1\ \, 0}\right)$ und $T=\left({1\, 1\atop 0\, 1}\right)$
hat man das folgende Transformationsverhalten:
\begin{eqnarray*}
\chi_V(\tau)\mid_T& =& e^{-\frac{2 \pi i c}{24}}\cdot\chi_V(\tau), \\
\chi_V(\tau)\mid_S& =& \chi_V(\tau).
\end{eqnarray*}
Die zweite Gleichung gilt wegen $S^2=\mbox{id}$ und 
$\chi_V(i)=\chi_V(i)\mid_S=
\sum_{n=0}^{\infty}\dim V_n \cdot e^{- 2 \pi (- c/24+n)}>0$.
Wegen der Relation $(ST)^3=\mbox{id}$ in $\pslz$ folgt 
$e^{2 \pi i \frac{c}{8}}=1$, d.h.~$c\in8\Z$. 
Daher ist $\chi_V(\tau)$ invariant unter $T^3$.

Die von $S$ und $T^3$ erzeugte Untergruppe $\Gamma$ hat Index $3$ in $\pslz$
und die zugeh"orige Fl"ache $\overline{\H/\Gamma}$ hat das Geschlecht $0$.
Den K"orper der meromorphen Funktionen bilden die rationalen Funktionen in
$\sqrt[3]{j}= \mbox{$q^{-(1/3)}(1+248 q +\cdots)$}\in q^{-(1/3)}\Z[[q]]$, der
$3$-ten Wurzel der absoluten Modulinvariante $j$ (siehe~\cite{schoen}, 
S.~131).
Die Funktion $\sqrt[3]{j}$ hat einen einfachen Pol in der Spitze $\infty$
(lokale Koordinate $q^{1/3}$). Der Charakter
\hbox{$\chi_V\in q^{-\frac{1}{3}\cdot \frac{c}{8}}\Z[[q]]$} ist nach 
Satz~\ref{trans-charvoa}
holomorph in der oberen Halbebene und hat einen Pol der Ordnung 
$\frac{c}{8}$ in $\infty$. Er ist daher ein Polynom
vom Grad $\frac{c}{8}$ in $\sqrt[3]{j}$, in dem wegen der $q$\,-Entwicklung nur
Exponenten kongruent $c/8 \pmod{3}$ vorkommen. Hieraus folgt der Satz.
Die Umformulierung gilt wegen $\chi_{V_{E_8}}=\sqrt[3]{j}$ und
$\chi_{\VM}=\left(\sqrt[3]{j}\right)^3-744$. \qed

{\bf Beispiele von selbstdualen \netten rationalen \VOAs:}
 
Positiv definite gerade selbstduale Gitter liefern nach 
Konstruktion~\ref{satzgittersvoa} die folgenden Beispiele von selbstdualen
unit"aren \VOAs: F"ur $c=8$: $V_{E_8}$, f"ur $c=16$: $V_{E_8}\otimes
V_{E_8}$, $V_{D_{16}^+}$ und f"ur $c=24$ die $24$ \VOAs, die zu den $24$ 
Niemeiergittern~\cite{niemeier} geh"oren. F"ur $c\geq 32$ gibt es mehr als 
$80$ Millionen Gitter (vgl.~\cite{serre}, Kap.~7) und daher mindestens ebenso 
viele \VOAs, d.h.~eine Klassifikation f"ur $c\geq 32$ ist praktisch nicht 
durchf"uhrbar.

Die Monster-\VOA $\VM$~(s.~\ref{monstermodul}) ist nach Dong~\cite{Do}
rational und selbstdual.

Schellekens~\cite{schellekens1,schellekens2} hat f"ur $c=24$ eine Lis\-te von 
$71$ selbstdualen 
\VOA-Kan\-di\-da\-ten an\-ge\-ge\-ben. Da\-von sind aller\-dings nur f"ur 
$39$ Kandidaten
($24$ Niemeiergitter-\VOAs \& $15$ "`$\Z_2$-Orbifolds"'~\cite{DGM}) die volle 
\VOA-Struktur konstruiert
und f"ur $25$ davon (Niemeiergitter-\VOAs{}~\cite{Do-gitter} und $\VM$) 
die Selbstdualit"at bewiesen.
Die Konstruktion der Liste beruht auf der Klassifikation der von $V_1$
erzeugten affinen Liealgebren. Da $V_1$ nach Lemma~\ref{lemmareduktiv} reduktiv
ist, falls $V$ \nett und unit"ar ist, fassen wir zusammen:
\begin{vermutung}\label{vermutungsdvoa}
Sei $V$ eine selbstduale unit"are \nette rationale \VOA
vom Rang \hbox{$c\leq24$.} Dann gilt
\begin{eqnarray*}
\mbox{f"ur\ } c=8\  & \!\!:& V\cong V_{E_8}, \\
\mbox{f"ur\ } c=16& \!\!:& V\cong V_{E_8}\otimes V_{E_8},\quad V_{D_{16}^+}, \\
\mbox{f"ur\ } c=24&\!\! :& V \mbox{\ ist isomorph zu einer der $71$ \VOAs 
   aus~\cite{schellekens1}.}
\end{eqnarray*}
\end{vermutung}
Bei Annahme der Eigenschaft \Llie{} (siehe~\ref{L1-lie}) k"onnen wir f"ur $c=8$
bzw.~$c=16$ die Vermutung fast beweisen. Aus Satz~\ref{charsdvoa} folgt 
$\dim V_1=248$ bzw.~$\dim V_1=496$. Die Gleichungen~(s.~\ref{satzlievoas})
$\dim V_1=\sum \dim {\bf g}_i$ und 
$c=\sum\frac{k_i\cdot\dim {\bf g}_i}{k_i+\check{g}_i)}$
liefern als einzige L"osungen \hbox{$\widetilde{V}_1=V_{E_{8,1}}$} f"ur $c=8$ 
bzw.~(a) $\widetilde{V}_1=V_{E_{8,1}}\otimes V_{E_{8,1}}$ oder (b) 
$\widetilde{V}_1=V_{D_{16,1}}$ f"ur $c=16$. 
Da $\chi_V=\chi_{V_{E_{8,1}}}$, folgt $V=V_{E_{8,1}}$ f"ur $c=8$. F"ur $c=16$ 
folgt (a) $V=V_{E_{8,1}}\otimes V_{E_{8,1}}$ oder (b) $V=V_{D_{16,1}}\oplus
V_{D_{16,1}}^{h=2}$. Die Gitter-\VOAs zu den Wurzelsystemen $A_n$, $D_n$,
$E_n$ sind nach Satz~\ref{satzgitterlieiso} isomorph zu den \VOAs,
die zu den fundamentalen Stufe $1$ Darstellungen der zugeh"origen affinen Kac-Moody-Algebren
assoziiert sind. Allerdings bleibt zu zeigen, da"s bei (b) die 
\VOA-Struktur stets mit der \VOA-Struktur der Gitter-\VOA "ubereinstimmt.

\section{Selbstduale \SVOAs{ }}

F"ur selbstduale Codes ist es offensichtlich, da"s sie nur f"ur gerade Dimensionen existieren. Der Rang eines Gitters ist nach Definition
immer ganz. Wir werden hier zeigen, da"s der Rang einer selbstdualen
\SVOA halbganz ist. Als Vorbereitung ben"otigen wir folgendes einfaches
\begin{lemma}\label{nennerlemma}
Sei $f(p)=1+p+a_2\,p^2+\cdots\in\Z[[p]]$. Dann werden die Nenner der 
Potenzreihe $f^r$ f"ur $r\in\Q\setminus \Z$ beliebig gro"s.
\end{lemma}
{\bf Beweis:} Sei $\pi$ ein Primteiler des Nenners des gek"urzten Bruches
$r=\frac{u}{v}$. Annahme:
$$f^r=\sum_{i=0}^{\infty}b_i\,p^i=1+\frac{u}{v}\,p+O\!\left(p^2\right)$$
hat beschr"ankte Nenner. Dann gibt es einen kleinsten Index $i_0\geq 1$,
f"ur den $b_{i_0}$ die maximale $\pi$-Potenz $\pi^l$ mit $l\geq 1$ in
den Nennern der $b_i$ hat. Der Koeffizient von $p^{v\cdot i_0}$ in 
$f^{r \cdot v}=f^u\in \Z[[p]]$ ist
$$\sum_{{{\scriptstyle \alpha_1,\dots,\alpha_v} \atop 
{\scriptstyle \alpha_1+\cdots+\alpha_v=v\cdot i_0}}}
b_{\alpha_1}\cdot\, \cdots\, \cdot b_{\alpha_v}.$$
Der in der Summe enthaltende Summand $b_{i_0}^v$ hat den Faktor $\pi^{v\cdot l}$
im Nenner. Nach Definition von $i_0$ ist dies aber auch der einzige
mit dieser Eigenschaft, im Widerspruch zur Ganzzahligkeit der Summe. \qed

\begin{satz}\label{rangsdsvoa}
Der Rang $c$ einer selbstdualen \netten rationalen \SVOA $V$ ist eine 
halbganze Zahl.
\end{satz}
{\bf Beweis:} Da $V$ selbstdual und rational ist, bildet nach 
Korollar~\ref{trans-charsvoa-korollar} der von dem Charakter $\chi_V$ 
aufgespannte $\C$-Vektorraum eine eindimensionale Darstellung von 
$\GT=\langle S, T^2 \rangle$. Wie beim Beweis von 
Satz~\ref{charsdvoa} hat man:
\begin{eqnarray}\label{transf}
\chi_V(\tau)\mid_{T^2}& =& e^{2 \pi i\frac{2 c}{24}}\cdot\chi_V(\tau), \\
\chi_V(\tau)\mid_S& =& \chi_V(\tau).
\end{eqnarray}

Bezeichne f"ur eine nat"urliche Zahl $N$
mit $\GTN$ den Kern des durch $\rho(S)=1$ und $\rho(T^2)=
e^{2 \pi i\cdot \frac{1}{N}}$ gegebenen Gruppenhomomorphismus
$\rho: \GT\longrightarrow\C^*$. Da $T^2$ in $\GT$ einen freien 
zyklischen Faktor erzeugt,
ist $\rho$ wohldefiniert. Die Fl"ache $\overline{\H/\GTN}$
ist dann eine $N$-bl"attrige "Uberlagerung von $\overline{\H/\GT}\cong
{\bf CP}_1$, die in den beiden Spitzen $\infty$ und $1$ von $\GT$ einen
$N$-fachen Verzweigungspunkt besitzt. 
Die Modulfunktion $j_{\theta}$ (s.~\ref{modulfktjt}) bildet 
$\overline{\H/\GT}$ biholomorph nach ${\bf CP}_1$ ab, d.h.~der
K"orper der meromorphen Funktionen von $\overline{\H/\GT}$ sind
die rationalen Funktionen in $\jt$.
Da $\jt$ einen einfachen Pol in
der Spitze $\infty$ und eine einfache Nullstelle in der Spitze $1$ besitzt,
ist $\C(\!\sqrt[N]{\jt})$ der K"orper der meromorphen Funktionen auf $\overline{\H/\GTN}$.

Sei $N$ der Nenner der rationalen Zahl $\frac{2c}{24}$.
Nach Satz~\ref{trans-charsvoa} ist $\chi_V$ holomorph in $\H$ und 
wegen~(\ref{transf}) ist $\chi_V$ invariant unter $\GTN$, definiert also
eine auf $\H/\GTN$ holomorphe Funktion. In der Spitze $\infty$ hat $\chi_V$ h"ochstens einen Pol. Dies gilt auch f"ur die andere Spitze:
Nach Definition der Rationalit"at ist auch die Unter-\VOA $V_{(0)}$
rational. Da $\chi_V=\chi_{V(0)}+\chi_{V(1)}$, erh"alt man mit  
Satz~\ref{trans-charvoa} f"ur die $q$-Entwicklung in der \hbox{Spitze $1$}
\begin{equation}\label{STent}
\chi_V\mid_{TS}=\sum_{i=1}^l r_i\, \chi_{M_i},
\end{equation}
wobei die $M_1$ bis $M_l$ die
irreduziblen Moduln von $V_{(0)}$ sind. Die $M_i$ haben nach 
Definition~\ref{SVOA-Modul} nach unten beschr"ankte $L_0$-Eigenwerte, 
d.h.~$\chi_{M_i}$ hat h"ochstens einen Pol in $\infty$. 
Insgesamt ist somit $\chi_V$ eine meromorphe Funktion auf $\overline{\H/\GTN}$, 
also $\chi_V\in{\bf C}(\!\sqrt[N]{\jt})$.

Die $24$-te Wurzel $\chi_{1/2}$ aus $\jt$ hat nach (\ref{modulfktjt})
noch eine ganzzahlige
$q$-Entwicklung:\footnote{Lemma~\ref{nennerlemma} zeigt, da"s $\sqrt[N]{\jt}$
genau dann eine $q$-Entwicklung mit beschr"ankten Nennern hat, falls $N$ ein
Teiler von $24$ ist. Da die Koeffizienten von Modulfunktionen zu $\Gamma(n)$
mit rationaler $q$-Entwicklung stets beschr"ankte Nenner besitzen, kann
nur f"ur $N\mid 24$ die Gruppe $\GTN$ eine Kongruenzuntergruppe sein.
Andererseits gilt $\Gamma(48)\subset\Gamma_{\theta,24}$, wie man mit
(\ref{modulfktjt}) und dem Transformationsverhalten der $\eta$-Funktion unter 
beliebigen Matrizen $A\in\slz$ (vgl.~\cite{Ra}, S.~163) leicht nachrechnet, 
d.h.~$\GTN$ ist f"ur $N\mid 24$ eine Kongruenzuntergruppe.}
  $$\chi_{1/2}=\sqrt[24]{\jt}=q^{-\frac{1}{48}}\,\prod_{n=0}^{\infty}
(1+q^{n+\frac{1}{2}})\in q^{-\frac{1}{48}}\Z[[q^{1/2}]]. $$
Sei $c=\bar{c}+r$ mit $\bar{c}\in \frac{1}{2}\Z$ und $0\leq r<\frac{1}{2}$.
Die Funktion $\chi_V\cdot\chi_{1/2}^{-2r}\in q^{-\bar c/24}\Q[[q^{1/2}]]$
multipliziert sich mit einer $24$-ten Einheitwurzel unter $T^2$
und liegt daher in ${\cal M}(\overline{\H/\Gamma_{\theta,24}})=\C(\chi_{1/2})$.
Es folgt
\begin{equation}\label{partent}
\chi_V\cdot\chi_{1/2}^{-2r}=\sum_{i=0}^{n_0}a_i\cdot\chi_{1/2}^{2\bar{c}-24i}
\end{equation}
mit rationa\-len Zahlen $a_i$. Sei $M$ der gemeinsame
Nen\-ner der $a_i$.
Da auch $\chi_V^{-1}\in q^{c/24}\Z[[q^{1/2}]]$, erh"alt man schlie"slich
$$ \chi_{1/2}^{-2r}\in q^{\frac{r}{24}}\frac{1}{M}\Z[[q^{1/2}]]. $$
Nach Lemma~\ref{nennerlemma} ist dies aber nur f"ur $r=0$ m"oglich, d.h.~es
gilt $c=\bar{c}\in\frac{1}{2}\Z$.\qed

Satz~\ref{rangsdsvoa} war von Verlinde in~\cite{verlinde}, S. 375 vermutet 
worden.
Setzen wir $V$ zus"atzlich noch als unit"ar voraus, so k"onnen wir den 
Charakter noch genauer bestimmen:

\begin{satz}\label{sdsvoachar}
Sei $V$ eine selbstduale \nette unit"are rationale \SVOA vom Rang $c$.
Dann ist ihr Charakter $\chi_V=q^{-\frac{c}{24}}\sum_{n=0}^{\infty}
\dim V_n\, q^n$ ein homogenes Polynom $P(x,y)$ vom Gewicht $c$ in den Reihen
$x=\sqrt{\frac{\Theta_{\Z}}{\eta}}=\chi_{1/2}$ (Gewicht $1/2$) 
und $y=\sqrt[3]{j}$ (Gewicht $8$), d.h.~in den Charakteren der 
\SVOAs $\VF$ vom Rang $1/2$ und $V_{E_8}$ von Rang $8$.
\end{satz}

{\bf Beweis:}
Wie im letzten Satz gezeigt, ist $c$ halbganz, und nach 
Gleichung~(\ref{partent}) gilt daher
\begin{equation}\label{nextent}
\chi_V=\sum_{i=0}^{n_0}a_i\,\chi_{1/2}^{2c-24i}.
\end{equation} 
Auf $\overline{\H/\Gamma_{\theta,24}}$ hat $\chi_{1/2}$ einen einfachen Pol  
in $\infty$ (lokale Koordinate $q^{1/48}$) und eine einfache Nullstelle in
$1$ (lokale Koordinate $q^{1/24}$).
Nach Gleichung~(\ref{STent}) gilt f"ur die Entwicklung von
$\chi_V$ \hbox{in $1$:}
$$\chi_{V}\mid_{TS}=\sum_{i=1}^l r_i\, \chi_{M_i} \qquad\mbox{mit\quad} \chi_{M_i}=
q^{-\frac{c}{24}+h_i}\sum_{n \geq 0}\dim M_i\, q^n.$$
Da $V$ als unit"ar vorausgesetzt wurde, ist nach~ Lemma~\ref{lemmaunitaer}
$h_i\geq 0$. 
Daher hat $\chi_{M_i}$ h"ochstens einen Pol der Ordnung
$c$ in $\infty$ bzw.~$\chi_V$ h"ochstens einen Pol der Ordnung
$c$ \hbox{in $1$},
und daher kann in~(\ref{nextent}) $n_0=\left[\frac{c}{8}\right]$
gesetzt werden.
Umformen liefert
$$\chi_V=\sum_{i=0}^{\left[\frac{c}{8}\right]}b_i\,\chi_{1/2}^{2c-16i}\cdot
 \left(\chi_{1/2}^{16}-16\,\chi_{1/2}^{-8}\right)^i $$
mit Zahlen $b_i$, die sich linear aus den $a_i$ ergeben.
Also ist $\chi_V$ --- wie behauptet --- ein homogenes Polynom vom Gewicht
$c$ in $\chi_{1/2}=\chi_{\VF}$ und
$(\chi_{1/2}^{16}-16\,\chi_{1/2}^{-8})=\sqrt[3]{j}=\chi_{V_{E_8}}$.
\phantom{xxxxxxxxxxxxxxxx}\qed

"Uber die Eigenschaften der Unter-\VOA $\V{0}$ einer selbstdualen
\netten rationalen \SVOA $V=\V{0}+\V{1}$ lassen sich weitere Aussagen machen.

Nach Satz~\ref{trans-charsvoa} ist der von dem Charakter 
$\chi_V=\chi_{\V{0}}+\chi_{\V{1}}$ aufgespannte \C-Vektorraum invariant unter
$\Gamma_{\theta}=\langle S, T^2 \rangle$. Die drei Funktionen
$\chi_V$, $\chi_V\mid_T$ und $\chi_V \mid_{TS}$ spannen daher einen
h"ochstens $3$-dimensionalen (gleich $3$ falls $\V{0}\not= 0$)
${\rm PSL}_2(\Z)$-Modul auf.
Wir w"ahlen eine neue Basis, auf der $T$ f"ur beliebige $c$ diagonal operiert:
$$\chi_{V(0)}=\frac{1}{2}(\chi_V+e^{2 \pi i \frac{c}{24}}\chi_V\mid_T),\quad
\chi_{V(1)}=\frac{1}{2}(\chi_V-e^{2 \pi i \frac{c}{24}}\chi_V\mid_T),\quad
\chi_{\rm Rest}=\frac{1}{\sqrt{2}}\cdot e^{2 \pi i \frac{c}{24}}\chi_V
\mid_{TS}.$$
F"ur die Operation der Erzeuger $S$ und $T$ von ${\rm PSL}_2(\Z)$ erhalten wir
unter Verwendung der Relation $(ST)^3={\rm id}$ in dieser Basis die Matrizen
\begin{equation}\label{ST-svoatrans}
\overline{T}=\left(\begin{array}{ccc}
e^{- 2 \pi i \frac{c}{24}} & 0 & 0 \\
0 &  e^{2 \pi i (-\frac{c}{24}+\frac{1}{2})} & 0 \\
0 & 0 & e^{2 \pi i (-\frac{c}{24}+\frac{c}{8})}
\end{array}\right),\qquad
\overline{S}=\left(\begin{array}{ccc}
\frac{1}{2} & \frac{1}{2} & \frac{1}{\sqrt{2}} \\
\frac{1}{2} & \frac{1}{2} & -\frac{1}{\sqrt{2}} \\
\frac{1}{\sqrt{2}} & -\frac{1}{\sqrt{2}} & 0
\end{array}\right).
\end{equation}
Die Funktion $\chi_{\rm Rest}$ wurde so normiert, da"s $\overline{S}$ unit"ar
ist. 
Nehmen wir nun an, da"s $\V{0}$ und $\V{1}$ irreduzible  $\V{0}$-Moduln sind,
die Verlindeformel~(\ref{verlindeformel}) gilt und die "`Quantendimensionen"'
 $\frac{\widetilde{S}_{i1}}{\widetilde{S}_{11}}\geq 1$ sind (s.~\cite{Fuchs}, S.~343),
so sind die m"oglichen $\widetilde{S}$-Matrizen (s.~\ref{trans-charvoa})
gegeben durch
\begin{equation}\label{smatrixsehrnett}
\hbox{Fall (a)\quad}\widetilde{S}=\left(\begin{array}{rrr}
\frac{1}{2} & \frac{1}{2} & \frac{1}{\sqrt{2}} \\
\frac{1}{2} & \frac{1}{2} & -\frac{1}{\sqrt{2}} \\
\frac{1}{\sqrt{2}} & -\frac{1}{\sqrt{2}} & 0
\end{array}\right),\qquad 
\hbox{Fall (b)\quad}  \widetilde{S}=\left(\begin{array}{rrrr}
\frac{1}{2} & \frac{1}{2} & \frac{1}{2} & \frac{1}{2} \\
\frac{1}{2} & \frac{1}{2} &-\frac{1}{2} &-\frac{1}{2} \\
\frac{1}{2} &-\frac{1}{2} &\mp\frac{i}{2} &\pm \frac{i}{2} \\
\frac{1}{2} &-\frac{1}{2} &\pm\frac{i}{2} &\mp\frac{i}{2} \\
\end{array}\right)
\end{equation}
$$
\hbox{ oder Fall (c)\quad}  \widetilde{S}=\left(\begin{array}{rrrr}
\frac{1}{2} & \frac{1}{2} & \frac{1}{2} & \frac{1}{2} \\
\frac{1}{2} & \frac{1}{2} &-\frac{1}{2} &-\frac{1}{2} \\
\frac{1}{2} &-\frac{1}{2} &\pm\frac{1}{2}&\mp\frac{1}{2} \\
\frac{1}{2} &-\frac{1}{2} &\mp\frac{1}{2}  &\pm\frac{1}{2} \\
\end{array}\right),
$$
und $\chi_{\rm Rest}$ ist der Charakter von einem 
bzw.~$\sqrt{2}\chi_{\rm Rest}$ der Charakter von zwei $\V{0}$-Moduln.
F"ur die zugeh"origen Fusionsalgebren ${\cal F}(\V{0})$ erg"abe sich aus
der Verlindeformel im
\begin{equation}\label{fusionsehrnett}
{\vbox{Fall (a)\ : \  3 Moduln $\V{0}$,  $\V{1}$ und  $\V{2}$,
Fusionsalgebra wie beim "`Isingmodell"', \newline
Fall (b)\ : \ 4 Moduln $\V{0}$,  $\V{1}$, $\V{2}$ und $\V{3}$,
Fusionsalgebra ${\cal F}(\V{0})\cong\Z[\Z/4\Z]$,  \newline
Fall (c)\  : \  4 Moduln $\V{0}$,  $\V{1}$, $\V{2}$ und  
$\V{3}$,
Fusionsalgebra ${\cal F}(\V{0})\cong\Z[\Z/2\Z\times\Z/2\Z]$.} }
\end{equation}
Wir machen daher die folgende 
\begin{definition}\label{sehrnett}
Eine selbstduale \nette rationale \SVOA $V$ hei"st \sehrnett, falls die
Fusionsalgebra der Unter-\VOA $\V{0}$ eine der drei F"alle (a), (b) oder (c)
in~(\ref{fusionsehrnett}) ist, die Operation von $S\in {\rm SL}_2(\Z)$ 
auf den "`Korellationsfunktionen auf dem Torus"' durch~(\ref{smatrixsehrnett}) 
gegeben ist und $\chi_{\rm Rest}$ der Charakter des Moduls $V_{(2)}$ 
(im Fall (a)) bzw.~$\sqrt{2}\chi_{\rm Rest}$ die Summe der
Charaktere von $V_{(2)}$ und $V_{(3)}$ (F"alle (b) und (c)) ist. Au"serdem
ist der Fall $V=\V{0}$ zugelassen, d.h.~$V$ ist sogar eine selbstduale \VOA.
\end{definition}
Vermutlich ist eine selbstduale \nette unit"are rationale \SVOA stets 
\sehrnett, hierzu m"u"ste man die Beziehung zwischen den Algebren $A(\V{0})$ 
und $A(V)$ genauer untersuchen (vgl.~die Beweisskizze von 
Satz~\ref{trans-charsvoa}).

Es besteht der folgende Zusammenhang zwischen dem nach Satz~\ref{rangsdsvoa} halbganzen Rang $c$ und der Struktur der Fusionsalgebra von $\V{0}$:
\begin{satz}\label{sdsvoafusion}
Sei $V$ eine selbstduale \sehrnette unit"are rationale \SVOA. In Abh"angigkeit
vom Rang $c$ ist die Struktur der Fusionsalgebra von $\V{0}$ gegeben durch
\begin{equation}\label{svoafusion}
{\cal F}(\V{0})=\cases{
\hbox{Ising (Fall (a)),}              & falls $c\in\frac{1}{2}\Z\setminus\Z$,\cr
\hbox{$\Z[\Z/4\Z]$ (Fall (b)),}       & falls $c\in\Z\setminus 2\Z$,\cr
\hbox{$\Z[\Z/2\Z\times\Z/2\Z]$ (Fall (c)),}& falls $c\in 2\Z$.\cr}
\end{equation}
\end{satz}
{\bf Beweis:} Nach Gleichung~(\ref{nextent}) gilt 
$\chi_V=\sum_{i=0}^{n_0}a_i\chi_{1/2}^{2c-24i}$ mit ganzen Zahlen $a_i$,
also 
$$\chi_{\rm Rest}=\frac{1}{\sqrt{2}}\,e^{2 \pi i \frac{c}{24}}\,\chi_V\mid_{TS}
=\frac{1}{\sqrt{2}}\sum_{i=0}^{n_0}a_i\,(-1)^{i}\,\left(e^{\frac{2 \pi i}{48}}
\chi_{1/2}\mid_{TS}\right)^{2c-24i}. $$
Da $\chi_{1/2}\mid_{TS}=e^{-\frac{2 \pi i}{48}}\sqrt{2}\,q^{\frac{1}{24}}
\prod_{n=1}^{\infty}
(1+q^n)$, hat $\chi_{\rm Rest}$ f"ur $c\in\frac{1}{2}\Z\setminus\Z$
eine rationale \hbox{$q$-Entwicklung} bzw.~$\sqrt{2}\chi_{\rm Rest}$ f"ur
$c\in\Z$, d.h.~f"ur $c\in\frac{1}{2}\Z\setminus\Z$ mu"s Fall (a) vorliegen
und f"ur $c\in\Z$ der Fall (b) oder (c).

Im Fall (b) oder (c) gilt $\widetilde{T}={\rm Diag}(
e^{- 2 \pi i \frac{c}{24}}, e^{2 \pi i (-\frac{c}{24}+\frac{1}{2})},
e^{2 \pi i (-\frac{c}{24}+\frac{c}{8})},e^{2 \pi i (-\frac{c}{24}+\frac{c}{8})}
)$. Man erh"alt aus $(ST)^6={\rm id}$ in $\slz$
im Fall (b) die Bedingung  $c\in\Z\setminus 2\Z$
bzw.~im Fall (c) die Bedingung $c\in 2 \Z$ . 
Die Vorzeichen in $\widetilde{S}$ sind durch die Restklasse $c$ mod $4$ 
eindeutig festgelegt.\qed

\label{twistsektor}
Definiert man auf einer \SVOA $V$ einen linearen Isomorphismus $\kappa$ durch
$$ \kappa(v)=\cases{ v, &f"ur $v\in V_{(0)}$, \cr -v, &f"ur $v\in V_{(1)}$,} $$
wobei $V_{(0)}$ und $V_{(1)}$ der gerade bzw.~ungerade Unterraum von $V$ sind, 
so ist $\kappa$ eine Involution und ein Automorphismus der \SVOA $V$.
Es ist daher nat"urlich, $\kappa$-getwistete \hbox{$V$-Moduln}
 zu betrachten (f"ur
eine Definition siehe~\cite{ffr}). In der physikalischen Literatur werden die
(ungetwisteten) $V$-Moduln als die {\it Neveu-Schwarz Sektoren f"ur $V$\/}
und die \hbox{$\kappa$-getwisteten} 
$V$-Moduln als die {\it Ramond Sektoren f"ur $V$\/}
bezeichnet. F"ur eine selbstduale \sehrnette \SVOA $V$ ist zu erwarten, da"s
die irreduziblen $\kappa$-getwisteten Moduln von $V$ gerade $V_{(2)}\oplus V_{(2)}$
($c\in \Z+\frac{1}{2}$) bzw.~$V_{(2)}\oplus V_{(3)}$ ($c \in \Z$) sind.

\pagebreak[2]
{\bf Beispiele von selbstdualen \sehrnetten unit"aren rationalen \SVOAs:}
\nopagebreak[2]

{\it Die \SVOA $\VF$:\/} Eine \SVOA mit dem kleinsten zul"assigen Rang
$c=\frac{1}{2}$ wurde in \cite{KacWang}, konstruiert.
Diese auch als "`freies Fermion"' bezeichnete \SVOA $\VF$ ist rational und
selbstdual (Satz 4.1 in \cite{KacWang}). F"ur die Eigenschaft unit"ar
vgl.~\cite{Go-mero}, Abschnitt 8. Die Unter-\VOA $(\VF)_{(0)}$ kann nach
Satz~\ref{L2-vir} nur die zur unit"aren \hbox{$c=\frac{1}{2}$} 
H"ochst\-gewichts\-darstellung
der Virasoro\-algebra assozierte \VOA $\La$ sein, deren Rationalit"at
in \cite{Wan} gezeigt wurde. Dort wurden auch die $3$ irreduziblen 
\hbox{$\La$-Moduln}
$\La$,  $\Lb$ und  $\Lc$ mit den konformen Gewichten $0$, $\frac{1}{2}$ und 
$\frac{1}{16}$ bestimmt.
F"ur die Fusionsalgebra von  $\La$ ergibt sich nach \cite{DoMaZhu},
Abschnitt 3, die in der physikalischen Literatur als "`Fusionsregeln des
Isingmodells"' bezeichnete Struktur:
$$  \Lb \times  \Lb =  \La, \qquad
 \Lb \times  \Lc = \Lc, $$
\begin{equation}\label{fusion-ising} \Lc \times  \Lc= \La +  \Lb \end{equation}
und $\La$ ist die Identit"at in der Fusionsalgebra.
Da f"ur den Charakter $\chi_{\VF}=\chi_{1/2}$ gilt
(Satz~\ref{sdsvoachar} oder (\ref{char-fermi})),
mu"s schlie"slich $\VF=
 L_{1/2}(0)\oplus  L_{1/2}(\frac{1}{2})$ gelten.
Nach Gleichung~(\ref{char-vir}) ist weiter $\frac{1}{\sqrt{2}}\chi_{\rm Rest}=\chi_{\Lc}$, d.h.~$\VF$ ist \sehrnett.

{\it Die zu ungeraden selbstdualen Gittern $L$ assoziierten \SVOAs $V_L$:\/}
Die Eigenschaft \sehrnett ergibt sich aus Satz~\ref{satzgittersvoa}, da f"ur 
das gerade Untergitter $L_0 \subset L=L_0 \cup (L_0+[2])$ in geraden 
Dimensionen $L_0^*/L_0 \cong \Z/2\Z \times \Z/2\Z$ und in ungeraden 
Dimensionen $L_0^*/L_0\cong \Z
/4\Z$ gilt und die Zerlegung der \SVOA $V_L=V_{L_0} \oplus V_{L_0 + [2]}$ 
in den geraden und den ungeraden Teil besteht.

{\it Die SVOA $V\!O^{\,\natural}=W_{NS}$:\/} Die mit der Monster \VOA $\VM$ 
in Beziehung stehende \SVOA
$W_{NS}$ ist das \VOA-Analogon zu dem ungeraden Golay Code $Z_{24}$~\cite{PlSl}
und dem ungeraden Leechgitter $O_{24}$~\cite{CoPa} in Dimension $24$ und
k"onnte daher auch als {\it ungerader Mondscheinmodul} $V\!O^{\,\natural}$
bezeichnet werden.
Die \SVOA $W_{NS}$ ist eine Unter-\SVOA einer in \cite{DGH-Monster} 
beschriebenen "`superkonformen"' \VOA $W$ vom Rang $24$:
In \cite{Hua} wird gezeigt, da"s die bei der Konstruktion des Monstermoduls 
$\VM=V_{\Lambda}^+\oplus(V_{\Lambda}^T)^+$ verwendete Unter-\VOA 
$V_{\Lambda}^+$ die $4$ irreduziblen Moduln $V_{\Lambda}^+$, $V_{\Lambda}^-$,
$(V_{\Lambda}^T)^-$ und $(V_{\Lambda}^T)^+$ mit den konformen Gewichten
$0$, $1$, $\frac{3}{2}$ und $2$ besitzt, und auf $W=V_{\Lambda}^+\oplus V_{\Lambda}^-\oplus (V_{\Lambda}^T)^- \oplus (V_{\Lambda}^T)^+ $
wird die Existenz einer "`superkonformen"' Struktur bewiesen.
Diese induziert dann auf $W_{NS}=V_{\Lambda}^+ \oplus (V_{\Lambda}^T)^-$ 
eine \SVOA-Struktur (die sogar $N=1$ supersymmetrisch im Sinne
von \cite{KacWang} ist). Die Unter-\VOA $(W_{NS})_{(0)}$ von $W_{NS}$
ist $V_{\Lambda}^+$ und hat die Fusionsalgebra $\Z[\Z/2\Z\times\Z/2\Z]$
(\cite{Hua}, Satz 3.7). Die Selbstdualit"at von $W_{NS}$ wurde nicht gezeigt,
ist aber zu vermuten. In Analogie zu der Situation bei Codes und Gittern ist
au"serdem zu vermuten, da"s $\VM$ und $W_{NS}$ die einzigen selbstdualen
\sehrnetten unit"aren rationalen \SVOAs mit $V_1=0$ sind. Das Beispiel
$W_{NS}$ zeigt, da"s die Charaktere von $V_{(2)}$ und $V_{(3)}$, deren Summe
die Beziehung $\sqrt{2}\chi_{\rm Rest}=\chi_{V_{(2)}}+\chi_{V_{(3)}}$ erf"ullt,
auch verschieden sein k"onnen.

Das Tensorprodukt zweier selbstdualer \sehrnetter \SVOAs sollte wieder
\sehrnett sein. Wir betrachten hier das im n"achsten Kapitel ben"otigte
{ \it Beispiel $\VF^{\otimes k}$\/}:\hfill\break
Zuerst sei daran erinnert, da"s die \VOA $V_{D_{l,1}}$, die zur fundamentalen
Stufe $1$ Darstellung der affinen Kac-Moody Algebra $\widetilde{D}_l\cong \widetilde{\bf so}(2l)$ assoziert ist,
 isomorph zur Gitter-\VOA $V_{D_l}$ von Rang $l$ ist
(Satz~\ref{satzgitterlieiso}), und daher nach Satz~\ref{satzgittersvoa}
die Fusionsalgebra ${\cal F}(V_{D_l})$ durch (\ref{svoafusion}) gegeben ist.
Die Charaktere der irreduziblen Moduln sind in 
Abschnitt~\ref{beispielcharaktere} beschrieben. Die zu $\widetilde{B}_l
\cong \widetilde{\bf so}(2l+1)$, $l\geq 1$ geh"orige \VOA $V_{B_l,1}$
hat als irreduzible Moduln die $3$ Stufe $1$
Darstellungen von $\widetilde{B}_l$. Die Fusionsalgebra wird f"ur affine
Kac-Moody \VOAs ganz allgemein in~\cite{FreZhu} beschrieben, und man pr"uft
nach, da"s sich die Ising Fusionsregeln ergeben.
Die Charaktere ergeben sich aus der Weyl-Kac-Charakterformel~\cite{Kac}. 
\begin{satz}\label{vsok}
Das $k$-fache Tensorprodukt $\VF^{\otimes k}$ ist eine selbstduale
\sehrnette unit"are rationale \SVOA von Rang $\frac{k}{2}$.
Die Unter-\VOA $V_{{\bf so}(k)}:=(\VF^{\otimes k})_{(0)}$ ist 
die zu den folgenden H"ochstgewichtsdarstellungen von 
affinen Kac-Moody Algebren assoziierte \hbox{\VOA:}
$$V_{{\bf SO}(k)}=\cases{V_{D_{l,1}} & falls $l=2k$, \cr 
V_{B_{l,1}} & falls $l=2k+1$,}$$
und der ungerade Teil $(\VF^{\otimes k})_{(1)}$ ist der 
(bzw.~einer, falls $k=8$) Stufe $1$ Modul mit dem konformen Gewicht
$h=\frac{1}{2}$.
\end{satz}
{\bf Beweis:} Nach Satz~\ref{tensorprodukt} ist $\VF^{\otimes k}$ eine
\nette unit"are rationale \SVOA. Zu zeigen bleibt die angegebene Struktur
des geraden und ungeraden Teils von $\VF^{\otimes k}$. Dann gelten wegen
der oben beschriebenen Eigenschaften von $V_{D_{l,1}}$ und $V_{B_{l,1}}$
alle die f"ur \sehrnett ben"otigten Eigenschaften.

Der Summe $V_{D_{l,1}}\oplus V_{D_{l,1}^{h=1/2}}$ kann die Struktur 
einer \SVOA gegeben werden (vgl.~\cite{tsukada,ffr,DoMa-orbi}), f"ur die
Summe $V_{B_{l,1}}\oplus V_{B_{l,1}^{h=1/2}}$ wurde dies 
in~\cite{DML-simcurr} (Beispiel 5.12)
gezeigt. Aufgrund der Eigenschaft \Lcliff{} (Satz~\ref{Leinhalb-cliff}) 
erzeugen 
die Vertexoperatoren des Gewicht $\frac{1}{2}$-Anteils von $\VF^{\otimes k}$
als auch von $V_{D_{l,1}}\oplus V_{D_{l,1}^{h=1/2}}$ 
bzw.~$V_{B_{l,1}}\oplus V_{B_{l,1}^{h=1/2}}$ eine zur unendlich
dimensionalen Cliffordalgebra ${\bf Cliff}(\Z+\frac{1}{2})$ assoziierte 
Unter-\SVOA mit dem Charakter $\chi_{1/2}^k$. Dies ist aber auch der Charakter
von  $V_{D_{l,1}}\oplus V_{D_{l,1}^{h=1/2}}$ 
bzw.~$V_{B_{l,1}}\oplus V_{B_{l,1}^{h=1/2}}$, d.h.~beide \SVOAs sind 
isomorph zur Clifford \SVOA.        \qed

Da nach Satz~\ref{sdsvoachar} der Charakter einer selbstdualen \sehrnetten 
unit"aren und rationalen \SVOA f"ur R"ange $c$ kleiner als $8$ durch
$\chi_{1/2}^{2c}$ gegeben ist, also ${\rm dim}\, V_{1/2}= 2c $ gilt, liefert 
voranstehender Beweis das folgende Klassifikationsresultat f"ur den Bereich 
$0\leq c<8$:

\begin{satz}\label{SVOAS1-8}
Eine selbstduale \sehrnette unit"are rationale \SVOA vom Rang $c$ mit 
\hbox{$0\leq c<8$} ist isomorph zu $\VF^{\otimes 2c}$. \qed
\end{satz}
Im n"achsten Kapitel, in dem ganz allgemein die Klassifikation der selbstdualen
\SVOAs auf diejenige der selbstdualen \VOAs zur"uckgef"uhrt wird, 
werden wir einen Klassifikationssatz f"ur den Bereich $8\leq c<16$
angeben. Der folgende Satz zeigt, da"s man sich auf \SVOAs mit $V_{1/2}=0$
beschr"anken kann.
\begin{satz}[vgl.~\cite{Go-mero}, Abschnitt 8]\label{satzabspalt}
Sei $V$ eine selbst\-duale \nette uni\-t"a\-re ra\-tiona\-le \SVOA 
vom Rang $c$ mit
$\dim V_{1/2} = k$. Dann ist $V$ isomorph zu $W \otimes \VF^{\otimes k}$, 
wobei $W$ selbst wieder eine selbstduale \nette unit"are rationale \SVOA 
vom Rang $c'=c-\frac{k}{2}$ ist.
\end{satz}
{\bf Beweis (Skizze):} 
Da wegen $\chi_{U\otimes V}=\chi_U\otimes\chi_V$, d.h.~$\dim (U\otimes V)_{1/2}=
\dim U_{1/2} + \dim V_{1/2}$, induktiv vorgegangen werden kann, gen"ugt es,
$V=\VF\otimes W$ zu zeigen. Dies folgt aus der Eigenschaft \Lcliff{}
mit $W={\rm Com}_V(\VF)$.
Nach~\cite[Prop.~2.7]{DoMaZhu} ist $W$ rational und 
nach~\cite[Lemma~2.8]{DoMaZhu} ist $W$ selbstdual. (Bei Verallgemeinerung der
jeweiligen S"atze dort auf \SVOAs, vgl.~\cite{KacWang} sowie~\cite{DoLe}.)
Die Eigenschaften \nett und unit"ar f"ur $W$ folgen aus denen von $V$ und $\La$.
\qed
%

\chapter{Die Beziehung zwischen selbstdualen \VOAs und \SVOAs{ } }

Die Klassifikation der ungeraden selbstdualen Codes bzw.~Gitter l"a"st sich
auf die Klassifikation der geraden selbstdualen Codes bzw.~Gitter in der 
n"achsth"oheren zul"assigen Dimension zur"uckf"uhren.
In diesem Kapitel wird ein analoges (Teil-)Resultat f"ur die Beziehung
zwischen selbstdualen unit"aren \sehrnetten rationalen \SVOAs und \VOAs 
bewiesen.

\medskip

Wir skizzieren hier kurz die Beziehung zwischen ungeraden und geraden
selbstdualen Codes (siehe \cite{CoPless}, S.~42).
Sei $C$ ein selbstdualer gerader Code der Dimension $n\in 8\Z$. Bezeichne
mit $d_k$ den Tetradencode der L"ange $k\in 2 \Z$, 
dies ist der gerade Untercode von
$c_2^{k/2}$, wobei $c_2=\{(0,0),(1,1)\}$. Der duale Code $d_k^{\perp}$
hat bez"uglich $d_{k}$ vier Nebenklassen $d_k^0=d_k$, $d_k^1$, $d_k^2$ 
und $d_k^3$, die jeweils Vektoren des Minimalgewichts $0$, $2$, $\frac{k}{2}$ 
und $\frac{k}{2}$ enthalten. 
Bei Wahl eines Untercodes $d_k \subset C$ zerlegt sich $C$ in
$$(w^0\oplus d_k^0)\cup (w^1\oplus d_k^1)\cup(w^2\oplus d_k^2)
\cup(w^3\oplus d_k^3).$$
Hierbei ist $w^0$ das Komplement von $d_k$ in $C$, d.h.~die Vektoren aus $C$
der Gestalt $(u_1,u_2,\ldots,u_{n-k},0,\ldots,0)$, wobei f"ur $d_k$ die
letzten $k$ Koordinaten verwendet wurden.
Die Vereinigung $W:=w^0\cup w^1$ bildet dann (falls $w^1\not =0$)
einen ungeraden selbstdualen Code der 
L"ange $n-k$. Es leicht einzusehen, da"s die Konstruktion umkehrbar
ist und eine Bijektion zwischen den Isomorphieklassen von ungeraden
selbstdualen Codes der L"ange $n-k$ und den Isomorphieklassen von
Paaren $(C,d_k)$ von selbstdualen Codes $C$ mit Untercode $d_k$ liefert.
(F"ur $k=4$ mu"s zus"atzlich noch die Wahl der Nebenklasse $d_k^1$ fixiert
werden, da das Minimalgewicht aller drei Nebenklassen $d_k^1$, $d_k^2$ und
$d_k^3$ gleich $2$ ist.
Umgekehrt ist eine Unterscheidung zwischen $d_k^2$ und $d_k^3$ nicht 
erforderlich,
da beide durch einen "`"au"seren"' Automorphismus von $d_k$ vertauscht werden
k"onnen.)

F"ur die analoge Konstruktion bei selbstdualen Gitter siehe die 
Arbeit~\cite{CoSl-ul23}. Die Rolle von $d_k$ "ubernimmt hier das Gitter $D_k$,
das gerade Untergitter von $\Z^k$.

Im Fall der Vertexoperatoralgebren schlie"slich ist $\VSO{k}$, die Unter-\VOA
der \SVOA $\VF^{\otimes k}$, das richtige Analogon. In Abschnitt~\ref{hinR}
beweisen wir einen Satz f"ur die Richtung \VOA nach \SVOA, in 
Abschnitt~\ref{rueckR} diskutieren wir die R"uckrichtung.

\section{Die Konstruktion von \SVOAs aus \VOAs{ }}\label{hinR}

In diesem Abschnitt konstruieren wir zu einer selbstdualen VOA $V$ und
einer Unteralgebra $\VSO{k}$ eine \SVOA $W$. Da beim gegenw"artigen Stand der
Theorie der \VOAs es im allgemeinen nicht bekannt ist, ob die Kommutante
von $\VSO{k}$ in $V$ rational ist --- was aber zu vermuten ist ---, machen
wir hier die etwas technische Voraussetzung, da"s die Kommutante zumindest eine
{\it rationale Unter-\VOA{}\/} mit gleichem Virasoroelement 
wie die Kommutante enth"alt.
Die Konstruktion von $W$ h"angt aber nicht von der Wahl dieser Unter-\VOA
ab. In den uns
interessierenden F"allen ${\rm rank}(V)\leq 24$ ist diese Voraussetzung auch
erf"ullt; man findet stets eine geignete Virasoro oder Kac-Moody Unter-\VOA.

Sei also $V$ eine \nette rationale \VOA vom Rang $c$ 
mit Virasoroelement $\omega$ zusammen mit einer
Unter-\VOA $L:=\VSO{k}$ vom Rang $\frac{k}{2}$ mit Virasoroelement $\omega''$.
Die Kommutante $\overline{W}_{(0)}=\COM_V(L):=\{v\in V\mid \omega_0'' v=0\}$
ist nach~\cite{FreZhu}, Satz 5.1 und 5.2 eine Unter-\VOA von $V$ vom Rang $c'=c-\frac{k}{2}$, wenn die Bedingung $\omega_2\omega''=0$ erf"ullt ist.
Weiter sei angenommen, da"s eine 
rationale Unter-\VOA $U\subset\overline{W}_{(0)}$ mit gleichem Virasoroelement
$\omega':=\omega-\omega''$ wie $\overline{W}_{(0)}$ existiert.
Nach~\cite{DoMaZhu}, Prop.~2.7 ist dann auch $U\otimes L$ rational, und wir 
haben eine direkte Summenzerlegung
\begin{equation}\label{dirsumme}
V=\bigoplus_{i \in I}M_{i}\otimes N_{i}
\end{equation}
mit irreduziblen $U$ bzw.~$L$ Moduln $M_{i}$ und $N_{i}$.
Nach Satz~\ref{vsok} besitzt $L$, in Abh"angigkeit vom Rang,
$3$ oder $4$ Typen von irreduziblen Moduln, die wir mit
$L=L_{(0)}$, $L_{(1)}$, $L_{(2)}$ (und $L_{(3)}$) bezeichnen und die die
konformen Gewichte $0$, $\frac{1}{2}$, $\frac{k}{8}$ (und $\frac{k}{8}$)
besitzen. Setzen wir 
$W_{(a)}:=\bigoplus_{i \in I_a}M_{i}$,
mit $I_a=\{i\in I\mid N_i\cong L_{(a)}\}$,
so ist (\ref{dirsumme}) "aquivalent zu 
\begin{equation}\label{dirtensorsumme}
V=\bigoplus_{a=0,1,2(,3)}W_{(a)}\otimes L_{(a)}.
\end{equation} 
Die Graduierung von $W_{(a)}$ und $L_{(a)}$ ist durch die Virasoroelemente
$\omega'\in \overline{W}_{(0)}$ bzw.~$\omega''\in L_{(0)}$ gegeben und sie 
ist kompatibel mit der von $V$.
\begin{lemma}\label{identifizierung}
Der $U$-Modul $W_{(0)}\otimes {\bf 1}$ stimmt mit der Unter-\VOA 
$\overline{W}_{(0)}$ "uberein.
\end{lemma}
{\bf Beweis:} Ein von Null verschiedener Vektor $m\otimes n\in M_i\otimes N_i$
liegt genau dann in der Kommutante $\overline{W}_{(0)}$ von $L$, wenn f"ur 
das Virasoroelement $\omega''={\bf 1}\otimes\omega''$ die Gleichung
$({\bf 1}\otimes \omega'')_0(m\otimes n)=0$ gilt. Nach Definition des 
Tensorproduktes ist $({\bf 1}\otimes \omega'')_0(m\otimes n)=m \otimes
\omega''_0 n$ und $\omega''_0 n=0$ hei"st, da"s $n$ ein vakuumartiges
Element ist. Da wir $V$ als \nett vorausgesetzt haben, ist dies 
nach~\cite{DoLiMa}, Abschnitt 2.4, gleichbedeutend mit $n\in \C\cdot {\bf 1}$,
insbesondere also $n\in N_i\cong L_{(0)}$ und 
damit $m\otimes n\in W_{(0)}\otimes {\bf 1}$.
\qed

Wir werden auf $W:=W_{(0)}\oplus W_{(1)}$ eine \SVOA-Struktur definieren, 
wobei auch $W_{(1)}=0$ zugelassen sein soll.
F"ur $k=8$ bedeutet die Wahl von $W_{(1)}$ zus"atzlich die Auswahl 
von $L_{(1)}$ unter den $3$ "aquivalenten $\VSO{8}$-Moduln 
$L_{(1)}$, $L_{(2)}$ und $L_{(3)}$.
Bisher haben wir f"ur $W$ eine $U$-Modulstruktur und
 --- da $U$ in $W_{(0)}$ liegt ---
auch Elemente  ${\bf 1}\in (W_{(0)})_0$ und  \hbox{${\omega' }\in (W_{(0)})_2$.}
Zu konstruieren bleibt der Vertexoperator
$\WY(.,z):W\longrightarrow {\rm End}(W)[[z^{-1},z]]$.
Dazu zerlegen wir den Vertexoperator $\VY(.,z):V\longrightarrow {\rm 
End}(V)[[z^{-1},z]]$ von $V$ bzgl.~$U\otimes L$:
$$\VY=\bigoplus_{i,j,k\in I} 
\VY_{ij}^k,\quad \hbox{wobei\ }\VY_{ij}^k
\in\left({\normalsize M_k\otimes N_k}\atop{\normalsize M_i\otimes N_i\quad 
 M_j\otimes N_j }\right)_{U\otimes L}.$$
Nach~\cite{DoMaZhu}, Prop.~2.10, besteht der Isomorphismus
$$\left({\normalsize M_k\otimes N_k}\atop{\normalsize M_i\otimes N_i\quad 
 M_j\otimes N_j }\right)_{U\otimes L}\cong \left({\normalsize M_k}
\atop{\normalsize M_i\quad M_j}\right)_{U}\otimes \left({\normalsize N_k}
\atop{\normalsize N_i\quad N_j}\right)_{L}.$$
Sei nun $(i,j,k)\in I_a\times I_b\times I_c$ mit $a$, $b$, $c$ $\in \{0,1\}$.
Nach den Fusionsregeln f"ur $\VSO{k}\cong L$ ist $N_{ab}^c=\dim\left(
{\normalsize L_{(c)}}\atop{\normalsize L_{(a)}}\quad L_{(b)}\right)_{L}
\in \{0,1\}$. 
Die SVOA-Struktur von $\VF^{\otimes k}\cong L_{(0)}\oplus L_{(1)}$ liefert
f"ur $a+b+c\equiv 0\pmod{2}$ ein von Null verschiedenes Element
$\LY_{ab}^c\in\left({\normalsize L_{(c)}}\atop{\normalsize L_{(a)}}\quad L_{(b)}
\right)_{L}\cong \left({\normalsize L_k}\atop{\normalsize L_i\quad
L_j}\right)_{L}$.
Wir k"onnen f"ur $(i,j,k)\in I_a\times I_b\times I_c$
den Vertexoperator von $V$ daher schreiben als
\begin{equation}\label{zeryv}
\VY_{ij}^k=\cases{
\UY_{ij}^k\otimes\LY_{ab}^c, & falls $a+b+c\equiv0\pmod{2}$, \cr
0,   & falls $a+b+c\equiv 1\pmod{2}$. \cr}
\end{equation}
F"ur $W=W_{(0)}\oplus W_{(1)}=\bigoplus_{i \in I_0\cup I_1} M_i  $ 
definieren wir daher:
\begin{equation}\label{defyw}
\WY:=\bigoplus_{(i,j,k) \in I_a\times I_b\times I_c}\UY_{ij}^k, 
\end{equation}
wobei wir f"ur  $(i,j,k)\in I_a\times I_b\times I_c$,\,
 $a+b+c\,\equiv\, 1\pmod{2}$ noch $\UY_{ij}^k=0$ gesetzt haben.\hfill 
\vspace{-2mm}
\begin{lemma}\label{svoawohldefiniert}
Der in (\ref{defyw}) definierte Vertexoperator $\WY: W\longrightarrow
{\rm End}(W)[[z,z^{-1}]]$ h"angt nicht von der Auswahl der rationalen
Unter-\VOA $U\subset W_{(0)}$ ab.
\end{lemma}
{\bf Beweis:}
Seien homogene Elemente $u$, $v$ in $W_{(0)}$ bzw.~$W_{(1)}$ und $w'$ in dem
eingeschr"ankten Dualraum $W'_{(0)}$ bzw.~$W'_{(1)}$ gegeben. Sei $\mid.\mid $
wie in Definition~\ref{SVOA}.
Falls $\mid w'\mid\not\equiv 
\mid u \mid + \mid v \mid \pmod{2}$, verschwindet nach Definition 
(\ref{defyw}) die Funktion $\langle w', \WY(u,z) v\rangle$.

Nach Satz~\ref{vsok} und Satz~\ref{sdsvoafusion} ist 
\hbox{$\LY_{ab}^{c}\in\left({\normalsize L_{(c)}}\atop{\normalsize 
L_{(a)}}\quad L_{(b)}\right)_{L}\cong \C$} f"ur $a+b+c\equiv 0 \pmod{2}$ von
Null verschieden, da nach Definition $L_{(0)}\oplus L_{(1)}$ eine zu
$\VF^{\otimes k}$ isomorphe \SVOA ist.
Falls $\mid w'\mid\equiv \mid u \mid + \mid v \mid \pmod{2} $, 
gibt es nach~\cite[Prop.~11.9]{DoLe} daher Elemente 
\hbox{$\ou$, $\ov \in L_{(0)}\oplus L_{(1)}$}
bzw.~$\overline{w}' \in L'_{(0)}\oplus L'_{(1)}$ mit
$|u|=|\ou|$, $|v|=|\ov|$ und $|w'|=|\ow'|$
und \pagebreak[2]
\begin{equation}\label{lemmabed}\vspace{-4mm}
\langle \ow',\LY(\ou,z)\ov\rangle\not\equiv 0.
\end{equation}
Nach Definition gilt aber 
\begin{equation}\label{lemmadef}
\langle w',\WY(u,z)v\rangle\cdot
\langle \ow',\LY(\ou,z)\ov\rangle =
\langle w'\otimes\ow',\VY(u\otimes \ou,z)v\otimes \ov\rangle.
\end{equation}
Aus (\ref{lemmabed}) folgt, da"s
$\langle w', \WY(u,z) v\rangle$ eine in $z$ rationale Funktion ist,
d.h.~die Laurent\-entwicklung von $\WY$ in $0$ ist wohldefiniert. \qed

Insbesondere stimmt $\WY|_{W_{(0)}}$ mit der $W_{(0)}$-Modulstruktur auf 
$W$ "uberein, welche auch ohne R"uckgriff auf die Intertwinerr"aume definiert
werden kann. W"urde man andererseits $\WY$ durch Gleichung (\ref{lemmadef}) definieren, w"are es schwieriger das folgende Resultat zu beweisen.
\begin{satz}\label{satzvoatosvoa}
Das Tupel $(W,\WY,{\bf 1},\omega')$ mit dem in~(\ref{defyw}) definierten 
Vertexoperator $\WY$ ist eine \nette \SVOA vom Rang $c-\frac{k}{2}$.
\end{satz}
{\bf Beweis:}
Da $\WY$ als direkte Summe von $U$-Intertwinern geschrieben ist, sind alle 
Axiome bis auf die Jacobi Identit"at klar.

Die Jacobi Identit"at f"ur \SVOAs ist bei Annahme der "ubrigen Axiome
nach~\cite{DoLe}, Prop.~7.16, "aquivalent zur verallgemeinerten Rationalit"at 
und Kommutativit"at und zu zwei Bedingungen an $L_0$ und $L_1$. Bis auf 
die Kommutativit"at sind diese offensichtlich erf"ullt.

Zu zeigen ist somit f"ur alle $w$, $u$, $v$ in $W_{(0)}$ oder $W_{(1)}$
bzw.~$w'$ in $W_{(0)}'$ oder $W_{(1)}'$ die Identit"at
\begin{equation}\label{Wkommutativitaet}
\langle w',\WY(u,z_1)\WY(v,z_2)w\rangle=(-1)^{|u||v|}
\langle w',\WY(v,z_2)\WY(u,z_1)w\rangle
\end{equation}
zwischen rationalen Funktionen in $z_1$ und $z_2$.
Wegen Definition~(\ref{defyw}) verschwinden 
beide Seiten von~(\ref{Wkommutativitaet}),
falls $|w'|\not\equiv |u|+|v|+|w| \pmod{2}$ ist.
Mit dem gleichen Argument wie beim Beweis der vorangehenden Lemmas
findet man Elemente \hbox{$\ow$, $\ou$, $\ov \in L_{(0)}\oplus L_{(1)}$}
bzw. \hbox{$\overline{w}' \in L'_{(0)}\oplus L'_{(1)}$} mit
\begin{equation}\label{bedi1}
|u|=|\overline{u}|,\quad |v|=|\overline{v}|,\quad
|w|=|\overline{w}|\hbox{\quad und\quad } |w'|=|\overline{w'}|
\end{equation}
und
\begin{equation}\label{bedi2}
\langle \ow',\LY(\ou,z_1)\LY(\ov,z_2)\ow\rangle=(-1)^{|\ou||\ov|}
\langle \ow',\LY(\ov,z_2)\LY(\ou,z_1)\ow\rangle\not=0.
\end{equation}
Sei also 
$|w'|\equiv |u|+|v|+|w| \pmod{2}$
und seien die
$\ou$, $\ov$, $\ow$, $\ow'$ so gew"ahlt, da"s sie die Bedingungen 
(\ref{bedi1}) und (\ref{bedi2}) erf"ullen. Wir erhalten:
\begin{eqnarray*}
& &\langle w',\WY(u,z_1)\WY(v,z_2)w\rangle\cdot
   \langle \ow',\LY(\ou,z_1)\LY(\ov,z_2)\ow\rangle\hfill \\
&=&\langle w'\otimes\ow',(\WY(u,z_1)\otimes\LY(\ou,z_1))
(\WY(v,z_2)\otimes\LY(\ov,z_2))w\otimes\ow\rangle   \\
\noalign{\hbox{nach Definition von $\WY$}\vspace{-5mm}}           \\
&=&\langle w'\otimes\ow',\VY(u\otimes\ou,z_1)
\VY(v\otimes\ov,z_2)w\otimes\ow\rangle             \\
\noalign{\hbox{wegen der Jacobi Identit"at von $V$}\vspace{-5mm}}  \\
&=&\langle w'\otimes\ow',\VY(v\otimes\ov,z_2)
\VY(u\otimes\ou,z_1)w\otimes\ow\rangle   \\
\noalign{\hbox{nach Definition von $\WY$}\vspace{-5mm}}           \\
&=&\langle w'\otimes\ow', (\WY(v,z_2)\otimes\LY(\ov,z_2))
(\WY(u,z_1)\otimes\LY(\ou,z_1))w\otimes\ow\rangle   \\
&=&\langle w',\WY(v,z_2)\WY(u,z_1)w\rangle\cdot
    \langle \ow',\LY(\ov,z_1)\LY(\ou,z_2)\ow\rangle \\
&=&\langle w',\WY(v,z_2)\WY(u,z_1)w\rangle\cdot (-1)^{|\ou||\ov|}
  \langle \ow',\LY(\ou,z_2)\LY(\ov,z_1)\ow\rangle.
\end{eqnarray*}
Hieraus erh"alt man schlie"slich wegen der Bedingungen (\ref{bedi1}) und 
(\ref{bedi2}) die zu zeigende Gleichung~(\ref{Wkommutativitaet}).
\qed
\begin{vermutung}\label{voa2svoa}
Ist die \VOA $V$ selbstdual, so ist die so erhaltene \SVOA $W$ rational, selbstdual und \sehrnett. 
Die Unter-\VOA $W_{(0)}$ besitzt (falls $W_{(1)}\not=0$)
zus"atzlich zu $W_{(0)}$ und $W_{(1)}$ noch 
die irreduziblen Moduln $W_{(2)}$ (und $W_{(3)}$, falls $k$ gerade). 
Ist weiter $V$ unit"ar, so ist es auch $W$.
\end{vermutung}
Liegt die $V$-Unter\VOA $L^k=\VSO{k}$ in einer $V$-Unter\VOA $L^l=\VSO{l}$
f"ur ein \hbox{$l>k$}
(von der nat"urlichen Einbettung $\SO{k}\subset\SO{l}$
herkommend), so l"a"st sich der Zusammenhang zwischen der wie oben mit
$L^k$ konstruierten \SVOA $W$ und der mit $L^l$ konstruierten \SVOA $W'$
leicht beschreiben. Dazu betrachte man die folgende Zerlegung von 
$\VF^{\otimes l}$ in den geraden und den ungeraden Anteil:
\begin{eqnarray*}
\VF^{\otimes l} & = & \VF^{\otimes l-k}\otimes \VF^{\otimes l}
 =L_{(0)}^l\oplus L_{(1)}^l \\
& = &
(L_{(0)}^{l-k}\otimes L_{(0)}^{k}\oplus L_{(1)}^{l-k}\otimes L_{(1)}^{k})
\oplus
(L_{(0)}^{l-k}\otimes L_{(1)}^{k}\oplus L_{(1)}^{l-k}\otimes L_{(0)}^{k}).
\end{eqnarray*}
Mit $W'_{(0)}={\rm Com}_V(L^l)$ k"onnen wir daher $V$ wie folgt zerlegen:
\begin{eqnarray*}
V &= & \bigoplus_{a=0,1,2(,3)} W'_{(a)}\otimes L_{(a)}^l \\
  &= & W'_{(0)} \otimes
(L_{(0)}^{l-k}\otimes L_{(0)}^{k}\oplus L_{(1)}^{l-k}\otimes L_{(1)}^{k})
\oplus W'_{(1)} \otimes
(L_{(0)}^{l-k}\otimes L_{(1)}^{k}\oplus L_{(1)}^{l-k}\otimes L_{(0)}^{k})
\oplus \cdots \\
 &= &(W'_{(0)}\otimes L_{(0)}^{l-k} \oplus W'_{(1)}\otimes L_{(1)}^{l-k})
\otimes L_{(0)}^k \oplus 
(W'_{(0)}\otimes L_{(1)}^{l-k} \oplus W'_{(1)}\otimes L_{(0)}^{l-k})
\otimes L_{(1)}^k \oplus \cdots\\
  &= & \bigoplus_{a=0,1,2(,3)} W_{(a)}\otimes L_{(a)}^k.
\end{eqnarray*}
Somit ergibt sich die Zerlegung
\begin{eqnarray}
W & = & W_{(0)}\oplus W_{(1)}   \nonumber \\
  & = & (W'_{(0)}\otimes L_{(0)}^{l-k} \oplus W'_{(1)}\otimes L_{(1)}^{l-k})
\oplus (W'_{(0)}\otimes L_{(1)}^{l-k} \oplus W'_{(1)}\otimes L_{(0)}^{l-k})
\nonumber \\ \label{Wunterschiedlich}
  & = & (W' \otimes \VF^{\otimes l-k})_{(0)}\oplus
 (W' \otimes \VF^{\otimes l-k})_{(1)}=
 W' \otimes \VF^{\otimes l-k}.
\end{eqnarray}\nopagebreak[2]
Es ist leicht einzusehen, da"s diese Zerlegungen alle mit der jeweiligen \OVOA-Struktur vertr"aglich sind.
\pagebreak[2]

\section{Die Konstruktion von \VOAs aus \SVOAs{ }}\label{rueckR}

Die umgekehrte Konstruktion von \VOAs aus \SVOAs ist nicht so einfach.
Im Gegensatz zu der Situation bei Codes und Gittern fehlt ein "`einbettender
Raum"' wie $\F_2^n$ oder $\R^n$, von dem der konstruierte Kandidat
in nat"urlicher Weise die Struktur einer \VOA erbt.
Der der Konstruktion zu Grunde liegende Vektorraum und die zugeh"orige
\VOA-Struktur lassen sich trotzdem leicht beschreiben.

Sei also $W$ eine selbstduale \sehrnette rationale \SVOA vom Rang $c'$,
die nicht schon eine \VOA ist. Betrachte die \SVOA $\VF^{\otimes k}=
\L{0}\oplus \L{1}$ vom Rang $\frac{k}{2}$, wobei $k$ so gew"ahlt
sei, da"s $c'+\frac{k}{2}\in 8\Z$. Die R"ange $c'$ und $\frac{k}{2}$
liegen dann entweder beide in $\frac{1}{2}\Z\setminus\Z$, $\Z\setminus 2\Z$
oder $2 \Z$ und die Fusionsalgebren von $L_{(0)}$ und $W_{(0)}$ sind 
isomorph (vgl.~Satz~\ref{sdsvoafusion}).
Die Summe $V:=\bigoplus_{a=0,1,2(,3)}W_{(a)}\otimes L_{(a)}$ ist dann ein
$W_{(0)}\otimes L_{(0)}$-Modul, der eine
{\it ganzzahlige} Graduierung durch die Eigenwerte von $\omega_1$ 
des Virasoroelement $\omega$ von $W_{(0)}\otimes L_{(0)}$ besitzt.
Die Intertwinerr"aume
$$\left( W_{(k)}\otimes L_{(k)} \atop W_{(i)}\otimes L_{(i)}\quad
W_{(j)}\otimes L_{(j)} \right)_{W_{(0)}\otimes L_{(0)}}\cong 
\left( W_{(k)} \atop W_{(i)}\quad W_{(j)} \right)_{W_{(0)} }\otimes
\left( L_{(k)} \atop L_{(i)}\quad L_{(j)} \right)_{L_{(0)}} $$
sind aufgrund der Struktur der Fusionsalgebra von $W$ und $L$
alle null- oder eindimensional.
Man kann daher geeignet normierte $\VY_{ij}^k\in \left( W_{(k)}\otimes L_{(k)}
\atop W_{(i)}\otimes L_{(i)}\quad W_{(j)}\otimes L_{(j)} \right)_{W_{(0)}
\otimes L_{(0)}}$ ($\VY_{ij}^k\not=0$, falls $N_{ij}^k=1$) w"ahlen
und schlie"slich $$\VY=\bigoplus_{i,j,k=0,1,2(,3)} \VY_{ij}^k $$
setzen.

Wenn die \SVOA $W$ schon eine \VOA ist, setzen wir $V=W\otimes V_{D_k^+}$,
wobei $V_{D_k^+}\cong \L{0}^k\oplus \L{2}^k$ die selbstduale \VOA zum
f"ur $k\in 8\Z$ geraden selbstdualen Gitter $D_k^+$ ist.
\begin{vermutung}\label{svoa2voa}
Das Tupel $(V,\VY,{\bf 1},\omega)$ wird so zu einer selbstdualen \netten
rationalen \VOA. Ist weiter $W$ unit"ar, so ist es auch $V$.
\end{vermutung}\nopagebreak[2]
Die Konstruktion von $V$ ist f"ur $c'\in\Z$ unabh"angig von der gemachten
Auswahl zwischen $W_{(2)}$ und $W_{(3)}$, da $L_{(2)}$ und $L_{(3)}$
durch einen "au"seren Automorphismus von $\so{2l}=D_l$, $l=\frac{k}{2}\in\Z$,
der auf der Fusionsalgebra von $L_{(0)}$ operiert, vertauscht werden.
\pagebreak[2]

Wenn wir die beiden Vermutungen~\ref{voa2svoa} und~\ref{svoa2voa} als
richtig voraussetzen, so ist klar, da"s die Konstruktionen von
\VOAs zu \SVOAs und umgekehrt von \SVOAs zu \VOAs
invers zueinander sind. Man erh"alt genauer eine
Bijektion zwischen den Isomorphieklassen von selbstdualen \SVOAs $W$ vom 
Rang $c'$ und den Isomorphieklassen von Paaren $(V,\VSO{k})$, wobei $V$
eine selbstduale \VOA vom Rang  $c=c'+\frac{k}{2}\in 8\Z$ und $\VSO{k}$ eine
Unter-\VOA vom Rang $\frac{k}{2}$ ist (zusammen mit der Auswahl des 
$\VSO{k}$-Moduls $L_{(1)}$, falls $k=8$ ist).
Isomorphieklasse des Paares $(V,\VSO{k})$ bedeutet Fixierung der
Isomorphieklasse von $V$ und dann Wahl der Unter-\VOA $\VSO{k}$ bis auf
"Aquivalenz unter $\rm{Aut(V)}$, den Automorphismen von $V$. (Falls $k=8$,
mu"s die  Operation von der Untergruppe von $\rm{Aut(V)}$, die $\VSO{8}$ 
fixiert, auf den drei $\VSO{8}$-Moduln betrachtet werden: Zwei Auswahlen von
$\VSO{8}$-Moduln sind  "aquivalent, falls ein Automorphismus sie ineinander
"uberf"uhrt.)
Wenn die \SVOA $W$ schon eine \VOA ist, so ist diese Beziehung offenbar auch 
richtig.

Verwendet man noch Satz~\ref{satzabspalt} und Zerlegung 
(\ref{Wunterschiedlich}), so kann man die Klassifikation der selbstdualen 
\SVOAs mit $W_{1/2}=0$ auf die Bestimmung der {\it maximalen\/} 
$\VSO{k}$-Unteralgebren von $V$ zur"uckf"uhren,
da diese das erste Glied in der Folge $W$, $W\otimes \VF$, 
$W \otimes \VF^{\otimes 2}$, $\ldots$ bilden.

Die Konstruktionen aus diesem Abschnitt sind in dieser Arbeit allerdings nur
skizziert, und Teilschritte sind noch zu beweisen. Evt.~m"ussen au"ser
selbstdual (sehr) \nett unit"ar und rational noch weitere technische
Forderungen an die \OVOAs gestellt werden. Sicherlich ist aber die Klasse
von \OVOAs, f"ur die die "Uberlegungen dieses Kapitels gelten, diejenige, die
interessant ist und die man betrachten sollte.

Wir formulieren den beschriebenen Zusammenhang noch als 
\begin{vermutung}\label{umkehr}
Es besteht eine umkehrbar ein\-deutige Be\-zie\-hung zwi\-schen den 
Iso\-mor\-phie\-klas\-sen von selbst\-dua\-len \sehrnetten uni\-t"aren 
ra\-tio\-nalen \SVOAs vom Rang $c'$ und
den Isomorphieklassen von selbstdualen \netten unit"aren rationalen \VOAs
vom Rang $c$ zusammen mit einer Unter-\VOA $\VSO{k}$ vom Rang 
$\frac{k}{2}=c-c'$ (und Auswahl eines $\VSO{k}$-Moduls falls $k=8$).
\end{vermutung}

Es ist einfach, die Berechnung der Auto\-mor\-phismen\-grup\-pe von $W$ auf das 
Stu\-dium der Ope\-ra\-tion der Auto\-mor\-phis\-men\-grup\-pe von $V$ auf der 
Un\-ter\-al\-ge\-bra $\VSO{k}$ zur"uck\-zuf"uhren.
Wir werden dies hier nicht genauer untersuchen, sondern im n"achsten Kapitel 
nur den dort ben"otigten Zusammenhang zwischen den Automorphismengruppen 
von $\VM$ und $\VB$ beschreiben.

{\bf Beispiele f"ur die Beziehung zwischen \VOAs und \SVOAs{}:}\newline
Wenn $k\geq 3$ ist, wird die \VOA $L=\VSO{k}$ von den Elementen aus 
$(\VSO{k})_{1}\cong\so{k}$ erzeugt. Die Klassifikation von \VOA-Unteralgebren
$\VSO{k}$ in $V$ kann daher im wesentlichen auf die Klassifikation
von Lieunteralgebren $\so{k}$ in der reduktiven Liealgebra $V_1=g$
zur"uckgef"uhrt werden (vgl.~\cite{Fuchs}, Abschnitte~1.8, 2.8, 3.8 und 
\cite{McPa,Dy}).

{\it Bereich $0\leq c< 8$:\/}
Wir erhalten so wieder das Resultat aus Korollar~\ref{SVOAS1-8}:
W"ahlen wir $\rang V=8$, so existiert nur die \VOA $V_{E_8}\cong
V_{E_{8,1}}$. In der Liealgebra $E_8$ gibt es bis auf Isomorphie nur
eine Klasse von regul"aren (konformen) Einbettungen der $\so{k}$, $k\geq 3$
(vgl.~\cite{McPa}). Die $\so{16}\cong D_8$ ist maximal.
F"ur $0\leq c \leq 6\frac{1}{2}$ gilt daher $W\cong\VF^{\otimes 2c}$.
F"ur $c=7$ und $c=7\frac{1}{2}$ ist es einfacher, die beiden \VOAs
vom Rang $16$ zu betrachten: $V_{E_8\times E_8}$ und $V_{D_{16}^+}$.
In $E_8\times E_8$ ist eine $\so{k}$, $k>16$ nicht einbettbar, und in
$D_{16}$ ist die $\so{32}$ auch maximal, d.h.~auch f"ur $c=7$ oder 
$c=7\frac{1}{2}$ gilt $W\cong\VF^{\otimes 2c}$.

{\it Bereich $8\leq c< 16$:\/}
F"ur den Bereich $8\leq c\leq 16$ ist es am einfachsten, Schellekens
Liste~\cite{schellekens1} der selbstdualen \VOAs mit Rang $c=24$ heranzuziehen. 
Denn f"ur $k\geq 9$ gibt es bis auf Isomorphie immer nur h"ochstens
eine regul"are (konforme)
maximale Einbettung von $\so{k}$ in eine einfache Liealgebra $g$ 
(in dem Sinne, da"s $\widetilde{V}_{{\bf SO}(k),1}$ 
\VOA-Unteralgebra von $\widetilde{V}_1$).
Wir erhalten die folgende Liste:
\begin{equation}\label{listeSVOAS8-16}
\begin{array}{l|ccccc}
c              &  8      & 12           & 14  & 15 &\frac{31}{2}  \\ \hline
{\rm \SVOA{ }} & V_{E_8} & V_{D_{12}^+} & V_{(E_7+E_7)^+} & 
V_{A_{15}^+} & V_{E_{8,2}^+} \\
\end{array}
\end{equation}
Man "uberpr"uft leicht, da"s die so erhaltenen Beispiele f"ur 
$c \in \{8,12,14,15\}$ mit den entsprechenden Gittertheorien (siehe 
Satz~\ref{satzgitterlieiso})
"ubereinstimmen ($V_{E_8}$ ist sogar \hbox{eine \VOA{}).}
F"ur $c=15\frac{1}{2}$ erh"alt man die "`neue"' \SVOA $V_{E_{8,2}^+}$
--- die zur Level $2$
Darstellung von $E_8$ geh"orende \VOA zusammen mit dem Modul mit konformen
Gewicht $\frac{1}{2}$ ---, die allerdings in der physikalischen Literatur schon
untersucht wurde~\cite{ScheYu} (zumindest die Fusionsalgebra); vgl.~auch
\cite{DML-simcurr}, Bemerkung 5.9.
\begin{vermutung}\label{SVOAS8-16}
Die Liste in~(\ref{listeSVOAS8-16}) enth"alt alle selbst\-dua\-len 
un\-zer\-leg\-ba\-ren SVOAs in dem Be\-reich \hbox{$8\leq c < 16$}. 
\end{vermutung}
Unsere obigen "Uberlegungen haben zumindest das folgende gezeigt:
\begin{satz}\label{satzSVOAS8-16}
Die Liste der in (\ref{listeSVOAS8-16}) angegebenen selbst\-dualen unzerleg\-baren \SVOAs f"ur den Bereich \hbox{$8\leq c < 16$} 
ergibt sich aus der Annahme der Umkehr\-vermutung~\ref{umkehr} und der Vermutung~\ref{vermutungsdvoa}
f"ur $c=16$ bzw.~$c=24$.
\end{satz}

Die Klassifikation der selbstdualen \SVOAs f"ur den Bereich
$16\leq c <24$ ist mit Hilfe von Schellekens Liste~\cite{schellekens1}
ebenfalls m"oglich.
F"ur $V\not=\VM$ ist dies, wie schon oben gesehen,
auf die Klassifikation von Einbettungen der Liealgebra $\so{k}$ in die 
Liealgebra $V_1$ zur"uckf"uhrbar. Dies wird in einer separaten Arbeit
genauer ausgef"uhrt~\cite{hoehn2}.
Im besonders interessanten Fall des Mondscheinmoduls $\VM$ ist $\VM_1=0$,
und daher ist eine maximale Unter-\VOA h"ochstens eine $\VSO{1}$.
Tats"achlich l"a"st sich eine "Aquivalenzklasse von zur 
$(c,h)=(\frac{1}{2},0)$-H"ochstgewichtsdarstellung der Virasoroalgebra 
assozierten Unter-\VOAs $\La\cong\VSO{1}$ finden. Die auf diese Weise erhaltene
Babymonster-\SVOA $\VB$ wird im n"achsten Kapitel genauer untersucht werden.

 \chapter{Die Babymonster Vertexoperator-Superalgebra}\label{babykapitel}

Das Ziel dieses Kapitels ist die Konstruktion einer \SVOA $\VB$ vom Rang
$23\frac{1}{2}$ auf der das Babymonster $B$ durch Automorphismen operiert.
Diese \SVOA ist das Analogon des k"urzeren Golay Codes $g_{22}$ 
(s.~\cite{PlSl}) bzw.~des k"urzeren Leechgitters $O_{23}$ 
(s.~\cite{CoSl-ul23}), welche eine zweifache Erweiterung der
Mathieu Gruppe $M_{22}$ bzw.~der Conway Gruppe $Co_2$ als 
Automorphismengruppe besitzen. So wie $g_{22}$ keine Tetraden
(Vektoren vom Gewicht $4$) und $O_{23}$ keine Wurzeln (Vektoren der 
Quadratl"ange $2$) besitzt, so enth"alt $\VB$ keine nichttriviale Lie Unteralgebra $V_1$
(Vektoren vom konformen Gewicht 1). Die \SVOA $\VB$ ist das
nat"urliche Objekt, auf dem $B\times\Z_2$ durch Automorphismen operiert,
und in diesem Sinne liefert sie die nat"urlichste Definition des 
Babymonsters.

In Abschnitt 4.2 wird $\VB$ aus dem Mondscheinmodul $\VM$ unter Verwendung
der Metho\-den des vorherigen Kapitels konstruiert.
Dazu wird in Abschnitt 4.1 eine Beschreibung von Virasoro 
Unter\-vertex\-operator-\-Algebren vom Rang $\frac{1}{2}$ in $\VM$ gegeben.

\section{Zerlegung des Mond\-schein\-moduls $\VM$ unter $\Lt$ }

Der Mondscheinmodul kann, da $\VM_1=0$, 
nicht in eine Summe von H"ochst\-gewichts\-darstellungen von affinen Kac-Moody Algebren zerlegt
werden. Wie Dong, Mason und Zhu in~\cite{DoMaZhu} gezeigt haben, kann man aber
$48$ paarweise zur $c=\frac{1}{2}$ unit"aren H"ochstgewichtsdarstellung der
Virasoroalgebra assozierte kommutierende Unter-\VOAs $\La$ vom Rang
$\frac{1}{2}$ finden. Die Monster-\VOA zerlegt sich dann als direkte
Summe von \VOA-Moduln von $\Lt$. Wir werden in diesem Abschnitt
zeigen, da"s dieses Resultat auch aus einer Arbeit~\cite{MeNe} von W.~Meyer und
W.~Neutsch "uber assoziative Unteralgebren der Griess Algebra
folgt. Weiter geben wir mit Hilfe von Invariantentheorie eine neue genauere
Beschreibung der Zerlegung von $\VM$.

\medskip

Die $196884$-dimensionale kommutative nicht assoziative Griess Algebra
${\cal B}$~\cite{Gr} ist der Gewicht $2$ Anteil $\VM_2$ des Mondscheinmoduls.
F"ur zwei Elemente $a$, $b \in {\cal B}$ ist das Algebraprodukt von ${\cal B}$
gegeben durch $a \times b=a_1b$ und die nichtsingul"are symmetrische
Bilinearform $\langle .,. \rangle$ durch $a_3 b=\langle a,b \rangle \cdot
{\bf 1}$. Sie stimmt mit der in Satz~\ref{bilinearform} beschriebenen
invarianten symmetrischen  Bilinearform "uberein, d.h.~es gilt
$$\langle a \times b,c \rangle = \langle a, b \times c \rangle $$
f"ur Elemente $a$, $b$, $c\in {\cal B}$. Das Einselement von ${\cal B}$
ist $e=\frac{1}{2}\omega$, die H"alfte des Virasoroelementes, und hat die
Norm $\langle e, e \rangle=\frac{1}{4}\cdot \omega_3\omega=\frac{c}{8}=3$
(vgl.~\cite{FLM}, Prop.~10.3.6). Da $\VM_1=0$, gilt noch $a_2 b=0$.

Die Struktur von assoziativen Unteralgebren von ${\cal B}$ wurde
in~\cite{MeNe} untersucht. Wir fassen die dort erhaltenen allgemeinen
Strukturresultate zusammen in
\begin{satz}[Meyer, Neutsch]\label{meyer-neutsch}
Sei ${\cal U}$ eine $k$-dimensionale assoziative Unteralgebra
(der reellen Form) der Griess Algebra ${\cal B}$.
Dann gilt
\begin{list}{}{}
\item
1) ${\cal U}$ ist isomorph zu einer direkten Summe von $k$ Kopien von
$\R${\rm :} ${\cal U} \cong \R^k$.
\item
2) ${\cal U}$ enth"alt eine Basis von $k$ paarweise sich anni\-hilieren\-den
idem\-poten\-ten Ele\-menten $a_i$ die ortho\-gonal zueinander sind:
$a_i \times a_j=0 $ und $\langle a_i, a_j \rangle=0$ f"ur alle $1\leq i,j
\leq k$ mit $i\not=j$.
\item
3) Die idempotenten Elemente von ${\cal U}$ sind genau die $2^k$ Elemente
$a_{i_1}+a_{i_2}+\cdots+a_{i_{\nu}}$, $1\leq i_1 <i_2 < \dots < i_v \leq k$,
wobei $\{ a_1, \dots a_k \}$ die einzige Orthogonal Basis unter ihnen ist.
Die $a_1$, $\dots$, $a_k$ hei"sen daher die Fundamental-Idempotenten von
${\cal U}$.
\item
4) ${\cal U}$ ist genau dann maximal assoziativ, wenn die folgenden beiden
Bedingungen erf"ullt sind:
\begin{list}{}{}
\item
a) Das Einselement $e$ von ${\cal B}$ liegt in ${\cal U}$.
\item
b) Die Fundamental-Idempotenten $a_1$, $\cdots$, $a_k$ sind unzerlegbar,
d.h.~sie lassen sich nicht als Summe von zwei oder mehr von Null
verschiedenen Idempotenten schreiben.
\end{list}
\end{list}
\end{satz}
Der Beweis beruht im wesentlichen auf der Nortonungleichung~\cite{Co-monster}
$$\langle a\times a,b\times b\rangle\geq\langle a\times b,a\times
 b\rangle$$
f"ur Elemente $a$, $b$ aus ${\cal B}$.

Eine Verbindung zwischen assoziativen Unteralgebren von ${\cal B}$ und
Unter-\VOAs von $\VM$ stellt nun der folgende Satz her.
\begin{satz}\label{ass2vir}\footnote{Der Zusammenhang zwischen Idempotenten
in ${\cal B}$ und Virasoroalgebren ist auch in~\cite{Mia} betrachtet worden.} 
Sei ${\cal U}$ eine $k$-dimensionale maximale assoziative Unter\-algebra
der Griess Algebra ${\cal B}$ mit den Fundamental\--Idempotenten $a_1$,
$\dots$, $a_k$.
Dann bildet der von den Komponenten $(2a_i)_n$ des
Vertex\-operators $Y(2a_i,z)$ vom Vakuum ${\bf 1}$ erzeugte Unter\-vektor\-raum
von $\VM$ eine zur Virasoro\-algebra assozierte Unter-\VOA mit Virasoroelement
$\omega_i=2 a_i$ vom Rang $c_i=8\langle a_i,a_i\rangle$. Diese $k$
Unter-\VOAs kommutieren paarweise, und ihr Tensor\-produkt ist eine Unter-\VOA
von $\VM$ mit dem gleichen Virasoro\-element $\omega=\sum_{i=1}^k \omega_i$
wie $\VM$.
\end{satz}
{\bf Beweis: }
F"ur ein Idempotent $a\in {\cal B}$ gilt $a_3 a=a$ und $a_2 a=0$. Nach
Satz~\ref{L2-vir} erzeugen die Koeffizienten von $Y(2a,z)$ daher eine
Unter-\VOA vom Rang $8\langle a,a\rangle$, die zu einer Darstellung der
Virasoroalgebra assoziert ist.
F"ur zwei Virasoroelemente $\omega_i\not=\omega_j$ gilt aus Dimensionsgr"unden
$(\omega_i)_2\omega_j=0$ und $(\omega_i)_n\omega_j=0$ f"ur $n\geq 4$.
Nach Voraussetzung und Satz~\ref{meyer-neutsch}, 2) gilt $(\omega_i)_1\omega_j=
(\omega_i)_3\omega_j=0$. Der Vektor $\omega_j$ ist daher ein 
H"ochstgewichtsvektor vom H"ochstgewichts $h=0$ f"ur die von den Koeffizienten
$L_i(n):=(\omega_i)_{n+1}$ von $Y(\omega_i,z)=\sum_{n\in\Z}(\omega_i)_n 
z^{-n-1}$ erzeugte Virasoroalgebra $Vir_i$. Der von $\omega_j$ erzeugte 
$Vir_i$-Untermodul von $\VM$ ist unit"ar und mu"s daher zu einem Quotienten
des Virasoromoduls $M_c={\cal U}(Vir^-){\bf 1}/\langle L_{-1}{\bf 1}\rangle$  
isomorph sein. F"ur diesen Vakuummodul gilt aber $L_i(-1)\omega_j=
(\omega_i)_0(\omega_j)=0$. Die Kommutatorformel (vgl.~\cite{FHL} (2.3.1)) 
liefert nun $[Y(\omega_i,z_1),Y(\omega_j,z_2)]={\rm res}_{z_0}z_2^{-1}
\delta(\frac{z_1-z_0}{z_2})Y(Y(\omega_i,z_0)\omega_j,z_2)=0$, da
$\sum_{n\geq 0}(\omega_i)_n\omega_j z_0^{-n-1}=0$ und nur diese singul"aren
Terme bei der Berechnung des Residuums ber"ucksichtigt werden. Die linearen
Abbildungen $(\omega_i)_n$ und $(\omega_j)_n$ kommutieren somit f"ur alle 
$m$, $n\in \Z$. Dann kommutieren aber auch alle Vertexoperatoren der von 
$\omega_i$ und $\omega_j$ erzeugten Unter-\VOAs, denn diese lassen sich nach Satz~\ref{L2-vir} durch die $(\omega_i)_m$ bzw.~$(\omega_j)_n$ ausdr"ucken.
\qed

Da $\VM$ unit"ar ist, sind die zu den Idempotenten geh"origen Unter-\VOAs
entweder die \VOAs zu den Virasoroh"ochstgewichtsdarstellungen $M_c=
{\cal U}(Vir^-){\bf 1}/\langle L_{-1}{\bf 1}\rangle$ ($c\geq 1$) oder die 
unit"aren
Darstellungen $L_c(0)=M_c/I_c$ der minimalen Serie, wobei hier f"ur den
Rang nur die Werte $c=1-\frac{6}{n(n+1)}$, $n=3$, $4$, $\dots$ zul"assig
sind~(s.~\cite{La} oder~\cite{min-unitaer}, Kap.~12).
Der kleinste m"ogliche Wert ist $c=\frac{1}{2}$, und wir erhalten das schon
in~\cite{MeNe} vermutete
\begin{korollar}
Die Dimension einer assoziativen Unteralgebra der Griess Algebra ist
h"ochstens $48$.
\end{korollar}
Die einzige bekannte Klasse von Idempotenten mit minimaler Norm
$\frac{1}{16}$ liefern die {\it Transpositions-Idempotenten}
(vgl.~\cite{Co-monster,atlas,MeNe}): Zu jeder Involution $\alpha$ vom Typ
$2A$ im Monster  ("`Transpositionen"') gibt es ein als axialen Vektor
$t_{\alpha}$ bezeichnetes Element in der Griess Algebra. Das zugeh"orige
Transpositions-Idempotent ist dann $i_{\alpha}=\frac{1}{64}t_{\alpha}$.
Auf der Menge $\{i_{\alpha}\}\subset{\cal B}$ der Transpositions-Idempotenten
operiert $M$ so, wie es auf der Menge $\{\alpha\}\subset M$ der
$2A$-Involutionen durch Konjugation operiert.
Zwei Transpositions-Idempotente $i_{\alpha}$ und $i_{\beta}$ annihilieren
sich gegenseitig genau dann, wenn $\alpha\beta$ eine Involution vom Typ $2B$
in $M$ ist.
Es gilt nun der
\begin{satz}[\cite{MeNe} und \cite{DoMaZhu}]\label{48idem}
Es gibt $48$-dimensionale assoziative Unter\-algeb\-ren in ${\cal B}$ mit
Transpositions-Idempotenten als Fundamental-Idempotente.
\end{satz}
Nach~\cite{No,GMS}
operiert das Monster als Rang $9$ Permutationsgruppe auf den
$2A$-Involutionen, und die Bahnen von Paaren von Involutionen sind durch
die $M$-Konjugationsklassen ihres Produktes gegeben, d.h.~alle Paare
$(\alpha,\beta)$ mit $\alpha\cdot \beta$ in $2B$ sind unter $M$ "aquivalent.
Daher sind auch alle Paare $(i_{\alpha},i_{\beta})$ von kommutierenden
Transpositions-Idempotenten $i_{\alpha}$ und $i_{\beta}$ "aquivalent.

Die Konstruktion von $48$-Tupel $S$ von Transpositions-Idempotenten wie in
Satz~\ref{48idem} bzw.~"aquivalent dazu $48$-Tupel von $2A$-Involutionen
in $M$, deren paarweises Produkt eine $2B$-Involution ist, ist modulo der
Operation von $M$ hingegen nicht eindeutig.
Betrachten wir den Monstergraphen, dessen Ecken die $2A$-Involutionen sind
und dessen Kanten mit der \hbox{$M$-Konjugationsklasse} des Produktes der 
Ecken gef"arbt sind, so interessieren wir uns f"ur die Orbiten
voll\-st"andiger Unter\-graphen vom Kantentyp $2B$. Diese gruppen\-theoretische
Um\-formulierung sollte mit gruppen\-theoretischen Methoden zu beantworten sein
(vgl.~\cite{No,Iv,CoMc}). Sei nun ein festes $48$-Tupel $S$ fixiert.

Wegen Satz~\ref{48idem}, Satz~\ref{ass2vir} und Satz~\ref{tensorprodukt}
zerlegt sich $\VM$ als direkte Summe von $\Lt$-Moduln:
\begin{equation}\label{monzer}
\VM=\bigoplus_{h_1,\ldots,h_{48}\in\{0,\frac{1}{2},\frac{1}{16}\}}
c_{h_1,\ldots,h_{48}}\, L_{1/2}(h_1,\ldots,h_{48}),
\end{equation}
wobei $L_{1/2}(h_1,\ldots,h_{48})=L_{1/2}(h_1)\otimes\ldots\otimes
L_{1/2}(h_{48})$ das $48$-fache Tensorprodukt von den drei irreduziblen
$\La$-Moduln $\La$, $\Lb$ und $\Lc$ ist.

Seien $i$, $j$ und $k$ nichtnegative Zahlen mit $i+j+k=48$. Wir sagen, da"s
der Modul $L_{1/2}(h_1,\cdots,h_{48})$ den Typ $(i,j,k)$ hat, falls
$\#\{\nu\mid h_{\nu}=0\}=i$, $\#\{\nu\mid h_{\nu}=\frac{1}{2}\}=j$ und
$\#\{\nu\mid h_{\nu}=\frac{1}{16}\}=k$ ist.
Die Vielfachheit $m_{i,j,k}$, mit der $\Lt$-Moduln vom
Typ $(i,j,k)$ in $\VM$ vorkommen, ist dann
\begin{equation}\label{mijk}
m_{i,j,k}=\sum_{\begin{array}{c}\scriptstyle
h_1,\ldots,h_{48}\in\{0,\frac{1}{2},\frac{1}{16}\} \\
\scriptstyle\#\{\nu\mid h_{\nu}=\frac{1}{2}\}=j \\
\scriptstyle\#\{\nu\mid h_{\nu}=\frac{1}{16}\}=k
\end{array}}
c_{h_1,\ldots,h_{48}}.
\end{equation}
Schlie"slich setzen wir noch
\begin{equation}\label{pvsabc}
P_{\VM}^S(a,b,c)=\sum_{{\scriptstyle i,j,k\in\Z_+}\atop {\scriptstyle
i+j+k=48}} m_{i,j,k}\, a^i b^j c^k
\end{equation}
und bezeichnen dieses homogene Polynom vom Grad $48$ als das
$\La$-Gewichts\-z"ahler\-polynom von $\VM$ bez"uglich $S$, da seine 
Koeffizienten die Anzahl der linear unabh"angigen
Virasoroh"ochstgewichtsvektoren von einem bestimmten Typ z"ahlen.
Dieses Polynom h"angt tats"achlich von der Auswahl von $S$ ab. Um im n"achsten
Abschnitt den Charakter der Babymonster-\SVOA berechnen zu k"onnen, ben"otigen
wir den
\begin{satz}\label{monsterpolynom}
Bei Wahl des Systems $S$ von $48$ Transpositions-Idempotenten wie
in~\cite{DoMaZhu} erh"alt man f"ur das $\La$-Gewichtsz"ahlerpolynom von $\VM$:

{\small $P_{\VM}^S(a,b,c)= 
   {a^{48}} + {b^{48}} +
   804\,({a^{44}}\,{b^4} + {a^4}\,{b^{44}})+
   10560\,({a^{42}}\,{b^6} + {a^6}\,{b^{42}}) +
   174306\,({a^{40}}\,{b^8} + {a^8}\,{b^{40}}) +$\hfill\newline$
   1615680\,({a^{38}}\,{b^{10}} + {a^{10}}\,{b^{38}})+
   16382612\,({a^{36}}\,{b^{12}} + {a^{12}}\,{b^{36}})+
   116707584\,({a^{34}}\,{b^{14}} + {a^{14}}\,{b^{34}})+$\hfill\newline$ 
   554455407\,({a^{32}}\,{b^{16}} + {a^{16}}\,{b^{32}})+
   1786512640\,({a^{30}}\,{b^{18}} + {a^{18}}\,{b^{30}})+               
   4077522504\,({a^{28}}\,{b^{20}} + {a^{20}}\,{b^{28}})+$\hfill\newline$
   6680893824\,({a^{26}}\,{b^{22}} + {a^{22}}\,{b^{26}} )+
   7891186524\,{a^{24}}\,{b^{24}} 
      +$\hfill\newline$
\Big( 1536\,({a^{37}}\,{b^3} + {a^3}\,{b^{37}})+
   155136\,({a^{35}}\,{b^5} + {a^5}\,{b^{35}})+
   4773888\,({a^{33}}\,{b^7} + {a^7}\,{b^{33}})+$\hfill\newline$
   70699008\,({a^{31}}\,{b^9} + {a^9}\,{b^{31}})+
   596299776\,({a^{29}}\,{b^{11}} + {a^{11}}\,{b^{29}})+
   3100876800\,({a^{27}}\,{b^{13}} + {a^{13}}\,{b^{27}})+$\hfill\newline$
   10370684928\,({a^{25}}\,{b^{15}} + {a^{15}}\,{b^{25}})+
   22879881216\,({a^{23}}\,{b^{17}} + {a^{17}}\,{b^{23}})+$\hfill\newline$
   33843588096\,({a^{21}}\,{b^{19}} + {a^{19}}\,{b^{21}})
 \Big)\, {c^8} 
    + $\hfill\newline$
\Big(  16512\,({a^{30}}\,{b^2} + {a^2}\,{b^{30}})+
   1112832\,({a^{28}}\,{b^4} + {a^4}\,{b^{28}})+
   28038528\,({a^{26}}\,{b^6} + {a^6}\,{b^{26}})+$\hfill\newline$
   325307904\,({a^{24}}\,{b^8} + {a^8}\,{b^{24}} )+
   1996192896\,({a^{22}}\,{b^{10}} + {a^{10}}\,{b^{22}})+
   6985020672\,({a^{20}}\,{b^{12}} + {a^{12}}\,{b^{20}})+$\hfill\newline$
   14585195904\,({a^{18}}\,{b^{14}} + {a^{14}}\,{b^{18}})+
   18596004864\,{a^{16}}\,{b^{16}} 
  \Big)\, {c^{16}}
       + $\hfill\newline$
 \Big( 168960\,({a^{23}}\,b + a\,{b^{23}} )+
   14306304\,({a^{21}}\,{b^3} + {a^3}\,{b^{21}} )+
   300432384\,({a^{19}}\,{b^5} + {a^5}\,{b^{19}} )+$\hfill\newline$
   2446205952\,({a^{17}}\,{b^7} + {a^7}\,{b^{17}} )+
   9241528320\,({a^{15}}\,{b^9} + {a^9}\,{b^{15}} )+
   17642698752\,({a^{13}}\,{b^{11}} + {a^{11}}\,{b^{13}} )
 \Big)\, {c^{24}} 
  + $\hfill\newline$
 \Big( 9024\,({a^{16}} + {b^{16}} )+
   941568\,({a^{14}}\,{b^2} + {a^2}\,{b^{14}} )+
   14445312\,({a^{12}}\,{b^4} + {a^4}\,{b^{12}} )+
   63361536\,({a^{10}}\,{b^6} + {a^6}\,{b^{10}} )+
   102007680\,{a^8}\,{b^8} 
 \Big)\, {c^{32}}
  +
 \Big(  135168\,({a^7}\,b + a\,{b^7} )+
   946176\,({a^5}\,{b^3} + {a^3}\,{b^5} )
 \Big)\,{c^{40}} 
  + 
   2048\,{c^{48}}. 
$}
\end{satz}
{\bf Beweis:}
Nach Satz~\ref{trans-charvoa} definiert der konforme Block einer \netten
rationalen \VOA mit $k$ irreduziblen Moduln eine $k$-dimensionale Darstellung
von $\slz$. Die Darstellung $\rho: \slz \longrightarrow {\rm End}(\C^3)$ 
f"ur die \VOA $\La$ war in (\ref{ST-svoatrans}) bzw.~(\ref{smatrixsehrnett}) beschrieben worden. F"ur $\Lt$ erh"alt man die Darstellung 
$\rho^{\otimes 48}: \slz\longrightarrow {\rm End}((\C^3)^{\otimes 48})$.
Der konforme Block von $\VM$ ist wegen der Selbstdualit"at des Mondscheinmoduls~\cite{Do} und Rang $\VM=24$ der triviale eindimensionale
$\slz$-Modul. Die Restriktion der Vertexoperatorabbildung
$Y(.,z):\VM\longrightarrow{\rm End}(\VM)[[z,z^{-1}]]$ auf
$\Lt$ definiert eine $\slz$-"aquivariante Einbettung des 
konformen Blockes von $\VM$ in den von $\Lt$, d.h.~wir erhalten
ein Element $\widetilde{P}$ in $\left((\C^3)^{\otimes 48}\right)^{\slz}$.
Das Gewichtz"ahlerpolynom $P_{\VM}^S(a,b,c)$ ist gerade die Projektion
von  $\widetilde{P}$ in den Grad-$48$-Anteil der $\slz$-Invarianten
der symmetrischen Algebra ${\rm Sym}^*(\rho)$. 
Um $P_{\VM}^S(a,b,c)$ zu berechnen, mu"s man also die unter der Operation von
$\rho(\slz)$ invarianten Polynome in den $3$ Variablen $a$, $b$ und $c$ betrachten. Die Gruppe $G:=\rho(\slz)=\langle \rho(S),\rho(T)\rangle$
hat die Ordnung $1152$. Die Dimension der $G$-invarianten Polynome
ist durch den folgenden Satz von Molien (vgl.~\cite{MacSl}) gegeben:
$$\rho_G(t):=\sum_{n=0}^{\infty} \dim\left( {\rm Sym}^n(\rho)^G\right)\,t^n=
\frac{1}{|G|}\sum_{g \in G} \frac{1}{\det (1-g\,t)}. $$
Man erh"alt\footnote{Z.B.~unter Verwendung der Computeralgebrasysteme 
{\it GAP}~\cite{GAP} und {\it Mathematica}.}
\begin{eqnarray*}
\rho_G(t) & =  & 1 +{t^3} + {t^6} + {t^9} + {t^{12}} + {t^{15}} + 
   {t^{18}} + {t^{21}} + \\
& &\qquad\qquad+ 3\,{t^{24}} +  3\,{t^{27}} + 3\,{t^{30}} + 
   3\,{t^{33}} + 3\,{t^{36}} + 3\,{t^{39}} + 
   3\,{t^{42}} + 3\,{t^{45}} + 7\,{t^{48}} + \cdots .
\end{eqnarray*}
Die folgenden Polynome sind invariant unter $G$:

{\small 
$p_1={a^2}\,c - {b^2}\,c$,

$p_2=- ({a^{23}}\,b  + a\,{b^{23}})
     + ({a^{21}}\,{b^3} + {a^3}\,{b^{21}})
     + 21\,({a^{19}}\,{b^5} +{a^5}\,{b^{19}})
     - 85\,({a^{17}}\,{b^7} +{a^7}\,{b^{17}})$\hfill\newline
\qquad $ + 134\,({a^{15}}\,{b^9} +{a^9}\,{b^{15}})
     - 70\,({a^{13}}\,{b^{11}} +{a^{11}}\,{b^{13}})
  +\Big(
   - 2\,({a^{16}}\,{c^8} +{b^{16}})
   - 240\,({a^{14}}\,{b^2} +{a^2}\,{b^{14}})$\hfill\newline
\qquad $   - 3640\,({a^{12}}\,{b^4} +{a^4}\,{b^{12}}) 
   - 16016\,({a^{10}}\,{b^6} + {a^6}\,{b^{10}}) 
   - 25740\,{a^8}\,{b^8} )
   \Big)\,{c^8}$\hfill\newline
\hfill $ +\Big( 256\,({a^7}\,b +a\,{b^7})
   + 1792\,{a^5}\,{b^3} + {a^3}\,{b^5})
 \Big)\,{c^{16}} $, 

$ p_3=    3\,({a^{23}}\,b + a\,{b^{23}})
        + 253\,({a^{21}}\,{b^3} +{a^3}\,{b^{21}}) 
        + 5313\,({a^{19}}\,{b^5} +{a^5}\,{b^{19}} )
        + 43263\,({a^{17}}\,{b^7} +{a^7}\,{b^{17}})$ \hfill\newline
\hfill $	+ 163438\,({a^{15}}\,{b^9}+ {a^9}\,{b^{15}})
 	+ 312018\,({a^{13}}\,{b^{11}}+{a^{11}}\,{b^{13}})
  	- 256\,{c^{24}}$,

$ p_4=\chi_{V_{E_8}^{\otimes 3}}=\left(\frac{1}{2}\Big((a+b)^{16}+(a-b)^{16}\Big)+128\,c^{16}
\right)^3$.\phantom{x}
\footnote{Die Vielfachheit von $L_{1/2}^{\otimes 16}(\frac{1}{16})$ in 
$V_{E_8}$ ist nicht $2^8$, wie in~\cite{DoMaZhu} angegeben, sondern richtig ist 
$2^7$.}
}

Die $7$ Polynome $p_1^{16}$, $p_1^8\,p_2$, $p_1^8\,p_3$, $p_2^2$, $p_3^2$,
$p_2\,p_3$ und $p_4$ sind, wie man leicht nachrechnet, linear unabh"angig,
bilden also eine Basis f"ur ${\rm Sym}^{48}(\rho)^G$.
In~\cite{DoMaZhu}, Satz~6.5 sind die Anzahl der verschiedenen Typen von
linear unabh"angigen von Virasoroh"ochstgewichtsvektoren 
in $\VM_0$, $\VM_1$ und $\VM_2$ angeben worden (f"ur $\VM_2$ allerdings
fehlerhaft). Die (richtigen) Anzahlen sind
\begin{eqnarray*}
\VM_0 & : & m_{48,0,0}=1, \\
\VM_1 & : & m_{46,2,0}=m_{39,1,8}=m_{32,0,16}=0, \\
\VM_2 & : & m_{44,4,0}={24 \choose 2}\cdot 3-24 ,\quad m_{37,3,8}= 24
\cdot 2^6,\quad m_{30,2,16}=258\cdot 2^6 ,\hfill \\
 & & m_{23,1,24}=24\cdot2^{11}+336\cdot 2^6+24\cdot
2^{12} ,\quad m_{16,0,32}= 141 \cdot 2^6. 
\end{eqnarray*}
Leider ergeben sich hieraus nur $6$ statt $7$ linear unabh"angige Bedingungen.
F"ur $\VM_3$ erh"alt man zumindest die Gleichung
\begin{equation}\label{relationinv3}
m_{7,1,40}=9\cdot 2^{14} -6\, m_{0,0,48}.
\end{equation}
Betrachtet man zus"atzlich die Dodekaden im Golay Code
(Vektoren vom Gewicht $12$),
so liefern diese unter anderem mindestens $2\cdot 2^{10}$ linear unabh"angige 
H"ochstgewichtsvektoren vom Typ $(0,0,48)$ in $\VM_3$ und
$132 \cdot 2^{10}$ vom Typ $(7,1,40)$, d.h.~es gilt 
$m_{0,0,48}\geq 2\cdot 2^{10}$ und $m_{7,1,40}\geq 132\cdot 2^{10}$.
Wegen Relation (\ref{relationinv3}) kann in beiden Ungleichungen nur das Gleichheitszeichen gelten. 
Damit ist eine $7$-te Bedingung gefunden, man kann $P_{\VM}^S(a,b,c)$  vollst"andig berechnen und erh"alt das oben angegebene Polynom. \qed

\section{Konstruktion der Babymonster-\SVOA $\VB$ }

In diesem Abschnitt wird nun mit Hilfe von Satz~\ref{satzvoatosvoa} aus dem 
Mondscheinmodul $\VM$ die Babymonster-\SVOA $\VB$ konstruiert.
Die hierzu ben"otigten rationalen Unteralgebren $U$ und $L$ liefert
die im letzten Abschnitt betrachtete Zerlegung von $\VM$ bez"uglich
$\Lt$. Wir vermuten, da"s die \SVOA $\VB$ selbstdual, \sehrnett, unit"ar
und rational ist und die einzige derartige \SVOA ist, die keine Clifford oder
Lie Unteralgebra besitzt.

\medskip

Sei $i_{\alpha}\in \VM_2$ das Transpositions-Idempotent zu einer 
$2A$-Involution $\alpha$ im Monster. Da alle Transpositions-Idempotente
$i_{\alpha}$ unter $M$ "aquivalent sind,
liegt $i_{\alpha}$ in einem $48$-Tupel $S$ von kommutierenden
Transpositions-Idempotenten, wie es in Satz~\ref{monsterpolynom} betrachtet
wurde. Die Numerierung der Idempotenten aus $S$ sei so gew"ahlt, da"s 
$\omega_{48}=2\cdot i_{\alpha}$ das Virasoroelement des $48$-ten Faktors
$\Lacht$ der zu $S$ geh"origen Unteralgebra $\Lt\subset\VM$ ist.
Wir setzen f"ur $l=0$, $1$ und $2$ mit der gleichen Notation wie in
(\ref{monzer})
\begin{equation}\label{kl}
K_{(l)}=\bigoplus_{{\scriptstyle h_1,\ldots, h_{48}}\atop {\scriptstyle
h_{48}=n_l}} c_{h_1,\ldots,h_{48}} L_{1/2}(h_1,\ldots,h_{48}),
\end{equation}
wobei $n_0=0$, $n_1=\frac{1}{2}$ und $n_2=\frac{1}{16}$.

Die Struktur der Fusionsalgebra von $\La$ zeigt, da"s $K_{(0)}$ eine 
Unteralgebra von $\VM$ ist, und der Mondscheinmodul die direkte Summe der
drei $K_{(0)}$-Moduln $K_{(0)}$, $K_{(1)}$ und $K_{(2)}$ ist:
\begin{equation}\label{vm-zerlegung}
\VM= K_{(0)}\oplus K_{(1)} \oplus K_{(2)}.
\end{equation}
Der Modul $K_{(l)}$ wird von den Eigenvektoren der Komponente 
\hbox{$(\omega_{48})_1\in {\rm End}(\VM)$} des Vertexoperators 
$Y(\omega_{48},z)$ aufgespannt, die einen Eigenwert $\lambda\equiv n_l
 \pmod{\Z}$ besitzen.

In Verallgemeinerung von~\cite{DoMaZhu}, Prop.~5.1 (2), gilt der folgende
\begin{satz}
Jede Unter-\VOA $R$ von $\VM$, die $\Lt$ enth"alt, ist einfach.
\end{satz}
{\bf Beweis:}
Sei $T$ ein von Null verschiedener $R$-Untermodul
von $R$.
Da $\Lt\subset R$ rational ist, zerlegt sich $T$ als 
direkte Summe von $\Lt$-Moduln. F"ur Elemente $u$ und $v$ eines in $T$
enthaltenen $\Lt$-Moduls $L(h_1,\ldots,h_{48})$ gilt wegen der Fusionsregeln
$Y(u,z)v\in S[[z,z^{-1}]]$, wobei 
$S=\left(\bigoplus_{f_i\in \{0,\frac{1}{2}\}} c_{f_1,\dots, f_{48}}\,L(f_1,\ldots,f_{48})\right)\cap R$ 
Unteralgebra von $R$ ist. Da $\VM$ als adjungierter $\VM$-Modul
irreduzibel ist, gilt nach $\cite{DoLe}$, Prop.~11.9, $Y(u,z)v\not=0$ 
f"ur von Null verschiedene $u$ und $v$, d.h.~es gibt einen $\Lt$-Untermodul 
$L(f_1,\ldots,f_{48})$ mit $f_i\in \{0,\frac{1}{2}\}$ in $T$.
Wendet man das gleiche Argument nochmal auf Elemente $u'$, $v' \in
L(f_1,\ldots,f_{48})$ an, folgt $\Lt\subset T$; insbesondere also ${\bf 1}\in
T$ und daher $T=R$, d.h.~$T$ ist irreduzibel. \qed

Die \VOA $K_{(0)}$, als adjungierter Modul aufgefa"st, ist somit irreduzibel.
Dar"uber hinaus gilt der
\begin{satz}
Der $K_{(0)}$-Modul $K_{(1)}$ ist irre\-duzi\-bel.
\end{satz}
{\bf Beweis: }Der Koeffizient von $b^{48}$ im $L$-Gewichtsz"ahlerpolynom
$P^S_{\VM}(a,b,c)$ ist nach Satz 4.1.5 gerade $1$, d.h.~$L(\frac{1}{2},\ldots,
\frac{1}{2})$ ist mit einfacher Vielfachheit in $K_{(1)}$ enthalten.
Sei $t$ ein $L:=\Lt$-H"ochstgewichtsvektor f"ur $L(\frac{1}{2},\ldots,
\frac{1}{2})$.

Wir zeigen als erstes: In jedem nicht\-trivialen $K_{(0)}$-Unter\-mo\-dul $N$ 
von 
$K_{(1)}=\bigoplus c_{h_1,\ldots,h_{47},1/2} 
L_{1/2}(h_1,\ldots,h_{47},\frac{1}{2})$
ist f"ur  $c_{h_1,\ldots,h_{47},1/2}>0$ ein $L$-Mo\-dul
$L(h_1,\ldots,h_{47},\frac{1}{2})$ in der iso\-typischen $(h_1,\ldots,h_{47},\frac{1}{2})$-Kompo\-nente von $N$ ent\-halten.

F"ur $0\not=v\in L(h_1,\ldots,h_{48})\subset\VM$ ist nach $\cite{DoLe}$, Prop.~11.9 $Y(t,z)v\not=0$. Setzen wir
$$h'_i:=\cases{\frac{1}{2} & f"ur $h_i=0$, \cr 0 & f"ur $h_i=\frac{1}{2}$, \cr
\frac{1}{16} & f"ur $h_i= \frac{1}{16} $, } $$
so erhalten wir 
wegen der Fusions\-regeln~(\ref{fusion-ising})
 einen Unter\-modul $L(h'_1,\ldots,h'_{48})$ in $\VM$, 
Zu Elementen $0\not=m\in L(h_1,\ldots,h_{47},\frac{1}{2})\subset
K_{(1)}$ und $0\not=n\in L(f_1,\ldots,f_{47},\frac{1}{2})\subset
K_{(1)}$ gibt es also Elemente  $0\not=m'\in L(h'_1,\ldots,h'_{47},0)\subset
K_{(0)}$ und $0\not=n'\in L(f'_1,\ldots,f'_{47},0)\subset K_{(0)}$. Da $K_{(0)}$
irreduzibel ist, gibt es nach~\cite{DoMa-quantumgalois}, Korollar 4.2,
 ein Element $w\in K_{(0)}$
mit $w_k m'= n'$ f"ur ein $k$. Wiederum aufgrund der Fusionsregeln (!) ist dann
$0\not= Y(w,z)m\in L(f_1,\ldots,f_{47},\frac{1}{2})[[z,z^{-1}]]$, d.h.~wir
erhalten ein Element $\overline{n}\not=0$ in einem $L$-Modul 
$L(f_1,\ldots f_{47},\frac{1}{2})\subset K_{(1)}$.

Zu zeigen bleibt noch: Ist ein Mo\-dul $L(h_{1},\ldots, h_{47},\frac{1}{2})$ in
einem $K_{(0)}$-Unter\-modul $N$ ent\-halten, so auch die ganze zu\-ge\-h"orige 
iso\-typische
$(h_1,\ldots,h_{47},\frac{1}{2})$-Kom\-po\-nente. \hfill\break
Sei $X$ bzw.~$X'$ der Vektor\-raum der Virasoro\-h"ochst\-gewichts\-vek\-toren 
vom Typ 
$(h_1,\ldots,h_{47},\frac{1}{2})$ bzw.~vom Typ $(h'_1,\ldots,h'_{47},0)$. Sei 
weiter
$W\subset K_{(0)}$ der $L$-Unter\-modul, dessen 
Vertex\-operatoren die Moduln vom
Typ $(h_1,\ldots,h_{47},\frac{1}{2})$ oder --- "aquivalent dazu --- 
vom Typ $(h'_1,\ldots,h'_{47},0)$ in sich abbilden.  
Bezeichne mit $A$ die von den Endomorphismen 
$w_{{\rm deg} w -1}$ zu Virasoroh"ochstgewichtsvektoren $w\in W$ erzeugte
Unteralgebra von ${\rm End}(\VM)$.
Das Korollar 4.2 aus \cite{DoMa-quantumgalois} und Lemma 
4.5 aus \cite{DoLiMa-baby} zeigen zusammen mit der Irreduzibilit"at von $K_{(0)}$:
 $X'$ ist ein irreduzibler $A$-Modul.

Aber auch der $A$-Modul $X$ ist irreduzibel:
Die durch das oben definierte Element $t$ gegebene Abbildung
$t_{{\rm deg} t-1}: X'\longrightarrow X$ ist ein $A$-Modul\-homomorphismus,
da aus der Assoziativit"at des Tensorproduktes folgt, da"s
$t_{{\rm deg} t-1}t_{{\rm deg} t-1}$ ein nichtverschwindender Skalar ist. \qed

Ein anderer Beweis ergibt sich mit der in~\cite{DoMa-quantumgalois} 
entwickelten Theorie:\footnote{Nach einem Hinweis von C.~Dong und G.~Mason.}
Die \VOA $K_{(0)}\oplus K_{(1)}$ ist als Fixpunktmenge der $2A$-Involution
$\alpha$ einfach (s.~Satz~\ref{2A-operation}).
Betrachte auf $K_{(0)}\oplus K_{(1)}$ den Automorphismus
$\sigma$, der mit $1$ auf $K_{(0)}$ und mit $-1$ auf $K_{(1)}$ operiert.
Es folgt dann aus Satz~3 in~\cite{DoMa-quantumgalois}, 
da"s $K_{(0)}$ eine einfache \VOA und $K_{(1)}$ ein irreduzibler 
$K_{(0)}$-Modul ist.

Vermutlich ist $K_{(2)}$ ebenfalls ein irreduzibler $K_{(0)}$-Modul.

In der Zerlegung~(\ref{kl}) l"a"st sich $L_{1/2}(n_l)$ abspalten: Setzt man
\begin{equation}\label{defbabymoduln}
\VB_{(l)}= \bigoplus_{h_1,\ldots, h_{47}} c_{h_1,\ldots,h_{47},n_l} 
L_{1/2}(h_1,\ldots,h_{47}),
\end{equation}
so gilt $K_{(l)}=\VB_{(l)}\otimes L_{1/2}(n_l)$.
Wir sind nun genau in der in Abschnitt 3.1 beschriebenen Situation,
wobei hier $U=L_{1/2}^{\otimes 47}(0)$ und $L=L_{1/2}^{(48)}(0)$ ist.

Gleichung (\ref{dirtensorsumme}) aus Kapitel 3 ist die Zerlegung
$$\VM=\bigoplus_{l=0,1,2}\VB_{(l)}\otimes L(n_l)$$
des Mondscheinmoduls.

Da $K_{(l)}=\VB_{(l)}\otimes L_{1/2}(n_l)$ f"ur $l=0$ oder $1$ 
irreduzibel ist, gilt dies auch f"ur die $\VB_{(0)}$-Moduln $\VB_{(0)}$ und 
$\VB_{(1)}$ (s.~\cite{FHL}, Prop.~4.7.2).

Nach Satz~\ref{satzvoatosvoa} wird durch (\ref{defyw}) auf $\VB:=\VB_{(0)}
\oplus \VB_{(1)}$ die Struktur einer \SVOA vom Rang $23\frac{1}{2}$ definiert.
Wir werden zeigen, da"s das Babymonster $B$ --- die zweitgr"o"ste sporadische
einfache endliche Gruppe, die von B.~Fischer gefunden und in~\cite{LeSi}
konstruiert worden ist --- eine Untergruppe der Automorphismengruppe 
von $\VB$ ist, und bezeichnen daher $\VB$ als 
den {\it k"urzeren Mondscheinmodul\/} oder als die {\bf Babymonster-\SVOA{}}.

Wie in Kapitel 3 gezeigt (Lemma~\ref{svoawohldefiniert}), 
h"angt die SVOA-Struktur nicht von der gew"ahlten
rationalen Unter-VOA $L\cong L_{1/2}^{\otimes 48}(0)$ ab.

Wir formulieren Vermutung~\ref{voa2svoa} f"ur die \SVOA $\VB$ nochmal 
f"ur sich alleine:
\begin{vermutung}\label{babysd}
Die Babymonster-\SVOA $\VB$ ist selbstdual \sehrnett unit"ar und rational.
\end{vermutung}
Einige der unter der Definition \sehrnett zusammengefa"sten Eigenschaften
werden wir im folgenden auch beweisen.

\label{auto}
Bevor wir die Operation von $B$ betrachten, sei an die Definition der
Auto\-morphis\-men\-grup\-pe ${\rm Aut}(V)$ 
einer $\OVOA$ $(V,Y,{\bf 1},\omega)$ erinnert.
Ein Automorphismus von $V$ ist eine invertierbare lineare Abbildung 
$f:V\longrightarrow V$, so da"s $f({\bf 1})={\bf 1}$, $f(\omega)=
\omega$ und 
$$ f\,Y(v,z)=Y(f(v),z)f $$ 
f"ur alle $v\in V$ gilt.

Ist $(M,Y_M)$ ein $V$-Modul, so induziert ein Automorphismus $f$ auf dem 
unterliegenden linearen Raum $M$ eine neue mit  $(M,Y_M^{(f)})$ bezeichnete
$V$-Modulstruktur, die durch
$$ Y_M^{(f)}(v,z)=Y_M(f(v),z)$$
f"ur alle $v\in V$ gegeben wird.

Der Auto\-morphismus $f$ hei"st {\it inne\-rer Auto\-morphis\-mus}, falls 
f"ur je\-den irre\-du\-zib\-len $V$-Mo\-dul $(M,Y_M)$ der $V$-Mo\-dul
$(M,Y_M^{(f)})$ zu $(M,Y_M)$ iso\-morph ist.
Die Unter\-gruppe ${\rm Inn}(V)$ der inne\-ren Au\-to\-mor\-phis\-men 
von $V$ ist ein Normal\-teiler von ${\rm Aut}(V)$; die Faktor\-gruppe 
operiert als 
Permutations\-gruppe auf der Menge der irreduziblen Moduln. Es ist zu vermuten,
da"s ${\rm Aut}(V)/{\rm Inn}(V)$ eine Unter\-gruppe (oder sogar gleich) der
Automorphismen\-gruppe der abstrakten Fusions\-algebra ${\cal F} (V)$ ist.

Auf jedem irreduziblen $V$-Modul existiert eine projektive Darstellung von 
${\rm Inn}(V)$ (siehe~\cite{DoMa-orbi}, Abschnitt 2).

Wir untersuchen nun, welche Untergruppe des Monsters $M$ auf den Moduln
$K_{(l)}$ bzw.~$\VB_{(l)}$ operiert. Der Zentralisator einer $2A$-Involution ist
${\rm Cent}_M(\alpha)\cong 2.B$, die maximale (nichtsplittende) zentrale 
Erweiterung des Babymonsters $B$. Die Zerlegung von $\VM_2$ unter
$2.B$ in irreduzible invariante Unterr"aume zeigt die folgende "Ubersicht
(vgl.~\cite{Co-monster} und~\cite{MeNe}):
\begin{equation}\label{griesszerlegung}
\begin{array}{l|cccccccccc}
\hbox{$2.B$ Darstellung:} & \underline{1} & \oplus & \underline{1} & \oplus &
\underline{96255} & \oplus & \underline{4371} & \oplus & 
\underline{96256}  \\ \hline
\hbox{Eigenwerte von $(\omega_{48})_1$:\phantom{\Huge{g}}}
 & 2 & & 0 & & 0 & & \frac{1}{2} & &
\frac{1}{16} \\ \hline 
\hbox{Eigenwerte von $\alpha$:} & +1 & & +1 & & +1 & & +1 & & -1 
\end{array}
\end{equation}
Hierbei ist $\omega_{48}=2\cdot i_{\alpha}$ das Vira\-soro\-ele\-ment zum 
Idem\-po\-tent $i_{\alpha}$. Die ersten bei\-den Kom\-ponen\-ten der Zer\-legung sind die von den Virasoro\-elementen $\omega_{48}$ von $L_{1/2}^{(48)}(0)$ und
$\omega-\omega_{48}$ von
$\VB_{(0)}$ auf\-ge\-spann\-ten Vektor\-r"aume, d.h. $2.B$ fixiert die
Virasoro\-elemen\-te von $L_{1/2}^{(48)}(0)$ und $\VB_{(0)}$. 
Daher sind auch die 
Unteralgebren $L_{1/2}^{(48)}(0)$ und $\VB_{(0)}$ \hbox{$={\rm Com}_{\VM}
(L_{1/2}^{(48)}(0))=\{v\in \VM \mid ( \omega_{48})_{0}v=0\}$}
invariant unter $2.B$ (vgl.~Lemma~\ref{identifizierung}).
Auch die Zerlegung (\ref{vm-zerlegung})
wird von $2.B$ respektiert, da die Zerlegung durch die Eigenwerte
$0$, $\frac{1}{2}$ und $\frac{1}{16} \pmod{\Z}$ von $(\omega_{48})_1$
gegeben ist. 
\begin{satz}\label{2A-operation}
Auf $K_{(0)}$ und $K_{(1 )}$ operiert das zentrale Element $\alpha \in 2.B$
mit $+1$, auf $K_{(2)}$ mit $-1$.
\end{satz}
{\bf Beweis:}
Der Mondscheinmodul $\VM$ ist ein irreduzibler Modul der affinen Griess 
\hbox{Algebra ${\hat{\cal B}}$}, d.h.~er
wird erzeugt von den Koef\-fizienten der 
Vertex\-operato\-ren von Ele\-menten aus $\VM_2$~(s.~\cite{FLM},
Satz 12.3.1.~(g))
Ein Element $v\in \VM$ ist also Linear\-kombination von Elementen
$$ w=(x_1)_{n_1}(x_2)_{n_2}\ldots(x_k)_{n_k}\,{\bf 1}\qquad
\hbox{mit $x_i\in (K_{(l)})_2$.}$$
Wegen der Fusionsregeln von $\La$ gilt:
\begin{equation}\label{invol1}
w\in\cases{
 K_{(0)},\ K_{(1)} & falls eine gerade Anzahl der $x_i\in (K_{(2)})_2$,\cr
 K_{(2)} & falls eine ungerade Anzahl der $x_i\in (K_{(2)})_2$.}
\end{equation}
F"ur die Involution $\alpha\in 2.B\subset {\rm Aut}(\VM)$ erh"alt man:
\begin{eqnarray}\label{invol2}
\alpha\,w & = & (\alpha x_1)_{n_1}(\alpha x_2)_{n_2}\ldots (\alpha x_k)_{n_k}
(\alpha {\bf 1})\nonumber \\
& =& (-1)^{\# \{i \mid x_i \in (K_{(2)})_2 \} }w,
\end{eqnarray}
denn nach (\ref{griesszerlegung}) gilt
$$\alpha\, x= \cases{+x& falls $x\in  (K_{(0)})_2$, $(K_{(1)})_2$, \cr
                   -x& falls $x\in  (K_{(2)})_2$.} $$
Die beiden Gleichungen (\ref{invol1}) und (\ref{invol2}) ergeben die Behauptung.
\qed

Wir erhalten somit eine $B$-Operation auf $K_{(0)}\oplus K_{(1)}$ und, da 
$B$ auf dem Virasoroelement von $L^{(48)}_{1/2}(0)$ trivial operiert,
eine Operation auf den Faktoren des Tensorproduktes: 
\begin{lemma}
Das Babymonster $B$ operiert durch Automorphismen auf den beiden $\VB_{(0)}$-Moduln $\VB_{(0)}$ und $\VB_{(1)}$.
\end{lemma}
{\bf Beweis:}
Da die $K_{(0)}$-Moduln $K_{(l)}$ ($l=0$ oder $1$) irreduzibel sind, k"onnen
wir die "`abstrakten"' $\VB_{(0)}$-Moduln $\VB_{(l)}$ als 
"`$\VB_{(0)}$-Untermoduln"' von $K_{(l)}$ auffassen (s.~\cite{FHL},
Beweis von Satz 4.7.4): Sei $w=x\otimes y$ ein zerlegbarer von Null
verschiedener Tensor in $(K_{(l)})_{\rm min}$, der kleinsten nichttrivialen
Komponente von $K_{(l)}$. Dann ist  $\VB_{(l)}$ der "`$\VB_{(0)}$-Untermodul"',
der von den Elementen der Form $(v \otimes {\bf 1})_n(x \otimes y)$ erzeugt
wird, wobei $v\in \VB_{(0)}$ als homogen vorausgesetzt werden kann. Hierbei 
haben wir noch Proposition 4.1 aus~\cite{DoMa-quantumgalois} benutzt. 
Die $\VB_{(0)}$-Modulstruktur ist offensichtlich.

Die $B$-Operation auf $\VB_{(0)}={\rm Span}((v\otimes {\bf 1})_n( {\bf 1}
\otimes  {\bf 1}))={\rm Span}(v_n {\bf 1}\otimes  {\bf 1})
\subset K_{(0)}$ ist gerade die Einschr"ankung der 
$B$-Operation von $K_{(0)}$ auf $\VB_{(0)}={\rm Com}(L^{(48)}_{1/2}(0))$.

Die $B$-Operation auf $K_{(1)}$ bildet f"ur $g\in B$ ein Ele\-ment 
$(v\otimes {\bf 1})_n (x \otimes y)$ auf das Ele\-ment 
$g((v\otimes {\bf 1})_n (x \otimes y))
=(g(v)\otimes {\bf 1})_n (g(x \otimes y))$ ab. Nun kann we\-gen 
${\rm dim\,}(K_{(1)})_2={\rm dim\,}(\VB_{(1)})_{3/2}\cdot {\rm dim\,}(L^{(48)}_{1/2}(\frac{1}{2}))_{1/2}$
und ${\rm dim\,}(L^{(48)}_{1/2}(\frac{1}{2}))_{1/2}=1$ die Komponente
$(K_{(1)})_2$ mit $(\VB_{(1)})_{3/2}$ identifiziert werden, 
d.h.~$g(x\otimes y)=:g(x)\otimes y$
liegt wieder in $\VB_{(1)}$. Somit ist $\VB_{(1)}$ ein $B$-invarianter
"`$\VB_{(0)}$-Unter\-modul"', wir haben also eine mit der 
$\VB_{(0)}$-Modul\-struktur vertr"agliche $B$-Operation auf $\VB_{(1)}$. \qed

Fassen wir $L_{1/2}(n_l)$ als trivialen $B$-Modul auf, dann ist f"ur
$l=0$ oder $1$ die $B$-Operation auf $\VB_{(l)}$ gerade so konstruiert, da"s
$K_{(l)}\cong \VB_{(l)}\otimes L_{1/2}(n_l)$ ein Isomorphismus von
$B$-Moduln ist.
\begin{satz}
Auf der Babymonster-\SVOA $\VB=\VB_{(0)}\oplus\VB_{(1)}$ 
operiert das Babymonster nichttrivial durch Automorphismen.
\end{satz}
{\bf Beweis:}
Die Effektivit"at der $B$-Operation folgt aus der Zerlegung 
(\ref{griesszerlegung}). Zu zeigen ist, da"s die Operation mit der vollen
\SVOA-Struktur, wie sie in (\ref{defyw}) definiert wurde, vertr"aglich ist. 
Dazu betrachten wir wie beim Beweis von 
Lemma~\ref{svoawohldefiniert} oder Satz~\ref{satzvoatosvoa} die 
Korrellationsfunktionen.

Seien $u$, $v$ Elemente in $\VB_{(0)}$ oder  $\VB_{(1)}$ bzw.~$w'$ ein Element
in dem eingeschr"ankten Dualraum ${\VB_{(0)}}'$ oder  ${\VB_{(1)}}'$, 
sei $g\in B$.
Auf  ${\VB_{(l)}}'$ operiert $g$ verm"oge
$\langle g(w'),v \rangle=\langle w',g^{-1}(v) \rangle$. Die Invarianz
des Vertexoperators $\BY:\VB\longrightarrow{\rm End}(\VB)[[z,z^{-1}]]$ ist
daher "aquivalent zu der Gleichheit
\begin{equation}\label{babyinvarianz}
\langle g(w'),\BY(g(u),z)g(v) \rangle=\langle w',\BY(u,z)v \rangle
\end{equation}
von rationalen Funktionen in $z$.
Wie beim Beweis von Lemma~\ref{svoawohldefiniert} k"onnen wir
$\mid w' \mid$ $\equiv\, \mid u \mid + \mid v \mid \pmod{2}$ voraussetzen und
dann zu $u$, $v$ und $w'$ Elemente $\ou$, $\ov\in L_{1/2}(n_l)$ 
bzw.~$\ow' \in L'_{1/2}(n_l) $ mit
$\mid u \mid =\mid \ou \mid$, $\mid v \mid =\mid \ov \mid$ und
$\mid w' \mid =\mid \ow' \mid$ sowie
\begin{equation}\label{ynonzero}
\langle \ow',\LVirY(\ou,z)\ov \rangle \not= 0
\end{equation}
finden. Man erh"alt
\begin{eqnarray*}
& &\langle w',\BY(u,z)v\rangle\cdot
   \langle \ow',\LVirY(\ou,z)\ov\rangle\hfill \\
&=&\langle w'\otimes\ow',(\BY(u,z)\otimes\LVirY(\ou,z))
v\otimes\ov\rangle   \\
&=&\langle w'\otimes\ow',\MY(u\otimes\ou,z)v\otimes\ov\rangle             \\
\noalign{\hbox{wegen der Invarianz von $\MY$ unter $2.B\subset M$}\vspace{-5mm}}  \\
&=&\langle g(w'\otimes\ow'),\MY(g(u\otimes\ou),z)g(v\otimes\ov)\rangle   \\
\noalign{\hbox{nach Definition der $B$-Operation auf $\VB_{(l)}\otimes 
L_{1/2}(n_l)$}\vspace{-5mm}}         \\
&=&\langle g(w')\otimes\ow',\MY(g(u)\otimes\ou,z)g(v)\otimes\ov\rangle   \\
&=&\langle g(w'),\BY(g(u),z)g(v)\rangle\cdot 
\langle \ow',\LVirY(\ou,z)\ov\rangle. \\
\end{eqnarray*}
Hieraus und aus (\ref{ynonzero}) folgt (\ref{babyinvarianz}).
\qed

Der Beweis l"a"st sich 
auf alle der in Kapitel 3 aus einem Paar $(V,L)$
konstruierten \SVOAs $W$ "ubertragen, falls die Gruppenoperation 
auf $W$ geeignet definiert wird.

J.~Tits hatte in~\cite{ti-monster} u.a.~gezeigt, da"s das Monster die volle 
Automorphismengruppe der Griessalgebra ${\cal B}$
ist. Dazu wird zuerst gezeigt, da"s 
die Automorphismengruppe endlich ist, und dann wird eine gruppentheoretische
Charakterisierung des Monsters durch S.~Smith verwendet. Da nach~\cite{FLM}
das Monster auf $\VM$ operiert, mu"s es wegen der Irreduzibilit"at von
$\VM$ als ${\hat{\cal B}}$-Modul
die volle Automorphismengruppe von 
$\VM$ sein. 

F"ur $\VB$ k"onnte man "ahnlich vorgehen, und versuchen die Arbeit~\cite{Bier} 
zu verwenden, in der das Babymonster dadurch charakterisiert wird, da"s der
Zentralisator einer Involution eine Erweiterung der Conwaygruppe $Co_2$
durch eine extraspezielle $2$-Gruppe ist.
\begin{satz}
Die Automorphismengruppe von $\VB_2$ der $96256$-dimensionalen Unteralgebra
der Griess Algebra ist eine endliche Gruppe $G>B$.
F"ur die Automorphismengruppe von \SVOA $\VB$ gilt 
$2\times B < {\rm Aut}(\VB) < 2 \times G$, wobei der Faktor $2$ dem in 
Seite~\pageref{twistsektor} betrachteten Automorphismus $\kappa$ entspricht, 
der auf $\VB_{(1)}$ mit $-1$ operiert.
\end{satz}
{\bf Beweis:} (Skizze)\hfill\break
(1) $G:={\rm Aut}(\VB_2)$ ist endlich: Die in $\VB_2$ gelegenen Transpositions-Idempotenten erzeugen $\VB_2$
als Vektorraum, da $B$ transitiv auf ihnen operiert. W"are $G$ unendlich,
so g"abe es auf der $96255$-dimensionalen Sph"are
ein Transpositions-Idempotent $i$ mit sehr nahem $G$-Bild $i'$:
$$i'=i+\epsilon\, u + {\rm O}(\epsilon^2),$$
wobei $u$ orthogonal zu $i$ und $\epsilon$ sehr klein.
Wir erhalten
$$i'_1i'=i_1i+2\epsilon \,i_1 u+  {\rm O}(\epsilon^2).$$
Nach~(\ref{griesszerlegung}) gilt aber $|| i_1u|| \leq \frac{1}{4} || u||$,
im Widerspruch zur Gleichung $i'_1i'=i'$.

(2) $B\subset {\rm Aut}(\VB_{(0)})\subset G$: Hierzu ist zu zeigen, da"s $\VB_{(0)}$ ein irreduzibler Modul "uber der affinen Algebra $\hat{VB}^{\VBE}_2$ ist.
Dazu verwende man die Zerlegung als $L_{1/2}^{\otimes 47}(0)$-Modul und
die Zerlegung als $B$-Modul, wie in Lemma~\ref{charzer} beschrieben.

(3) Sei $\alpha\in  {\rm Aut}(\VB)$, $\alpha=\alpha_0\oplus\alpha_1$.
Dann ist auch $\overline{\alpha}=\alpha_0\oplus -\alpha_1\in {\rm Aut}(\VB)$.
Umgekehrt ist $\alpha_1$ durch $\alpha_0$ bis auf das Vorzeichen
festgelegt:
Sei   $\alpha'=\alpha_0\oplus \alpha_1'$; dann gilt $\alpha'\alpha^{-1}=
{\rm id}\oplus \alpha_1'\alpha_1^{-1}$, d.h.~$\alpha_1'\alpha_1^{-1}$ ist ein
$\VB_{(1)}$-Modulhomomorphismus. Wegen der Irreduzibilit"at von $\VB_{(1)}$ 
ist er ein Skalar $s$ 
und somit $s=\pm 1$. (Sei $\alpha={\rm id}\oplus s\cdot {\rm id}$, f"ur $x\in 
\VB_{(1)}$ gilt $s^2\cdot x_nx=(sx)_n(sx)=x_nx$, also  $s=\pm 1$.)\qed

Vermutlich gilt die Gleichheit $G=B$. Ein anderer Beweis f"ur 
${\rm Aut}(\VB)=2\times B$ mit Hilfe von \VOA-Theorie erg"abe sich aus 
Vermutung~\ref{babysd} und der zu vermutenden (projektiven) Erweiter\-barkeit
von Auto\-morphismen einer \netten rationalen \VOA zur vollen
"`Inter\-twiner\-algebra"' (d.h.~den Inter\-twiner\-operatoren zwischen den
irreduziblen Moduln; siehe auch S.~\pageref{auto}): Jeder Auto\-morphismus
von $\VB_{(0)}$ lieferte einen Auto\-morphismus von
$\VM=\VB_{(0)}\otimes\La\oplus\VB_{(1)}\otimes\Lb\oplus\VB_{(2)}\otimes\Lc$
der $\La$ fixiert, d.h.~er m"u"ste wegen ${\rm Aut}(\VM)=M$ in $2.B$ liegen.

\begin{satz}[Charakter von $\VB$]\label{babychar}
F"ur den Charakter der Babymonster-\SVOA $\VB$ gilt
$$\chi_{\VB}=\chi_{1/2}^{47}-47\cdot\chi_{1/2}^{23}=
q^{-\frac{47}{48}}(1+ 4371\, q^{\frac{3}{2}}+ 96256\,q^2+
1143745\, q^{5/2}+ 9646891\, q^3+\cdots ).$$
Im einzelnen erh"alt man f"ur die Charaktere der 
$\VB_{(0)}$-Moduln:
\begin{eqnarray}\label{babychargleich}
\chi_{\VB_{(0)}} & = &  {q^{-{{47}\over {48}}}}\,
(1 + 96256\,{q^2} + 9646891\,{q^3} + 366845011\,{q^4} + 8223700027\,{q^5} + 
   \cdots)
, \nonumber \\
\chi_{\VB_{(1)}} & = & {q^{-{{47}\over {48}}}}\,
( 4371\,{q^{{3\over 2}}} + 1143745\,{q^{{5\over 2}}} + 
   64680601\,{q^{{7\over 2}}} + 1829005611\,{q^{{9\over 2}}} + 
 \cdots)
,  \\
\chi_{\VB_{(2)}} & = & 
{q^{-{{47}\over {48}}}}\,
(96256\,{q^{{{31}\over {16}}}} + 10602496\,{q^{{{47}\over {16}}}} + 
   420831232\,{q^{{{63}\over {16}}}} + 9685952512\,{q^{{{79}\over {16}}}} + 
     \cdots).  \nonumber
\end{eqnarray}
\end{satz}
Wenn wir Vermutung \ref{babysd} schon bewiesen h"atten,
erg"abe sich der Charakter von $\VB$ wegen $\VB_{1/2}=\VB_{1}=0$ unmittelbar 
aus Satz~\ref{sdsvoachar}. Stattdessen gehen wir direkter vor und verwenden 
Satz~\ref{monsterpolynom}.

{\bf Beweis von Satz~\ref{babychar}:} 
Fixiert man ein Idempotent $a_x$, $x\in\{1,\dots,48\}$ aus dem in 
Satz~\ref{monsterpolynom} verwendeten System $S$ von 
Transpositions-Idempotenten,
so entspricht $S_x:=S\setminus\{a_x\}$ einem System von $47$ 
Idempotenten f"ur die mittels des Idempotents
$a_x$ wie in (\ref{defbabymoduln}) konstruierte Version von $\VB$.
Wir definieren analog zu (\ref{mijk}) und (\ref{pvsabc}) die 
$\La$-Gewichtsz"ahlerpolynome
von $\VB_{(0)}$, $\VB_{(1)}$ und $\VB_{(2)}$ bzgl.~$S_x$ durch
$$P_{\VB_{(l)}}^{S_x}=\sum_{{\scriptstyle i,j,k\in\Z_+}\atop {\scriptstyle
i+j+k=47}} m_{i,j,k}^{x,l}\, a^ib^jc^k,$$
wobei $m_{i,j,k}^{x,l}$ die Anzahl von $\Lt$-Moduln in der Zerlegung
(\ref{monzer}) ist, f"ur die $h_x=0$ (falls $l=0$), $h_x=\frac{1}{2}$
 (falls $l=1$), bzw.~$h_x=\frac{1}{16}$ (falls $l=2$) ist.
Vermutlich h"angen die Koeffizienten $m_{i,j,k}^{x,l}$ von der Position
$x$ ab, so da"s wir sie nicht direkt ausrechnen k"onnen.  Durch Einsetzen 
der Charaktere der $\La$-Moduln $a=\chi_{\La}$, $b=\chi_{\Lb}$ und $c=\chi_{\Lc}$ in $P_{\VB_{(l)}}^{S_x}(a,b,c)$
erh"alt man die Charaktere  $\chi_{\VB_{(l)}}$.
Diese sind sicher unabh"angig von der Auswahl von $a_x$, da das Monster ja
transitiv auf den Idempotenten zu $2A$-Involutionen operiert. Daher 
erh"alt man
\begin{equation}\label{mittel}
\chi_{\VB_{(l)}}=\frac{1}{48}\sum_{x\in\{1,\ldots,48\}}P_{\VB_{(l)}}^{S_x}
(\chi_{\La},\chi_{\Lb},\chi_{\Lc}).
\end{equation}
Die Koeffizienten $w_{i,j,k}$ von 
$\sum_{x\in\{1,\ldots,48\}}P_{\VB_{(l)}}^{S_x}(a,b,c)$ ergeben sich aber aus
denen von $P_{\VM}^{S}(a,b,c)$ durch doppeltes Abz"ahlen:
$$w_{i,j,k}=
\cases{ 
(i+1)\,m_{i+1,j,k}, & falls $l=0$, \cr
(j+1)\,m_{i,j+1,k}, & falls $l=1$, \cr
(k+1)\,m_{i,j,k+1}, & falls $l=2$. \cr }$$
Die Gleichungen~(\ref{babychargleich}) folgen dann aus (\ref{mittel}) durch
Einsetzen und Feststellen der "Ubereinstimmung der ersten Koeffizienten der 
$q$-Entwicklung in $\Z[[q^{1/48}]]$, denn beide Seiten sind
meromorphe Funktionen auf einer Fl"ache $\overline{\H/\Gamma}$
mit endlich vielen Polstellen beschr"ankter
Ordnung.\ \qed

F"ur den Mondscheinmodul $\VM$ wird in~\cite{FLM} vermutet, da"s er die einzige selbstduale rationale \VOA $V$ vom Rang $24$ mit $V_1=0$ ist. In Analogie zu 
der Charakterisierung des k"urzeren Golaycodes und des k"urzeren Leechgitters
sollte folgendes gelten:
\begin{vermutung}\label{vermbabyext}
Die Babymonster-\SVOA $\VB$ ist bis auf Isomorphie die einzige selbstduale, 
unit"are, \sehrnette und rationale \SVOA $V$ vom Rang $23\frac{1}{2}$ mit
$V_{1/2}=V_1=0$.
\end{vermutung}
Dieses Resultat w"urde aus der Eindeutigkeit des Mondscheinmoduls $\VM$ folgen,
wenn gezeigt wird, da"s die Transpositions-Idempotenten die einzigen 
Idempotenten der Griess Algebra mit Norm $\frac{1}{16}$ sind und 
Vermutung~\ref{umkehr} gilt. 
Die Eindeutigkeit der Norm-$\frac{1}{16}$-Idempotenten sollte sich aus der Bemerkung in~\cite{Mia} ergeben, da"s n"amlich
ein solches Idempotent aufgrund der Struktur der Fusionsalgebra von $L_{1/2}(0)$
stets eine nichttriviale Involution in ${\rm Aut}(\VM)=M$ liefert. 
Um umgekehrt die Eindeutigkeit des 
Mondscheinmoduls aus der Eindeutigkeit
von $\VB$ herzuleiten, m"u"ste man zus"atzlich zeigen,
da"s die Algebra $V_2$ einer extremalen \VOA vom Rang $24$ (d.h.~einer \VOA
mit $V_1=0$, s.~Kap.~5) stets Idempotente der Norm $\frac{1}{16}$ besitzt.

\medskip
Die Zerlegung der homogenen Komponenten von $\VM$ in irreduzible $M$-Moduln
(s.~\cite{CoNo,Bo-lie}) liefert durch Restriktion auch eine Zerlegung in
irreduzible $2.B$-Moduln. Der Teil, auf dem das zentrale Element von $2.B$ mit
$-1$ operiert, bestimmt den $K_{(2)}$-Anteil von $\VM$ und man erh"alt so
die Zerlegung der Komponenten von $\VB{(2)}$ in irreduzible $2.B$-Moduln.
F"ur die Aufteilung der irreduziblen $B$-Moduln von $K_{(0)}\oplus K_{(1)}$
auf die Summanden $K_{(0)}$ bzw.~$K_{(1)}$ ben"otigt man weitere Informationen,
die sich z.B.~aus den verallgemeinerten Nortonvermutungen~\cite{no-general} und
einer expliziten Beschreibung des $2A$-getwisteten Sektors $\VM(2A)$ ergeben
w"urden (s.~\cite{DoLiMa-baby}). F"ur die ersten Komponenten von $\VB$ erh"alt
man die Aufteilung allerdings eindeutig aus dem in Satz~\ref{babychar} 
berechneten Charakter und der Zerlegung als Virasoromodul (vgl.~\cite{DoLiMa}).

\begin{lemma}\label{charzer}
Die beiden nachfolgenden Tabellen beschreiben die Zerlegung der er\-sten Kompo\-nenten von
$\VB$ bzw. $\VB_{(2)}$ in irreduzible $B$- 
bzw.~$2B$-Moduln. Die Eintr"age geben die jeweilige Vielfachheit des 
Charakters an.

\medskip
$
\begin{array}{l|rrrrrrrrrrrrrr}
 &\chi_1&\chi_2&\chi_3&\chi_4&\chi_5&\chi_6&\chi_7&\chi_8 &
\chi_9&\chi_{10}&\chi_{11}&\chi_{12}&\chi_{13}&\chi_{14}
\\ \hline
\VB_0 & 1 \\
\VB_1 & 0 \\
\VB_2 & 1 & 0 & 1 \\
\VB_3 & 1 & 0 & 1 & 0 & 0 & 1\\
\VB_4 & 2 & 0 & 2 & 0 & 1 & 1 & 0 & 1\\
\VB_5 & 2 & 0 & 3 & 0 & 1 & 3 & 0 & 1 & 0 & 0 & 1 & 0 & 0 &1\\
%
\\ \hline  
\VB_{1/2} & 0 \\
\VB_{3/2} & 0 & 1 \\
\VB_{5/2} & 0 & 1 & 0 & 1 \\
\VB_{7/2} & 0 & 2 & 0 & 1 & 0 & 0 & 1\\
\VB_{9/2} & 0 & 3 & 0 & 2 & 0 & 0 & 1 & 0 & 1 &1 \\
\end{array}
$

\medskip
$
\begin{array}{l|rrrrrrrrrrrrrrrr}
 &\chi_{185}&\chi_{186}&\chi_{187}&\chi_{188}  &\chi_{189}&\chi_{190}&\chi_{191}&\chi_{192}  &\chi_{193}&\chi_{194}&\chi_{195}&\chi_{196} 
\\ \hline
(\VB_{(2)})_{0} & 1 \\
(\VB_{(2)})_{1} & 1 & 1 \\
(\VB_{(2)})_{2} & 2 & 1 & 1\\
(\VB_{(2)})_{3} & 3 & 2 & 2 & 1 \\
(\VB_{(2)})_{4} & 5 & 4 & 4 & 2 & 1 & 0 & 0 &1 \\
(\VB_{(2)})_{5} & 8 & 7 & 8 & 5 & 2 & 0 & 0 &2 & 1 & 1 \\
(\VB_{(2)})_{6} & 13 & 12 & 15 & 10 & 5 & 0 & 0 &5 & 3 & 3 & 1 &1 \\

\end{array}
$
\end{lemma}\qed

In den n"achsten nicht mehr in der Tabelle aufgef"uhrten Komponenten von
$\VB_{(0)}$ und $\VB_{(1)}$ finden sich erstmals gleiche $B$-Charaktere. 
Die Charaktere $\chi_{190}$ und $\chi_{191}$ k"onnen, da nicht reell, 
in $\VB_{(2)}$ nicht vorkommen.

\chapter{Extremale selbstduale Vertexoperator-Superalgebren}

Auch dieses Kapitel ist durch die Analogie zwischen Codes, Gittern und \VOAs
moti\-viert. Ausgehend von dem Problem m"oglichst dichte Kugel\-packungen
in den $2$-fach homogenen R"aumen $\F_2^n$ (Hammingschema, s.~\cite{Del})
bzw.~$\R^n$ (euklidischer Raum) zu konstruieren, betrachtet man in
der Codierungstheorie lineare Codes mit m"oglichst gro"sem Minimalgewicht 
bzw.~Gitter mit gro"ser Minimalnorm (vgl.~\cite{CoSl}, Kap.~9).
Es zeigt sich, da"s in der Klasse der selbstdualen (geraden) Codes 
bzw.~Gitter, asymptotisch gute Packungen existieren, die den theoretischen 
Schranken sehr nahe kommen 
(vgl.~\cite{CoSl}, Kap.~7 u.~\cite{MiHu}, Kap.~II, \S 9). F"ur
kleine $n$ ($n\,{{\normalsize <}\atop{\normalsize \sim}}\, 24$) sind diese 
Codes und Gitter h"aufig die eindeutig bestimmte optimale L"osung des 
Packungsproblemes und weisen interessante Symmetriegruppen auf.

\medskip

Wir definieren das Minimalgewicht einer \OVOA als das minimale 
konforme Gewicht eines
vom Vakuum verschiedenen Virasoro\-h"ochstgewichts\-vektors.
V"ollig analog zu der Definition von
selbstdualen Codes und Gittern (vgl.~\cite{CoSl}, Kap.~7),
lassen sich nun extremale selbstduale \OVOAs definieren und "ahnliche
Aussagen beweisen; insbesondere wird eine Charakterisierung der extremalen 
selbstdualen \SVOAs (Satz~\ref{extsvoasatz}) gegeben. 

Offen bleibt die Frage nach der Existenz eines den 
Ku\-gel\-packungs\-prob\-lemen
ver\-wand\-ten "`geo\-metri\-schen Extremal\-problems"', f"ur wel\-ches diese
extremalen \OVOAs die (eindeutigen ?) L"osungen darstellen. Gibt es eine
"`nichtmeromorphe"' verallgemeinerte \VOA ${\cal U}_c$ vom Rang $c$, in die
sich alle \OVOAs gleichen Ranges einbetten lassen ? Sie w"are das Analogon
zu  $\F_2^n$ bzw.~$\R^n$. Was ist das \VOA-Analogon zu Blockpl"anen~\cite{BJL}
bzw.~sph"arischen Blockpl"anen~\cite{DGS} ?

\section{Das Minimalgewicht einer \OVOA{ }}

Das Minimalgewicht eines linearen Codes ist das kleinste Hamminggewicht
von allen Codew"ortern au"ser dem Nullvektor. Die Minimalnorm eines Gitters 
ist die kleinste Quadratl"ange von allen Gittervektoren au"ser dem Ursprung.
Das richtige Analogon bei \VOAs ist das kleinste konforme Gewicht eines 
Virasoro\-h"ochstgewichts\-vektors au"ser dem Vakuum.

Eine \OVOA $V$ ist ein Modul "uber der von dem Virasoro\-element
erzeugten Virsoro\-algebra. 
\begin{definition}\label{defmingew} Das Minimalgewicht $\mu(V)$ einer \SVOA $V$
ist der kleinste $L_0$-Eigenwert eines Virasoroh"ochstgewichtsvektors
aus $V$ ungleich dem Vakuum $\C \cdot {\bf 1} \subset V_0$. 
Falls au"ser dem Vakuum keine weiteren Virasoroh"ochstgewichtsvektoren 
existieren, sei $\mu(V)=\infty$ gesetzt.
\end{definition}
{\it Anmerkung:\/} W"urde man sich in der Definition nicht auf 
Virasoroh"ochstgewichtsvektoren beschr"anken, sondern beliebige
$L_0$-Eigenwerte zulassen, so w"are wegen $\omega\in V_2$ stets
$\mu(V)\leq 2$. 

Sei nun $V$ als \nett vorausgesetzt. Nach Definition~\ref{defnett} 
existiert dann eine Zerlegung von $V$ als direkte Summe von 
H"ochstgewichtsdarstellungen ihrer Virasoroalgebra:
$$ V=\bigoplus_i M_{c,h_i}. $$
Da dort weiter $\dim(V_n)=0$ f"ur $n<0$ und $\dim(V_0)=1$ vorausgesetzt werden, 
ist das Minimalgewicht $\mu(V)$ f"ur \VOAs eine positive ganze Zahl bzw.~eine
positive halbganze Zahl bei \SVOAs.

F"ur R"ange $c>1$ ist der vom Vakuum erzeugte Modul $M_c:={\cal U}
(\hbox{Vir}^{-})/ \langle L_{-1} {\bf 1} \rangle $ 
irreduzibel und hat den Charakter 
$\chi_{M_c}=q^{-c/24}\prod_{n=2}^{\infty}\frac{1}{1-q^n}$ (vgl.~\cite{Wan}).
Die Moduln $M_{c,h}$, \hbox{$c>1$}, \hbox{$h>0$} besitzen den Charakter
$\chi_{M_c}=q^{-c/24}\prod_{n=1}^{\infty}\frac{1}{1-q^n}$.
Insgesamt schreibt sich der Charakter von $V$ daher als
\begin{equation}\label{partgeneral}
\chi_V=q^{-\frac{c}{24}}\left(\prod_{n=2}^{\infty}\frac{1}{1-q^n}+
\prod_{n=1}^{\infty}\frac{1}{1-q^n}\left(\sum_{i\geq\mu(V)}P_i\cdot q^i
\right)\right),
\end{equation}
wobei $P_i$ die Dimension des Vektorraumes der Virasoroh"ochstgewichte
vom Gewicht $(c,i)$ ist,
d.h.~bis zum Term $q^{-c/24+\mu(V)}$ stimmen die Charaktere von $V$ 
und $M_c$ "uberein: $\chi_V=\chi_{M_c}+q^{-c/24}\cdot O(q^{\mu(V)})$.

Im Gegensatz zu der Situation bei Codes und Gittern, wo man durch 
Um\-skalie\-ren
beliebig gro"ses Minimal\-gewicht er\-h"alt, ist es bei \VOAs f"ur
\hbox{$c>1$} nicht so einfach,
uni\-t"are ratio\-nale Bei\-spiele mit gro"sem Mini\-mal\-gewicht zu fin\-den.
Die zu Virasoro\-h"ochst\-gewichts\-darstel\-lungen asso\-ziierten \VOAs
(s.~\ref{satzminimalemodelle}) 
sind nur f"ur $c=1-\frac{6}{n(n+1)}$, $n=3$, $4$, $\dots$
rational. Die zu H"ochstgewichtsdarstellungen von Kac-Moody Algebren
assoziierten \VOAs (s.~\ref{satzlievoas}) und die \VOAs zu Gittern 
(s.~\ref{satzgittersvoa})
haben wegen $V_1\not=0$ das Minimalgewicht $1$.
Die $\Z_2$-Orbifolds von Gitter-\VOAs haben h"ochstens das Minimalgewicht $2$
(vgl.~Gleichung (\ref{extremalorbi})). Wenn, wie vermutet
(vgl.~\cite{orbicft}), 
Orbifoldkonstruktionen von rationalen \VOAs wieder rational sind,
so haben wir folgendes Beispiel mit $c=24$ und $\mu(V)\geq 10$:
Wir betrachten die Unter-\VOA $(\VM)^M$ von $\VM$, die aus den unter der 
Operation des
Monsters invarianten Vektoren des Mondscheinmoduls besteht.
Die Zerlegung der "`Head-Charaktere"' in irreduzible \hbox{$M$-Darstellungen}
(s.~\cite{CoNo}, \cite{atlas}, S.~321) zeigt, da"s bis einschlie"slich
$\VM_{10}$ nur Vektoren im Vakuummodul $M_{24}$ invariant unter $M$ sind,
d.h.~es gilt $\mu((\VM)^M)\geq 10$. Tats"achlich ist  $\mu((\VM)^M)=12$,
wie man z.B.~aus der Tabelle in~\cite{HaLa} ablesen kann.

Wenn die Menge aller \netten unit"aren rationalen \VOAs "ahnliche Eigenschaften 
wie die der Codes und Gitter hat, ist folgendes zu erwarten: F"ur festen
Rang $c$ gibt es \VOAs mit beliebig gro"sem Minimalgewicht, die zugeh"orige
Fusionsalgebra wird dann aber auch beliebig "`gro"s"'. Genauer sollte  
f"ur {\it selbstduale} \VOAs die folgende Situation vorliegen: Bezeichne mit
$$\mu_c:=\sup_V \mu(V) $$
--- wobei $V$ alle selbstdualen \VOAs vom Rang $c$ durchl"auft ---
 das gr"o"ste Minimalgewicht einer selbstdualen \VOAs vom Rang $c$.
Zu erwarten ist, da"s positive Konstanten $C_1$ und $C_2$ existieren, so
da"s
$$C_1\,\leq\,\liminf_{c\longrightarrow\infty}\,\frac{\mu_c}{c}\,\leq\,
\limsup_{c\longrightarrow\infty}\,  \frac{\mu_c}{c}\,\leq\, C_2. $$
Im n"achsten Abschnitt werden wir zeigen, da"s zumindest die obere Absch"atzung
mit $C_2=\frac{1}{24}$ erf"ullt ist.

\section{Extremale \sde \VOAs{ }}

Codes und Gitter mit gro"sem Minimalgewicht finden sich unter den selbstdualen
Codes bzw.~Gittern (siehe \cite{CoSl}, Kap.~7.3--7.7). Insbesondere die 
extremalen
geraden (oder Typ II) selbstdualen Codes und Gitter besitzen interessante
Eigenschaften. In diesem Abschnitt betrachten wir das Analogon f"ur selbstduale \VOAs.

Sei in diesem Abschnitt $V$ eine \nette rationale selbstduale \VOA vom
Rang $c$. Nach Satz~\ref{charsdvoa} ist $c$ ein ganzzahliges Vielfaches 
von 8, und der Charakter von $V$ ist ein Polynom in
$\chi_8=\sqrt[3]{j}=q^{-1/3}(1+248\,q+4124\,q^2+\cdots)$:
\begin{equation}\label{eins}
\chi_V=\sum_{r=0}^k a_r\cdot\chi_8^{c/8-3r};\quad k=\left[\frac{c}{24}\right],
\end{equation}
mit eindeutig bestimmten ganzen Zahlen $a_0$, $\dots$, $a_k$.\pagebreak[2]
\begin{definition}
Sind die $a_r$ so gew"ahlt, da"s
\begin{equation}\label{extvoachar}
\chi_V=\chi_{M_c}\cdot(1+A_{k+1}q^{k+1}+A_{k+2}q^{k+2}+\cdots),
\end{equation}
so hei"st~(\ref{extvoachar}) der extremale Charakter in Rang $c$ und eine
\sde \VOA mit diesem Charakter eine extremale \sde \VOA.
\end{definition}

Der Charakter einer extremalen \VOA ist insbesondere eindeutig bestimmt.
Wegen Gleichung~(\ref{partgeneral}) gilt f"ur das Minimalgewicht einer extremalen
\VOA $\mu(V)\geq[\frac{c}{24}]+1$.

"Ahnlich wie in~\cite{MaOdSl} beweist man den
\begin{satz}\label{extvoasatz}
In dem extremalen Charakter~(\ref{extvoachar}) sind f"ur alle R"ange $c$
der Koeffizient $A_{k+1}$ und die Differenz $A_{k+2}-A_{k+1}$ positiv.
\end{satz}
Mit~(\ref{partgeneral}) erhalten wir aus~(\ref{extvoachar}) f"ur die 
Di\-men\-sion $P_i$ des Rau\-mes
der Virasoro\-h"ochst\-gewichts\-vektoren vom Gewicht $i$: 
$P_0=1$, $P_1=\ldots=
P_k=0$, $P_{k+1}=A_{k+1}$, $P_{k+2}=A_{k+2}-A_{k+1}$.

Wegen $P_{k+1}=A_{k+1}\not=0$ folgt unmittelbar das
\begin{korollar}\label{korminVOA}
Das Minimalgewicht einer selbstdualen VOA $V$ erf"ullt
\begin{equation}\label{voamin}
\mu(V)\leq\left[\frac{c}{24}\right]+1.
\end{equation}
\end{korollar}

Anders als im Fall von Codes und Gittern k"onen wir aus Satz~\ref{extvoasatz}
{\it nicht\/} schlie"sen,
da"s f"ur gro"se $c$ keine extremalen selbstdualen \VOAs existieren, denn
$P_{k+2}=A_{k+2}-A_{k+1}$ wird f"ur gro"se $c$ nicht negativ.

{\bf Beweis von Satz~\ref{extvoasatz}:}
Wir entwickeln $\chi_{M_c}\cdot\chi_8^{-c/8}$ in Potenzen von
$\phi:=\chi_8^{-3}=j^{-1}=q-744\,q^2+O(q^3)$:
\begin{equation}\label{drei}
\chi_{M_c}\cdot\chi_8^{-c/8}=\sum_{r=0}^{\infty}\alpha_r\cdot\phi^r.
\end{equation}
Hierbei gilt f"ur den Koeffizienten $\alpha_r$ nach dem B"urmanschen
Satz~\cite{WiWo} die Formel
\begin{equation}\label{vier}
\alpha_r=\frac{1}{r!}\frac{d^{r-1}}{d q^{r-1}}\left\{
\frac{d(\chi_{M_c}\cdot\chi_8^{-c/8})}{dq}\left(\frac{q}{\phi}\right)^r
\right\}_{q=0}.
\end{equation}
Unter Verwendung von~(\ref{eins}), (\ref{extvoachar}) und (\ref{drei})
erh"alt man
\begin{eqnarray}
\sum_{r=0}^{k}a_r\phi^r & = & \chi_V\cdot\chi_8^{-c/8}
=\chi_V\cdot\chi^{-1}_{M_c}\cdot\chi_{M_c}\cdot\chi_8^{-c/8}\nonumber\\
 & = & \left(1+\sum_{n=k+1}^{\infty}A_n q^n\right)
 \left(\sum_{r=0}^k\alpha_r\phi^r+\sum_{r=k+1}^{\infty}\alpha_r\phi^r\right).
\end{eqnarray}\pagebreak
Koeffizientenvergleich liefert $a_r=\alpha_r$ f"ur $0\leq r \leq k$, 
und aus den Koeffizienten von $q^{k+1}$ und $q^{k+2}$ ergibt sich
\begin{eqnarray*}
A_{k+1} & = & -\alpha_{k+1}, \\
A_{k+2} & = & -\alpha_{k+2}+ 744(k+1)\alpha_{k+1}-\alpha_1\cdot A_{k+1}.
\end{eqnarray*}
F"ur $r=1$ liefert Gleichung~(\ref{vier}) $\alpha_1=-248\cdot\frac{c}{8}$.
Somit ist
$$A_{k+2}=
-\alpha_{k+2}+744(1+\left[\frac{c}{24}\right]-\frac{c}{24})\alpha_{k+1}.$$
Gleichung~(\ref{vier}) f"ur die $\alpha_r$ l"a"st sich 
unter Verwendung der Abk"urzung $\phantom{}\frac{d}{dq}=\phantom{}'$
weiter umformen zu
$$\alpha_r=\frac{1}{r!}\frac{d^{r-1}}{d q^{r-1}}\left\{
q^r\left[\chi'_{M_c}\cdot\chi_8-\frac{c}{8}\cdot
\chi_{M_c}\cdot\chi'_8\right]\cdot\chi_8^{3r-c/8-1}
\right\}_{q=0}. $$
F"ur $r\geq k+1$ ist $3r-c/8-1\geq3([\frac{c}{24}]+1)-c/8-1\geq 0$,
und daher sind die Koeffizienten von
$\chi_8^{3r-c/8-1}$ alle positiv.
Da $\chi_8$ der Charakter der \VOA $V_{E_8}$ ist, folgt aus (\ref{partgeneral}),
da"s $\frac{\chi_8}{\chi_{M_c}}=q^{(c-8)/24}\left(1+\frac{1}{1-q}
\left(\sum_{i\geq 1}P_i\,q^i\right)\right)$ positive Koeffizienten hat. Dann
hat aber auch $\chi_{M_c}^2\left(\frac{\chi_8}{\chi_{M_c}}\right)'=\chi_{M_c}
\chi_8'-\chi_{M_c}'\chi_8$ positive Koeffizienten und  
$\chi_{M_c}' \chi_8-\frac{c}{8}\chi_{M_c}\chi_8'$ hat f"ur $c\geq 8$ negative
Koeffizienten. Insgesamt ist somit $A_{k+1}=-\alpha_{k+1}$ stets positiv. 

Mit der Rademacherschen Kreismethode~\cite{Ra} kann man nun das Wachstum der
Koeffizienten $u_n$ und $v_n$ der Reihen $q^{c/24}\chi_{M_c}$
und $\chi_8$ absch"atzen (vgl.~\cite{Ap}): Es gilt 
$u_n=e^{\pi\sqrt{(2/3)n}+{\rm O}(\log n)}$ und  
$v_n=e^{\pi\sqrt{(16/3)n}+{\rm O}(\log n)}$. Somit verh"alt sich der
$n$-te Koeffizient von $q^r\left(\chi_{M_c} \chi_8'-\frac{c}{8}\chi_{M_c}'\chi_8
\right)\chi_8^{3r-c/8-1}$ wie $-e^{\pi\sqrt{((2/3)+(16/3)l)n}+{\rm O}(\log n)}$
mit $l=3r-\frac{c}{8}$. Setzen wir $x=3\left(\left(\frac{c}{24}\right)-
\left[\frac{c}{24}\right]\right)$, so strebt daher in Abh"angigkeit der 
Restklasse $c \pmod{24}$ der Quotient 
$\alpha_{k+2}/\alpha_{k+1}=e^{\pi(
\sqrt{((2/3)+(16/3)(6-x))}-\sqrt{((2/3)+(16/3)(3-x))}\,)\sqrt{k}
+{\rm O}(\log k)}$
gegen unendlich.
Insbesondere bleibt $A_{k+2}-A_{k+1}=-\alpha_{k+2}+
(744(1+[c/24]-c/24)+1)\alpha_{k+1}$ stets positiv.\newline
\phantom{fffffffffffff}\hfill \qed

\medskip
{\bf Beispiele von extremalen \sdn \VOAs:} 

F"ur $c=8$ ist $V_{E_8}$ extremal,
f"ur $c=16$ sind $V_{E_8}^{\otimes 2}$ und $V_{D_{16}^+}$ extremal und f"ur
$c=24$ ist die Monster-\VOA $\VM$ extremal. Es wird vermutet, da"s
dies f"ur diese Werte von $c$ die einzigen extremalen \VOAs sind, zumindest
dann, wenn man zus"atzlich unit"ar voraussetzt. Vergleiche hierzu 
die Diskussion in Kapitel 2 bei Vermutung~\ref{vermutungsdvoa}.
F"ur $c=24$ ist dies die Eindeutigkeitsvermutung von Frenkel, Lepowsky und
Meurmann~(vgl.~\cite{FLM}, Einleitung S.~xxxiii).

Zu unimodularen geraden Gittern kann man au"ser der Gitter-\VOA $V_L$ 
stets auch eine $\Z_2$-\-Orbifold-\-\VOA $V_L^{\rm twist}$ konstruieren
(vgl.~\cite{DGH-Monster,DGM,DGM-twist}). Die Konstruktion ist 
mathe\-matisch exakt formulierbar (s.~\cite{Le-Remark}, Abschnitt 4).
Der Charakter in Abh"angig\-keit vom Rang $c$ ist
\begin{eqnarray}\label{extremalorbi}
\chi_{V_L^{\rm twist}}& =& q^{-\frac{c}{24}}\cdot\frac{1}{2}\left(
\frac{\Theta_L}{\prod_{n=1}^{\infty}\left(1-q^n\right)^c}+
\frac{1}{\prod_{n=1}^{\infty}\left(1+q^n\right)^c}\right)+ \nonumber \\
& &\quad
q^{-\frac{c}{24}+\frac{c}{16}}\cdot 2^{\frac{c}{2}}\cdot\frac{1}{2}\left(
\frac{1}{\prod_{n=1}^{\infty}\left(1-q^{n-\frac{1}{2}}\right)^c}+ 
(-1)^{\frac{c}{8}}
\frac{1}{\prod_{n=1}^{\infty}\left(1+q^{n-\frac{1}{2}}\right)^c}\right).
\end{eqnarray}\pagebreak[2]
F"ur $c=32$ und $c=40$ sind daher die $\Z_2$-Orbifold-\VOAs zu den 
extremalen geraden Gittern in diesen Dimensionen extremal,
d.h.~es gilt $V_1=0$. Die Selbstdualit"at ist zu vermuten.
Da in den Dimensionen 32 und 40 nichtisomorphe extremale Gitter existieren
(s.~\cite{CoSl}, Kap.~7.7), gibt es auch nichtisomorphe extremale \VOAs.

F"ur $c \geq 48$ sind keine extremalen \VOAs bekannt.

\medskip 
Tabelle 5.1 enth"alt den Beginn der extremalen Charaktere f"ur
die R"ange $c=8$, $16$, $\ldots$, $48$ und $72$.
\begin{table}\label{extvoatab}
\caption{Charaktere der extremalen \VOAs f"ur die R"ange $c=8$, $16$, $\ldots$,
$48$ und $72$}{\small
$$\begin{array}{|ll|}
\hline
c=8\ \,:& q^{-1/3}(\,1 + 248\, q + 4124\, q^2 + 34752\, q^3 + 
213126\,q^4+1057504\, q^5 +   4530744\,q^6 +\cdots \,) \\  
c=16: &  q^{-2/3}(\,1 + 496\,q + 69752\,q^2 + 2115008\,q^3 + 34670620\,q^4 + 
394460000\,q^5 +  \cdots   \,)  \\ 
c=24: &  q^{-1}(\,1 + 196884\,q^2 + 21493760\,q^3 + 864299970\,q^4 + 
20245856256q^5 + \cdots    \,)  \\ 
c=32: &  q^{-4/3}(\,1 + 139504\,q^2 + 69332992\,q^3 + 6998296696\,q^4 + 
330022830080\,q^5 +   \cdots   \,) \\ 
c=40: &  q^{-5/3}(\, 1 + 20620\,q^2 + 86666240\,q^3 + 24243884350\,q^4 + 
2347780456448\,q^5 +  \cdots  \,)  \\ 
c=48: &  q^{-2}(\,1 + q^2 + 42987520\,q^3 + 40491909396\,q^4 + 
8504046600192\,q^5 +   \cdots   \,)  \\ 
c=72: & q^{-3}(\, 1 + q^2 + q^3 + 2593096794 \,q^4 + 12756091394048 \,q^5 + 
   9529321553850114 \,q^6 +\cdots  \,)   \\  \hline
\end{array}$$}
\end{table}

\section{Extremale selbstduale \SVOAs{ }}

Die extremalen ungeraden selbstdualen Codes 
(s.~\cite{MaSl,Wa})
und Gitter (s.~\cite{CoOdSl}) lassen sich alle explizit angeben. 
In diesem Abschnitt
beweisen wir ein analoges Resultat f"ur extremale \SVOAs, wobei die
vermutete Eindeutigkeit der angegebenen \SVOAs allerdings offen
bleibt.

\medskip
Sei $V$ in diesem Abschnitt eine selbstduale \sehrnette unit"are 
und rationale \SVOA vom Rang $c$. Nach Satz~\ref{sdsvoachar} ist $c$ eine
halbganze Zahl und ihr Charakter ist ein Laurentpolynom in $\chi_{1/2}=
\sqrt[24]{j_{\theta}}=q^{-1/48}(1+q^{1/2}+q^{3/2}+q^2+\cdots)$:
\begin{equation}\label{einss}
\chi_V=\sum_{r=0}^k a_r \chi_{1/2}^{2c-24r},\quad k=\left[\frac{c}{8}\right],
\end{equation}\nopagebreak[2]
mit eindeutigen ganzen Zahlen $a_0$, $\ldots$, $a_k$.
\pagebreak[2]

\begin{definition}\label{extsvoadef} 
Sind die $a_r$ so gew"ahlt, da"s
\begin{equation}\label{extsvoachar}
\chi_V=\chi_{M_c}\cdot(1+A_{k+1}\cdot q^{\frac{k+1}{2}}+
      A_{k+2}\cdot q^{\frac{k+2}{2}}+\cdots),
\end{equation}
so hei"st~(\ref{extsvoachar}) der extremale Charakter in Rang $c$ und eine
\SVOA mit diesem Charakter eine extremale \SVOA.
\end{definition}
Eine extremale selbstduale \SVOA kann sogar auch schon eine \VOA sein.

Wegen Gleichung~(\ref{partgeneral}) gilt f"ur das Mini\-mal\-ge\-wicht ei\-ner
ext\-re\-ma\-len \SVOA \phantom{xxxx} \mbox{$\mu(V)\geq 
\frac{1}{2}\left[\frac{c}{8}\right]+\frac{1}{2}$}.

Hauptresultat dieses Kapitels ist der
\begin{satz}\label{extsvoasatz}
Extremale selbstduale \sehrnette unit"are rationale \SVOAs existieren genau
f"ur die R"ange $c \in E := \{0,\frac{1}{2},1,\ldots,\frac{15}{2},8,12,14,15,
\frac{31}{2},\frac{47}{2},24\}$. F"ur jeden dieser Werte ist genau eine
\SVOA bekannt:
$$\begin{array}{l|cccccccc}
c  &  0\hbox{--}\frac{15}{2} & 8 & 12 & 14 & 15 &\frac{31}{2} & \frac{47}{2} &
 24 \\ \hline
{\rm \SVOA{ }} & \VF^{\otimes 2c} & V_{E_8} & V_{D_{12}^+} & V_{(E_7+E_7)^+} & 
V_{A_{15}^+} & V_{E_{8,2}^+} & \VB & \VM \\
\end{array}$$
Der Satz wird allerdings nur unter den folgenden Einschr"ankungen 
bewiesen:\newline
--- F"ur $c\in\badset$
mu"s die Vollst"andigkeit der Liste der
selbstdualen \VOAs in Vermutung~\ref{SVOAS8-16} vorausgesetzt werden, \newline
--- f"ur $c=\frac{31}{2}$ ist die \SVOA-Struktur nicht vollst"andig
konstruiert,\newline
--- und f"ur $c=\frac{31}{2}$ oder $\frac{47}{2}$ ist die Selbstdualit"at
nicht gezeigt.
\end{satz}

{\it Anmerkungen:} \newline
--- Es wird vermutet, da"s die angegebenen \SVOAs f"ur den jeweiligen
Rang eindeutig sind. 
F"ur $0\leq c <8 $ ist dies Satz~\ref{SVOAS1-8}.
F"ur $8\leq c < 16$ folgt die angegebene Struktur der von den
Gewicht $1$ Vektoren erzeugten Unter-\VOA aus der Annahme der Eigenschaft
\Llie{} (s.~\ref{L1-lie}); siehe dazu auch Vermutung~\ref{SVOAS8-16} und
Satz~\ref{satzSVOAS8-16}.
F"ur $\VB$ wurde die Eindeutigkeit schon in Kapitel~\ref{babykapitel}
vermutet und f"ur $\VM$ ist es wieder die Eindeutigkeitsvermutung aus~\cite{FLM}.\newline
--- Die \SVOAs $V_{E_8}$ und $\VM$ sind sogar extremale \VOAs.

Tabelle 5.2 fa"st alle extremalen ungeraden Codes und
Gitter bzw.~\SVOAs in einem "Ubersichtsdiagramm zusammen. Die Pfeile von
Codes zu Gittern sym\-bo\-lisie\-ren darin die auf S.~\pageref{KA}
beschriebene Konstruktion $C\rightarrow L_C$;
die Pfei\-le von Git\-tern zu \SVOAs sym\-bolisieren die Kon\-struk\-tion der 
Gitter-\SVOA zu diesem Gitter (s.~Satz \ref{satzgittersvoa}).
\begin{table}\label{codesgittersvoas}
\caption{Extremale ungerade Codes, Gitter und \SVOAs{ }} {
$$ \begin{array}{|l|*{15}{c}|}\hline
\mbox{Rang}& \frac{1}{2} & 1 & \frac{3}{2} & 2 &\frac{5}{2} & 3 &\frac{7}{2} 
& 4 &  \frac{9}{2} & 5 & \frac{11}{2} & 6 &\frac{13}{2} & 7 &\frac{15}{2}  \\ \hline
\mbox{Codes} &   &   &   &c_2&   &   &   & c_2^2 &   &   &   &c_2^3&  &  &  \\
             &   &   &   &\da&   &   &   & \da   &   &   &   &\da  &  &  &  \\
\mbox{Gitter}&   &\Z &   &\Z^2&  &\Z^3&  &\Z^4   &   &\Z^5&  &\Z^6 & &\Z^7& \\
             &   &\da&   &\da &  &\da &  &\da    &   &\da&   &\da  &  &\da& \\
\mbox{\SVOAs{ }}& V_F & V_F^2 & V_F^3 & V_F^4  & V_F^5  & V_F^6 &
             V_F^7& V_F^8 & V_F^9 &V_F^{10}&V_F^{11}&V_F^{12}&
              V_F^{13}&V_F^{14}&V_F^{15}\\ \hline
\end{array}$$

$$ \begin{array}{|l|*{9}{c}|}\hline
\mbox{Rang} & 8& 12 & 14 & 15 & \frac{31}{2} & 22 & 23 & \frac{47}{2} & 24 \\ \hline 
\mbox{Codes} & e_8 & d_{12}^+ &(e_7+e_7)^+ & &   & g_{22} & &  & g_{24}\\
             &\da  &\da       & \da        & &   &        & &  &   \\
\mbox{Gitter} &E_8  &
          D_{12}^+&(E_7+E_7)^+& A_{15}^+& & & O_{23} &  & \Lambda_{24} \\
            &\da & \da & \da & \da &   &   &  &    &   \\
\mbox{\SVOAs{ }} &V_{E_8}&
        V_{D_{12}^+}&V_{(E_7+E_7)^+}&V_{A_{15}^+}&V_{E_{8,2}^+}& & &\VB&\VM \\
\hline
\end{array}$$}
\end{table}

Aus Satz~\ref{extsvoasatz} folgt unmittelbar das
\begin{korollar}
Das Minimalgewicht einer selbstdualen \SVOA $V$ erf"ullt
$\mu(V)\leq\frac{1}{2}\left[\frac{c}{8}\right]+\frac{1}{2}$.
\end{korollar}
Die $\mu(V)\leq\frac{1}{2}[\frac{c}{8}]+\frac{1}{2}$ entsprechenden
Absch"atzungen f"ur das Minimalgewicht von Codes bzw.~Gittern sind
asymptotisch schlecht (siehe~\cite{CoSl-JN,CoSl-JN-korrekt,cosl-f1}). Vermutlich
gilt ein "ahnliches Resultat wie dort auch f"ur \SVOAs. Wir definieren 
extremal trotzdem wie in Definition~\ref{extsvoadef}, um mit der Notation bei
Codes und Gittern konform zu sein.

{\bf Beweis von Satz~\ref{extsvoasatz}:}
Die Beweisidee ist v"ollig analog zum Beweis der entsprechenden S"atze
f"ur ungerade Codes bzw.~Gitter.

{\em Existenz der extremalen \SVOAs f"ur $c \in E$:}

F"ur $0\leq c <8$ ist der extremale Charakter $\chi_{1/2}^{2c}$.
Nach (\ref{vsok}) ist $\VF^{2c}$ eine selbstduale \SVOA mit diesem
Charakter.

Im Bereich $8\leq c <16$ beginnt der extremale Charakter mit
$\chi_V=q^{-c/24}(1+A_2 q+ \cdots )$. F"ur die in Satz~\ref{SVOAS8-16}
angegebenen \SVOAs gilt $V_{1/2}=0$, und die Selbstdualit"at f"ur
$V_{E_8}$, $V_{D_{12}^+}$, $V_{(E_7+E_7)^+}$ und $V_{A_{15}^+}$ folgt aus
der Konstruktion als Gitter-\SVOA. F"ur $V_{E_{8,2}^+}$ sind Existenz
und Selbstdualit"at nicht vollst"andig gezeigt.

F"ur $c=47/2$ ist die in Kapitel 4 konstruierte Babymonster-\SVOA $\VB$
extremal, da nach Satz~\ref{babychar} 
$\chi_{\VB}=\chi_{1/2}^{47}-47\chi_{1/2}^{23}
=q^{-47/48}\,(1+4371\,q^{3/2}+\cdots)$; die
Selbstdualit"at wurde vermutet (Vermutung~\ref{vermbabyext}).

Schlie"slich ist f"ur $c=24$ die Monster-\VOA $\VM$ extremal. Die
Selbstdualit"at wurde von C.~Dong in~\cite{Do} gezeigt.

\begin{table}\label{tabcharsvoaext}
\caption{Charaktere der extremalen \SVOAs{ }}
{\scriptsize $$ \begin{array}{|c|c|l|}\hline
V & c & \hbox{Charaktere} \\ \hline\hline
{\scriptstyle \bf 1} & 0 &
1 \\ 
& &
{\rm VOA} \\ \hline
\VF &{1\over 2} & q^{- {1 \over 48}}(
1 + {q^{{1\over 2}}} + {q^{{3\over 2}}} + {q^2} + {q^{{5\over 2}}} + {q^3} + 
  {q^{{7\over 2}}} + 2\,{q^4} + 2\,{q^{{9\over 2}}} + 2\,{q^5} + 
  2\,{q^{{{11}\over 2}}} +\cdots \,)\\
& & q^{- {1 \over 48}}(
{q^{{1\over {16}}}} + {q^{{{17}\over {16}}}} + {q^{{{33}\over {16}}}} + 
  2\,{q^{{{49}\over {16}}}} + 2\,{q^{{{65}\over {16}}}} + 
  3\,{q^{{{81}\over {16}}}} + 4\,{q^{{{97}\over {16}}}}+\cdots \,)
\\ \hline
\VF^{\otimes 2} &1 & q^{- {2 \over 48}}(
1 + 2\,{q^{{1\over 2}}}  + q + 2\,{q^{{3\over 2}}} + 4\,{q^2} + 
  4\,{q^{{5\over 2}}} + 5\,{q^3} + 6\,{q^{{7\over 2}}} + 9\,{q^4} + 
  12\,{q^{{9\over 2}}} + 13\,{q^5} + 16\,{q^{{{11}\over 2}}} +\cdots \,)\\
& & q^{- {2 \over 48}}(
{q^{{1\over 8}}} + 2\,{q^{{9\over 8}}} + 3\,{q^{{{17}\over 8}}} + 
  6\,{q^{{{25}\over 8}}} + 9\,{q^{{{33}\over 8}}} + 14\,{q^{{{41}\over 8}}} + 
  22\,{q^{{{49}\over 8}}}+\cdots \,)
\\ \hline
\VF^{\otimes 3} &{3\over 2} & q^{- {3 \over 48}}(
1 + 3\,{q^{{1\over 2}}}  + 3\,q + 4\,{q^{{3\over 2}}} + 9\,{q^2} + 
  12\,{q^{{5\over 2}}} + 15\,{q^3} + 21\,{q^{{7\over 2}}} + 30\,{q^4} + 
  43\,{q^{{9\over 2}}} + 54\,{q^5} + 69\,{q^{{{11}\over 2}}} + \cdots \,)\\
& & q^{- {3 \over 48}}(
2\,{q^{{3\over {16}}}} + 6\,{q^{{{19}\over {16}}}} + 
  12\,{q^{{{35}\over {16}}}} + 26\,{q^{{{51}\over {16}}}} + 
  48\,{q^{{{67}\over {16}}}} + 84\,{q^{{{83}\over {16}}}} + 
  146\,{q^{{{99}\over {16}}}} +\cdots \,)
\\  \hline
\VF^{\otimes 4} &2 & q^{- {4 \over 48}}(
1 + 4\,{q^{{1\over 2}}}  + 6\,q + 8\,{q^{{3\over 2}}} + 17\,{q^2} + 
  28\,{q^{{5\over 2}}} + 38\,{q^3} + 56\,{q^{{7\over 2}}} + 84\,{q^4} + 
  124\,{q^{{9\over 2}}} + 172\,{q^5} +\cdots \,) \\
& & q^{- {4 \over 48}}(
2\,{q^{{1\over 4}}} + 8\,{q^{{5\over 4}}} + 20\,{q^{{9\over 4}}} + 
  48\,{q^{{{13}\over 4}}} + 102\,{q^{{{17}\over 4}}} + 
  200\,{q^{{{21}\over 4}}} + 380\,{q^{{{25}\over 4}}}+\cdots \,)
\\  \hline
\VF^{\otimes 5} &{5\over 2} & q^{- {5 \over 48}}(
1 + 5\,{q^{{1\over 2}}}  + 10\,q + 15\,{q^{{3\over 2}}} + 30\,{q^2} + 
  56\,{q^{{5\over 2}}} + 85\,{q^3} + 130\,{q^{{7\over 2}}} + 205\,{q^4} + 
  315\,{q^{{9\over 2}}} + 465\,{q^5} +\cdots \,)\\
& & q^{- {5 \over 48}}(
4\,{q^{{5\over {16}}}} + 20\,{q^{{{21}\over {16}}}} + 
  60\,{q^{{{37}\over {16}}}} + 160\,{q^{{{53}\over {16}}}} + 
  380\,{q^{{{69}\over {16}}}} + 824\,{q^{{{85}\over {16}}}}+\cdots \,)
\\ \hline
\VF^{\otimes 6} &3 & q^{- {6 \over 48}}(
1 + 6\,{q^{{1\over 2}}} + 15\,q + 26\,{q^{{3\over 2}}} + 51\,{q^2} + 
  102\,{q^{{5\over 2}}} + 172\,{q^3} + 276\,{q^{{7\over 2}}} + 453\,{q^4} + 
  728\,{q^{{9\over 2}}} + 1128\,{q^5}   +\cdots \,)\\
& & q^{- {6 \over 48}}(
4\,{q^{{3\over 8}}} + 24\,{q^{{{11}\over 8}}} + 84\,{q^{{{19}\over 8}}} + 
  248\,{q^{{{27}\over 8}}} + 648\,{q^{{{35}\over 8}}} + 
  1536\,{q^{{{43}\over 8}}}+\cdots \,)
\\ \hline
\VF^{\otimes 7} &{7\over 2}& q^{- {7 \over 48}}(
1 + 7\,{q^{{1\over 2}}} + 21\,q + 42\,{q^{{3\over 2}}} + 84\,{q^2} + 
  175\,{q^{{5\over 2}}} + 322\,{q^3} + 547\,{q^{{7\over 2}}} + 931\,{q^4} + 
  1561\,{q^{{9\over 2}}} + 2527\,{q^5} +\cdots \,)\\
& & q^{- {7 \over 48}}(
8\,{q^{{7\over {16}}}} + 56\,{q^{{{23}\over {16}}}} + 
  224\,{q^{{{39}\over {16}}}} + 728\,{q^{{{55}\over {16}}}} + 
  2072\,{q^{{{71}\over {16}}}} + 5320\,{q^{{{87}\over {16}}}}+\cdots \,)
\\ \hline
\VF^{\otimes 8} &4& q^{- {8\over48}}(
1 + 8\,{q^{{1\over 2}}} + 28\,q + 64\,{q^{{3\over 2}}} + 134\,{q^2} + 
  288\,{q^{{5\over 2}}} + 568\,{q^3} + 1024\,{q^{{7\over 2}}} + 1809\,{q^4} + 
  3152\,{q^{{9\over 2}}} + \cdots \,)\\
& & q^{- {8 \over 48}}(
8\,{q^{{1\over 2}}} + 64\,{q^{{3\over 2}}} + 288\,{q^{{5\over 2}}} + 
  1024\,{q^{{7\over 2}}} + 3152\,{q^{{9\over 2}}} 
+\cdots \,)\\ \hline
\VF^{\otimes 9} &{9\over 2}& q^{- {9 \over 48}}(
1 + 9\,{q^{{1\over 2}}} + 36\,q + 93\,{q^{{3\over 2}}} + 207\,{q^2} + 
  459\,{q^{{5\over 2}}} + 957\,{q^3} + 1827\,{q^{{7\over 2}}} + 3357\,{q^4} + 
  6061\,{q^{{9\over 2}}} +\cdots \,)\\
& & q^{- {9 \over 48}}(
16\,{q^{{9\over {16}}}} + 144\,{q^{{{25}\over {16}}}} + 
  720\,{q^{{{41}\over {16}}}} + 2784\,{q^{{{57}\over {16}}}} + 
  9216\,{q^{{{73}\over {16}}}} + 27216\,{q^{{{89}\over {16}}}}+\cdots \,)
\\ \hline
\VF^{\otimes 10} &5& q^{- {10 \over 48}}(
1 + 10\,{q^{{1\over 2}}} + 45\,q + 130\,{q^{{3\over 2}}} + 310\,{q^2} + 
  712\,{q^{{5\over 2}}} + 1555\,{q^3} + 3130\,{q^{{7\over 2}}} + 
  5990\,{q^4} + 11190\,{q^{{9\over 2}}} + \cdots \,) \\
& & q^{- {10 \over 48}}(
16\,{q^{{5\over 8}}} + 160\,{q^{{{13}\over 8}}} + 880\,{q^{{{21}\over 8}}} + 
  3680\,{q^{{{29}\over 8}}} + 13040\,{q^{{{37}\over 8}}} + 
  40992\,{q^{{{45}\over 8}}}+\cdots \,)
\\ \hline
\VF^{\otimes 11} &{{11}\over 2}& q^{- {11 \over 48}}(
1 + 11\,{q^{{1\over 2}}} + 55\,q + 176\,{q^{{3\over 2}}} + 451\,{q^2} + 
  1078\,{q^{{5\over 2}}} + 2453\,{q^3} + 5181\,{q^{{7\over 2}}} + 
  10329\,{q^4} + 19954\,{q^{{9\over 2}}}  +\cdots \,)\\
& & q^{- {11 \over 48}}(
32\,{q^{{{11}\over {16}}}} + 352\,{q^{{{27}\over {16}}}} + 
  2112\,{q^{{{43}\over {16}}}} + 9504\,{q^{{{59}\over {16}}}} + 
  35904\,{q^{{{75}\over {16}}}} + 119680\,{q^{{{91}\over {16}}}}+\cdots \,) 
\\ \hline
\VF^{\otimes 12} &6& q^{- {12 \over 48}}(
1 + 12\,{q^{{1\over 2}}} + 66\,q + 232\,{q^{{3\over 2}}} + 639\,{q^2} + 
  1596\,{q^{{5\over 2}}} + 3774\,{q^3} + 8328\,{q^{{7\over 2}}} + 
  17283\,{q^4} + 34520\,{q^{{9\over 2}}} +\cdots \,) \\
& & q^{- {12 \over 48}}(
32\,{q^{{3\over 4}}} + 384\,{q^{{7\over 4}}} + 2496\,{q^{{{11}\over 4}}} + 
  12032\,{q^{{{15}\over 4}}} + 48288\,{q^{{{19}\over 4}}} + 
  170112\,{q^{{{23}\over 4}}} +\cdots \,) 
\\ \hline
\VF^{\otimes 13} &{{13}\over 2}& q^{- {13 \over 48}}(
1 + 13\,{q^{{1\over 2}}} + 78\,q + 299\,{q^{{3\over 2}}} + 884\,{q^2} + 
  2314\,{q^{{5\over 2}}} + 5681\,{q^3} + 13052\,{q^{{7\over 2}}} + 
  28158\,{q^4} + 58136\,{q^{{9\over 2}}}  +\cdots \,) \\
& & q^{- {13 \over 48}}(
64\,{q^{{{13}\over {16}}}} + 832\,{q^{{{29}\over {16}}}} + 
  5824\,{q^{{{45}\over {16}}}} + 29952\,{q^{{{61}\over {16}}}} + 
  127296\,{q^{{{77}\over {16}}}} + 472576\,{q^{{{93}\over {16}}}}+\cdots \,)
\\ \hline
\VF^{\otimes 14} &7& q^{- {14 \over 48}}(
1 + 14\,{q^{{1\over 2}}} + 91\,q + 378\,{q^{{3\over 2}}} + 1197\,{q^2} + 
  3290\,{q^{{5\over 2}}} + 8386\,{q^3} + 20008\,{q^{{7\over 2}}} + 
  44800\,{q^4} +\cdots \,)\\
& & q^{- {14 \over 48}}(
64\,{q^{{7\over 8}}} + 896\,{q^{{{15}\over 8}}} + 6720\,{q^{{{23}\over 8}}} + 
  36736\,{q^{{{31}\over 8}}} + 164864\,{q^{{{39}\over 8}}} + 
  643328\,{q^{{{47}\over 8}}}+\cdots \,)
\\ \hline
\VF^{\otimes 15} &{{15}\over 2}& q^{- {15 \over 48}}(
1 + 15\,{q^{{1\over 2}}} + 105\,q + 470\,{q^{{3\over 2}}} + 1590\,{q^2} + 
  4593\,{q^{{5\over 2}}} + 12160\,{q^3} + 30075\,{q^{{7\over 2}}} + 
  69780\,{q^4}  +\cdots \,)\\
& & q^{- {15 \over 48}}(
128\,{q^{{{15}\over {16}}}} + 1920\,{q^{{{31}\over {16}}}} + 
  15360\,{q^{{{47}\over {16}}}} + 88960\,{q^{{{63}\over {16}}}} + 
  420480\,{q^{{{79}\over {16}}}} + 1720704\,{q^{{{95}\over {16}}}}+\cdots \,)
\\ \hline
V_{E_8}& 8& q^{- {1\over3}}(
1 + 248\,q + 4124\,{q^2} + 34752\,{q^3} + 213126\,{q^4} + 1057504\,{q^5} + 
  4530744\,{q^6}+ \cdots \,)\\
& &
{\rm VOA}
\\ \hline
V_{D_{12}^+}&12 & q^{- {24 \over 48}}(
1 + 276\,q + 2048\,{q^{{3\over 2}}} + 11202\,{q^2} + 
  49152\,{q^{{5\over 2}}} + 184024\,{q^3} + 614400\,{q^{{7\over 2}}} + 
  1881471\,{q^4}  \cdots \,) \\
& & q^{- {24 \over 48}}(
12\,{q^{{1\over 2}}} + 2048\,{q^{{3\over 2}}} + 49152\,{q^{{5\over 2}}} + 
  614400\,{q^{{7\over 2}}} + 5373952\,{q^{{9\over 2}}} + 
  37122048\,{q^{{{11}\over 2}}}+\cdots \,)
\\ \hline
\!\!\!V_{\!(\!E_7\!+\!E_7\!\,)^{\!+}}\!\!\!\!&14 &  q^{- {28 \over 48}}(
1 + 266\,q + 3136\,{q^{{3\over 2}}} + 21035\,{q^2} + 
  108416\,{q^{{5\over 2}}} + 468846\,{q^3} + 1777472\,{q^{{7\over 2}}} + 
  6094557\,{q^4} + \cdots \,) \\
& & q^{- {28 \over 48}}(
56\,{q^{{3\over 4}}} + 8416\,{q^{{7\over 4}}} + 229936\,{q^{{{11}\over 4}}} + 
  3327296\,{q^{{{15}\over 4}}} + 33491752\,{q^{{{19}\over 4}}} + 
  264189408\,{q^{{{23}\over 4}}}+ \cdots \,) \\ \hline
V_{A_{15}^+}&15 & q^{- {30 \over 48}}(
1 + 255\,q + 3640\,{q^{{3\over 2}}} + 27525\,{q^2} + 
  154056\,{q^{{5\over 2}}} + 713850\,{q^3} + 2878920\,{q^{{7\over 2}}} + 
  10432650\,{q^4} + \cdots \,) \\
& & q^{- {30 \over 48}}(
120\,{q^{{7\over 8}}} + 17104\,{q^{{{15}\over 8}}} + 
  494040\,{q^{{{23}\over 8}}} + 7626000\,{q^{{{31}\over 8}}} + 
  81775600\,{q^{{{39}\over 8}}} + 685224960\,{q^{{{47}\over 8}}}+ \cdots \,)
\\ \hline
V_{E_{8,2}^+}&{{31}\over 2}& q^{- {31 \over 48}}(
1 + 248\,q + 3875\,{q^{{3\over 2}}} + 31124\,{q^2} + 
  181753\,{q^{{5\over 2}}} + 871627\,{q^3} + 3623869\,{q^{{7\over 2}}} + 
  13496501\,{q^4} + \cdots \,) \\
& & q^{- {31 \over 48}}(
248\,{q^{{{15}\over {16}}}} + 34504\,{q^{{{31}\over {16}}}} + 
  1022752\,{q^{{{47}\over {16}}}} + 16275496\,{q^{{{63}\over {16}}}} + 
  179862248\,{q^{{{79}\over {16}}}} + \cdots \,)
\\ \hline
\VB&{{47}\over 2}& q^{- {47 \over 48}}(
1 + 4371\,{q^{{3\over 2}}} + 96256\,{q^2} + 1143745\,{q^{{5\over 2}}} + 
  9646891\,{q^3} + 64680601\,{q^{{7\over 2}}} + 366845011\,{q^4} + 
   \cdots \,)
  \\
& & q^{- {47 \over 48}}(
96256\,{q^{{{31}\over {16}}}} + 10602496\,{q^{{{47}\over {16}}}} + 
  420831232\,{q^{{{63}\over {16}}}} + 9685952512\,{q^{{{79}\over {16}}}} + 
  156435924992\,{q^{{{95}\over {16}}}}+ \cdots \,)
\\ \hline
\VM & 24 & q^{- 1}(
1 + 196884\,{q^2} + 21493760\,{q^3} + 864299970\,{q^4} + 20245856256\,{q^5} + 
  333202640600\,{q^6} + \cdots \,)\\
& &{\rm VOA} \\
 \hline
\end{array}$$}
\end{table}
In Tabelle~\ref{tabcharsvoaext} sind alle Charaktere der extremalen
\SVOAs aufgelistet. Die erste Spalte gibt die \SVOA, die zweite den Rang an. 
In der dritten Spalte zeigt
die erste Zeile eines Eintrages den Charakter $\chi_V=\chi_{V_{(0)}}+
\chi_{V_{(1)}}$, die zweite den Charakter $\chi_{V_{(2)}}=\chi_{V_{(3)}}$
der irreduziblen $V_{(0)}$-Moduln $V_{(2)}$ und $V_{(3)}$ (falls existent).

\medskip
{\em Nichtexistenz der extremalen \SVOAs f"ur $c \not\in E:$}

Falls $V$ sogar eine selbstduale \VOA ist, gilt nach Korollar~\ref{korminVOA}
$\mu(V)\leq [c/24]+1$ und $8\vert c$. Da f"ur $V$ als extremale \SVOA
$\mu(V)\geq \frac{1}{2}[\frac{c}{8}]+\frac{1}{2}$ ist, folgt $c=8$ oder $24$,
d.h.~$c\in E$. Wir k"onnen daher voraussetzen, da"s $V$ eine 
\SVOA mit $V_{(1)}\not=0$ ist.

Sei $\V0$ die gerade Unter-\VOA von $V$. Da wir $V$ als \sehrnett 
vorausgesetzt
haben, gilt $V=\V0\oplus\V1$ mit irreduziblen $\V0$-Moduln $\V0$ und
$\V1$. Zus"atzlich besitzt $\V0$ f"ur $c\in\Z+\frac{1}{2}$ einen dritten
irreduziblen Modul $\V2$ bzw.~f"ur $c\in\Z$ zwei irreduzible Moduln $\V2$
und $\V3$. F"ur den Charakter gilt 
$\chi_V=\chi_{\V0}+
\chi_{\V1}$. Bezeichnen wir mit $\widetilde{\chi}_V$ die Entwicklung von
$\chi_V$ in der Spitze $1$, so erhalten wir nach Definition~\ref{sehrnett}
und Satz~\ref{sdsvoafusion} den Charakter der Moduln $\V2$ bzw.~$\V3$:
\begin{equation}\label{plus}
e^{2c\cdot\frac{2 \pi i}{48}}\,\widetilde{\chi}_V=\cases{\chi_{\V2}+\chi_{\V3},
 &   falls $c\in\Z$, \cr
\sqrt{2}\chi_{\V2}, & falls $c\in\Z+\frac{1}{2}$.}
\end{equation}
Die Koeffizienten der Charaktere $\chi_{\V0}$, $\chi_{\V1}$, $\chi_{\V2}+\chi_{\V3}$ m"ussen positive ganze Zahlen sein. 
Wir haben somit die
folgenden Nichtexistenzargumente:\newline
{\em Argument $N$:} Mindestens ein Koeffizient von 
$\chi_{V_{(2)}}+\chi_{V_{(3)}}$
ist negativ.\newline
{\em Argument $G$:} Mindestens ein Koeffizient von 
$\chi_{V_{(2)}}+\chi_{V_{(3)}}$ ist nicht ganz.\newline
In den F"allen $c\in \badset$
reichen diese Argumente nicht aus, und wir ben"otigen\newline
{\em Argument $L$:} 
Die Liste der in (\ref{listeSVOAS8-16}) angegebenen selbstdualen \SVOAs im 
Bereich \hbox{$8\leq c < 16$} ist vollst"andig (Vermutung~\ref{SVOAS8-16},
vgl.~die dortige Diskussion und Satz~\ref{satzSVOAS8-16}).

Wir fassen die Rechnungen f"ur den Bereich $8<c<48$ in zwei Tabellen
zusammen.\newline
Tabelle~5.4 behandelt die F"alle $8<c<24$, $c\not\in E$:
In Spalte 1 ist der Rang $c$ angegeben, Spalte 2 und 3 geben die ersten
Koeffizienten der $q$-Entwicklung von $q^{c/24}\,(\chi_{\V0}+\chi_{\V1})$
bzw.~$q^{c/24}\,(\chi_{\V2}+\chi_{\V3})$ an und in Spalte 4 sind das oder die 
Nichtexistenzargumente angegeben.\newline
\begin{table}\label{c24}
\caption{Extremale Charaktere f"ur $c<24$, $c\not\in E$}
{\scriptsize $$ \begin{array}{|r|l|l|l|}\hline
{{17}\over 2} &
1 + 255\,q + 221\,{q^{{3\over 2}}} + 4216\,{q^2} + 4114\,{q^{{5\over 2}}} + 
  35666\,{q^3} + \cdots
&
{{17}\over {16}} \,{q^{{1\over {16}}}}+ 
  {{3977}\over {16}}\,{q^{{{17}\over {16}}}} + 
  {{69989}\over {16}}\,{q^{{{33}\over {16}}}} + \cdots & G
 \\
9 &
1 + 261\,q + 456\,{q^{{3\over 2}}} + 4500\,{q^2} + 8424\,{q^{{5\over 2}}} + 
  40641\,{q^3} +\cdots 
&
{{9}\over 4}\,{q^{{1\over 8}}} + {{997}\over 2}\,{q^{{9\over 8}}} + 
  {{36999}\over 4} \,{q^{{{17}\over 8}}} + 
\cdots &G
 \\
{{19}\over 2} &
1 + 266\,q + 703\,{q^{{3\over 2}}} + 4997\,{q^2} + 13091\,{q^{{5\over 2}}} + 
  49989\,{q^3}+\cdots &
{{19}\over 8}\,{q^{{3\over {16}}}} + 
  {{4001}\over 8}\,{q^{{{19}\over {16}}}} + 
  {{39007}\over 4}\,{q^{{{35}\over {16}}}} + 
 \cdots &G
  \\
10 &
1 + 270\,q + 960\,{q^{{3\over 2}}} + 5725\,{q^2} + 18304\,{q^{{5\over 2}}} + 
  64150\,{q^3} + \cdots&
5\,{q^{{1\over 4}}} + 1004\,{q^{{5\over 4}}} + 20510\,{q^{{9\over 4}}} + 
\cdots & L \\
{{21}\over 2} &
1 + 273\,q + 1225\,{q^{{3\over 2}}} + 6699\,{q^2} + 24276\,{q^{{5\over 2}}} + 
  83727\,{q^3} + \cdots&
{{21}\over 4}\,{q^{{5\over {16}}}} + 
  {{4033}\over 4}\,{q^{{{21}\over {16}}}} + 
  {{86079}\over 4}\,{q^{{{37}\over {16}}}} + \cdots& G \\
11 &
1 + 275\,q + 1496\,{q^{{3\over 2}}} + 7931\,{q^2} + 31240\,{q^{{5\over 2}}} + 
  109516\,{q^3} + \cdots
&
11\,{q^{{3\over 8}}} + 2026\,{q^{{{11}\over 8}}} + 
  45067\,{q^{{{19}\over 8}}} + \cdots & L
\\
{{23}\over 2} &
1 + 276\,q + 1771\,{q^{{3\over 2}}} + 9430\,{q^2} + 39445\,{q^{{5\over 2}}} + 
  142531\,{q^3} + \cdots
&
{{23}\over 2}\,{q^{{7\over {16}}}} + 
  {{4073}\over 2}\,{q^{{{23}\over {16}}}} +
 47104\,{q^{{{39}\over {16}}}} + 
\cdots& G
\\
{{25}\over 2} &
1 + 275\,q + 2325\,{q^{{3\over 2}}} + 13250\,{q^2} + 
  60630\,{q^{{5\over 2}}} + 235500\,{q^3} +\cdots&
 25\,{q^{{9\over {16}}}} + 4121\,{q^{{{25}\over {16}}}} + 
  102425\,{q^{{{41}\over {16}}}} +\cdots& L
\\
13 &
1 + 273\,q + 2600\,{q^{{3\over 2}}} + 15574\,{q^2} + 
  74152\,{q^{{5\over 2}}} + 298727\,{q^3} +\cdots &
 52\,{q^{{5\over 8}}} + 8296\,{q^{{{13}\over 8}}} + 
  213148\,{q^{{{21}\over 8}}} +\cdots &L
\\
{{27}\over 2} &
1 + 270\,q + 2871\,{q^{{3\over 2}}} + 18171\,{q^2} + 
  89991\,{q^{{5\over 2}}} + 375741\,{q^3} + \cdots&
54\,{q^{{{11}\over {16}}}} + 8354\,{q^{{{27}\over {16}}}} + 
  221508\,{q^{{{43}\over {16}}}} +\cdots& L
\\
{{29}\over 2} &
1 + 261\,q + 3393\,{q^{{3\over 2}}} + 24157\,{q^2} + 
  129688\,{q^{{5\over 2}}} + 580609\,{q^3} + \cdots &
116\,{q^{{{13}\over {16}}}} + 16964\,{q^{{{29}\over {16}}}} + 
  476876\,{q^{{{45}\over {16}}}}+ \cdots & L
\\
16  &
1 + 7936\,{q^{{3\over 2}}} + 2296\,{q^2} + 412672\,{q^{{5\over 2}}} + 
  65536\,{q^3}+\cdots
&
-{{15}\over {16}} + 527\,q + {{139039\over 2}\,{q^2}} + 2116124\,{q^3}
\cdots& N,G
\\
{{33}\over 2} &
1 + 7766\,{q^{{3\over 2}}} + 11220\,{q^2} + 408507\,{q^{{5\over 2}}} + 
  515251\,{q^3} + \cdots&
{{-231}\over {256}}\,{q^{{1\over {16}}}} + 
  {{138633}\over {256}}\,{q^{{{17}\over {16}}}} + 
  {{17969473}\over {256}}\,{q^{{{33}\over {16}}}} 
 +\cdots & N, G
 \\
17 &
1 + 7582\,{q^{{3\over 2}}} + 19907\,{q^2} + 413678\,{q^{{5\over 2}}} + 
  956573\,{q^3}+\cdots
&
{{-221}\over {128}} \,{q^{{1\over 8}}}+ 
  {{71179}\over {64}}\,{q^{{9\over 8}}} + 
  {{18149745}\over {128}}\,{q^{{{17}\over 8}}}  
+ \cdots & N, G
 \\
{{35}\over 2} &
1 + 7385\,{q^{{3\over 2}}} + 28315\,{q^2} + 427987\,{q^{{5\over 2}}} + 
  1398635\,{q^3} +\cdots &
{{-105}\over {64}} \,{q^{{3\over {16}}}}+ 
  {{73045}\over {64}}\,{q^{{{19}\over {16}}}} + 
  {{4584449}\over {32}}\,{q^{{{35}\over {16}}}} +\cdots& N, G
\\
18 &
1 + 7176\,{q^{{3\over 2}}} + 36405\,{q^2} + 451152\,{q^{{5\over 2}}} + 
  1850520\,{q^3} + \cdots &
{{-99}\over {32}}\,{q^{{1\over 4}}} + {{18729}\over 8}\,{q^{{5\over 4}}} + 
  {{4633405}\over {16}}\,{q^{{9\over 4}}} + 
\cdots& N, G
  \\
{{37}\over 2} &
1 + 6956\,{q^{{3\over 2}}} + 44141\,{q^2} + 482813\,{q^{{5\over 2}}} + 
  2321121\,{q^3} +\cdots &
{{-185}\over {64}}\,{q^{{5\over {16}}}} + 
  {{153587}\over {64}}\,{q^{{{21}\over {16}}}} + 
  {{18737217}\over {64}}\,{q^{{{37}\over {16}}}} + 
\cdots& N, G
 \\
19 &
1 + 6726\,{q^{{3\over 2}}} + 51490\,{q^2} + 522538\,{q^{{5\over 2}}} + 
  2819011\,{q^3} +\cdots &
 {{-171}\over {32}} \,{q^{{3\over 8}}}+ 
  {{78679}\over {16}} \,{q^{{{11}\over 8}}}+ 
  {{18948593}\over {32}}\,{q^{{{19}\over 8}}} + 
\cdots& N, G
\\
{{39}\over 2} &
1 + 6487\,{q^{{3\over 2}}} + 58422\,{q^2} + 569829\,{q^{{5\over 2}}} + 
  3352323\,{q^3} + \cdots &
{{-39}\over 8}\,{q^{{7\over {16}}}} + 
  {{40287}\over 8}\,{q^{{{23}\over {16}}}} + 
  {{1197985}\over 2}\,{q^{{{39}\over {16}}}} + 
\cdots& N, G
\\
20 &
1 + 6240\,{q^{{3\over 2}}} + 64910\,{q^2} + 624128\,{q^{{5\over 2}}} + 
  3928640\,{q^3} +\cdots
&
{{-35}\over 4}\,{q^{1\over 2}} + 10310\,{q^{{3\over 2}}} + 
  1212171\,{q^{{5\over 2}}} +\cdots& N, G
\\
{{41}\over 2} &
1 + 5986\,{q^{{3\over 2}}} + 70930\,{q^2} + 684823\,{q^{{5\over 2}}} + 
  4554895\,{q^3} + \cdots&
 {{-123}\over {16}}\,{q^{{9\over {16}}}} + 
  {{168797}\over {16}}\,{q^{{{25}\over {16}}}} + 
  {{19629545}\over {16}}\,{q^{{{41}\over {16}}}} + 
\cdots& N, G
\\
21 &
1 + 5726\,{q^{{3\over 2}}} + 76461\,{q^2} + 751254\,{q^{{5\over 2}}} + 
  5237281\,{q^3}+ \cdots&
 {{-105}\over 8}\,{q^{{5\over 8}}} + 
{{86331}\over 4} \,{q^{{{13}\over 8}}}+ 
  {{19872217}\over 8}\,{q^{{{21}\over 8}}} + 
\cdots& N, G
 \\
{{43}\over 2} &
1 + 5461\,{q^{{3\over 2}}} + 81485\,{q^2} + 822719\,{q^{{5\over 2}}} + 
  5981171\,{q^3} + \cdots
&
{{-43}\over 4} \,{q^{{{11}\over {16}}}}+ 
  {{88279}\over 4}\,{q^{{{27}\over {16}}}} + 
  {{5030697}\over 2}\,{q^{{{43}\over {16}}}} + 
\cdots& N, G 
 \\
22 &
1 + 5192\,{q^{{3\over 2}}} + 85987\,{q^2} + 898480\,{q^{{5\over 2}}} + 
  6791048\,{q^3} +\cdots&
{{-33\over 2}}\,{q^{{3\over 4}}} +
 45122\,{q^{{7\over 4}}} + 
  5095325\,{q^{{{11}\over 4}}} + \cdots& N, G
\\
{{45}\over 2} &
1 + 4920\,{q^{{3\over 2}}} + 89955\,{q^2} + 977769\,{q^{{5\over 2}}} + 
  7670445\,{q^3} +\cdots &
{{-45}\over 4}\,{q^{{{13}\over {16}}}} + 
  {{184455}\over 4}\,{q^{{{29}\over {16}}}} + 
  {{20922345}\over 2}\,{q^{{{23}\over 8}}}
\cdots & N, G
 \\
23 &
1 + 4646\,{q^{{3\over 2}}} + 93380\,{q^2} + 1059794\,{q^{{5\over 2}}} + 
  8621895\,{q^3} + \cdots &
{{-23}\over 2}\,{q^{{7\over 8}}} + 94231\,{q^{{{15}\over 8}}} + 
  {{20922345}\over 2}\,{q^{{{23}\over 8}}}
\cdots & N, G 
\\   \hline
\end{array} $$ }
\end{table}
\hbox{$\!\!\!$ Tabelle~5.5}
 behandelt die F"alle $24<c<48$: Aufgelistet sind der Rang
und der erste Koeffizient $B^*=B_0$ ($c \in\Z$) bzw.~$B^*=B_0/\sqrt{2}$
($c\in\Z+1/2$) von $\chi_{\V2}+\chi_{\V3}$ bzw.~$\chi_{\V2}$. 
Er ist nie ganz, und daher kann nach
Argument $G$ f"ur diese R"ange keine extremale \SVOA existieren.\newline
\begin{table}\label{c48}
\caption{Erster Koeffizient des extremalen Charakters
bei Entwicklung in der Spitze~$1$ f"ur $24<c<48$, $c\not\in E$}
{\footnotesize $$
\begin{array}{|*{11}{c|}}\hline
\mbox{Rang} & 49/2 & 25 & 51/2& 26& 53/2 & 27 & 55/2& 28 & 57/2& 29  \\ 
\hline
B^* 
&
{{1911}\over {2048}}
&
{{1775}\over {1024}}
&
{{3281}\over {2048}}
&
{{377}\over {128}}
&
{{689}\over {256}}
&
{{1251}\over {256}}
&
{{2255}\over {512}}
&
{{63}\over 8}
&
{{893}\over {128}}
&
{{783}\over {64}}
\\ \hline\hline
\mbox{Rang} & 59/2& 30 & 61/2 & 31 & 63/2 & 32 &65/2  & 33 & 67/2 & 34 \\ 
\hline
B^*
&
{{1357}\over {128}}
&
{{145}\over 8}
&
{{61}\over 4}
&
{{403}\over {16}}
&
{{651}\over {32}}
&
-{{73967}\over {65536}}
&
-{{67989}\over {65536}}
&
-{{31135}\over {16384}}
&
-{{56815}\over {32768}}
&
-{{12907}\over {4096}}
\\ \hline\hline
\mbox{Rang}  & 69/2& 35 & 71/2 & 36 & 73/2 & 37 &75/2  & 38 & 77/2 & 39 \\ 
\hline
B^*
&
-{{5839}\over {2048}}
&
-{{42069}\over {8192}}
&
-{{9425}\over {2048}}
&
-{{33605}\over {4096}}
&
-{{29783}\over {4096}}
&
-{{3279}\over {256}}
&
-{{22949}\over {2048}}
&
-{{9965}\over {512}}
&
-{{8585}\over {512}}
&
-{{14663}\over {512}}

\\ \hline\hline
\mbox{Rang}  & 79/2& 40 & 81/2 & 41 & 83/2 & 42 &85/2  & 43 & 87/2 & 44 \\ 
\hline
B^*
&
-{{6201}\over {256}}
&
{{86155}\over {65536}}
&
{{627945}\over {524288}}
&
{{285213}\over {131072}}
&
{{1033171}\over {524288}}
&
{{29145}\over {8192}}
&
{{839041}\over {262144}}
&
{{376073}\over {65536}}
&
{{335859}\over {65536}}
&
{{74689}\over {8192}}
\\ \hline
\end{array}$$\vskip-7.0mm $$ 
\begin{array}{|*{8}{c|}}\hline
\mbox{Rang}  & 89/2& 45 & 91/2 & 46 & 93/2 & 47 & 95/2  \\ \hline
B^*
&
{{132319}\over {16384}}
&
{{7293}\over {512}}
&
{{409677}\over {32768}}
&
{{44723}\over {2048}}
&
{{310803}\over {16384}}
&
{{134231}\over {4096}}
&
{{28811}\over {1024}}
\\ \hline
\end{array}\qquad\qquad\qquad\qquad\qquad\qquad\quad\quad\quad\ \,
$$}\end{table}
\hbox{$\!$ Die} Eintr"age in den Tabellen~5.4 und~5.5 ergeben sich leicht
durch Gleichsetzen von~(\ref{einss}) und~(\ref{extsvoachar}), sowie
durch~(\ref{plus}) und der Entwicklung von $\chi_{1/2}$ in der
Spitze $1$ in der lokalen Koordinate $q^{1/24}$ (vgl.~(\ref{spitze-1})):
\begin{equation}\label{chi12}
\widetilde{\chi}_{1/2}=e^{-\frac{2\pi i}{48}}\sqrt{2}\cdot q^{\frac{1}{24}}
\prod_{n=1}^{\infty}(1+q^n)=e^{-\frac{2\pi i}{48}}\sqrt{2}\cdot
q^{\frac{1}{24}}(1+q+q^2+\cdots).
\end{equation}

\medskip
Zu erledigen bleiben noch die F"alle $c\geq 48$. Wir werden zeigen:
In der Entwicklung $e^{2c\cdot\frac{2 \pi i}{48}}\widetilde{\chi}_{V}=
q^{c/12-[c/8]}(B_0+B_1\, q +\cdots)$ von $\chi_V$ in der Spitze $1$ 
ist entweder $B_0$ oder $B_1$
negativ --- Widerspruch nach Argument $N$ bei Verwendung von (\ref{plus}).

Dazu setzen wir $p=q^{1/2}$ und entwickeln, "ahnlich wie beim Beweis
von Satz~\ref{extvoasatz}, $\chi_{M_c}\cdot\chi_{1/2}^{-2c}$ in Potenzen
von $\phi:=\chi_{1/2}^{-24}=j_{\theta}^{-1}=p+O(p^2)$:
\begin{equation}\label{dreis}
\chi_{M_c}\cdot\chi_{1/2}^{-2c}=\sum_{r=0}^{\infty}\alpha_r\phi^r,
\end{equation}
wobei
\begin{eqnarray}
\alpha_r & =& \frac{1}{r!}\frac{d^{r-1}}{dp^{r-1}}\left\{\frac
{d(\chi_{M_c}\cdot\chi_{1/2}^{-2c})}{dp}\left(\frac{p}{\phi}\right)^r
\right\}_{p=0}  \nonumber \\ 
 & = & \frac{1}{r!}\frac{d^{r-1}}{dp^{r-1}}\left\{
 p^r\cdot \chi_{1/2}^{24r-2c-3}\left[\chi_{M_c}'\chi_{1/2}^3
-2c\cdot\chi_{M_c}\chi_{1/2}^2 \chi_{1/2}'\right]\right\}_{p=0}.\label{viers}
\end{eqnarray}
Gleichungen (\ref{einss}), (\ref{extsvoachar}) und (\ref{dreis}) liefern
\begin{eqnarray*}
\sum_{r=0}^{k}a_r\phi^r & = & \chi_V\cdot\chi_{1/2}^{-2c} =
\chi_V\cdot\chi_{M_c}^{-1}\cdot\chi_{M_c}\cdot\chi_{1/2}^{-2c} \\
 & = & \left(1+\sum_{n=k+1}^{\infty}A_n p^n\right)\left(\sum_{r=0}^k
 \alpha_r\phi^r+\sum_{r=k+1}^{\infty}\alpha_r\phi^r\right).
\end{eqnarray*}
Koeffizientenvergleich zeigt $a_r=\alpha_r$ f"ur $0\leq r\leq k$.

Man betrachte nun die Entwicklung in der Spitze $1$. Mit der lokalen
Koordinate $q^{1/24}$ erh"alt man unter Verwendung von (\ref{chi12})
f"ur den Charakter
$$e^{2c\cdot \frac{2 \pi i}{48}}\widetilde{\chi}_V=
e^{2c\cdot \frac{2 \pi i}{48}}
\sum_{r=0}^{k}a_k\widetilde{\chi}_{1/2}^{2c-24r}\!=\!
q^{\frac{s}{24}}(\sqrt{2})^s \left( (-1)^ka_k(1+sq)+(-1)^{k-1}a_{k-1}
(\sqrt{2})^{24}q+O(q^2)\right),$$
wobei $s=2c-24[\frac{c}{8}]<0$. Zusammenfassen ergibt
\begin{equation}\label{fuenfs}
 =\underbrace{(-1)^k(\sqrt{2})^s a_k}_{B_0}\cdot q^{\frac{s}{24}}+
\underbrace{(-1)^{k-1}(\sqrt{2})^s\left(-a_k s+
a_{k-1}\cdot2^{12}\right)}_{B_1}\cdot q^{\frac{s}{24}+1}+\cdots. 
\end{equation}

F"ur $c\geq 48$ und $r=k-1=[\frac{c}{8}]-1$ oder $r=k=[\frac{c}{8}]$ ist
$24 r -2c -3\geq 0$, und $\chi_{1/2}^{24r-2c-3}$ hat daher positive
Koeffizienten. 
Da $\chi_{1/2}^3$ der Charakter der \SVOA $\VF^{\otimes 3}$ vom Rang 
$\frac{3}{2}$ ist, hat $\frac{\chi_{1/2}^3}{\chi_{M_c}}=q^{c/24-3/48}
\left(1+\frac{1}{1-q}\left(\sum_{i\geq 1/2}P_i\, q^i\right)\right)$ nach
(\ref{partgeneral}) positive Koeffizienten.
Es folgt, da"s $\chi_{M_c}^2\cdot\left(\frac{\chi_{1/2}^3}{\chi_{M_c}}\right)'
=3\cdot \chi_{M_c}\chi_{1/2}^2\chi_{1/2}'-\chi_{M_c}'\chi_{1/2}^3$
positive und $\chi_{M_c}'\chi_{1/2}^3-2c\cdot \chi_{M_c}\chi_{1/2}^2\chi_{1/2}'$
f"ur $c\geq \frac{3}{2}$ negative Koeffizienten hat.
Wegen (\ref{viers}) sind also $a_k=\alpha_k$ und $a_{k-1}=\alpha_{k-1}$
f"ur $c\geq 48$ beide negativ, d.h.~in (\ref{fuenfs}) ist entweder $B_0$
oder $B_1$ negativ. Dies vervollst"andigt den Beweis von Satz~\ref{extsvoasatz}.
\qed

Alternativ kann man f"ur $8<c<112$ die Charaktere explizit ausrechnen
und dann zeigen, da"s der Koeffizient $A_{k+2}$ ab $c\geq  112$ negativ
wird.

\medskip
Eine diesem Abschnitt gewisserma"sen komplement"are Fragestellung ist die
Untersuchung von selbstdualen \SVOAs $V$ mit m"oglichst gro"sem Minimalgewicht
von $V_{(2)}$ bzw.~$V_{(2)}\oplus V_{(3)}$ (dem "`Schatten"').
Man erh"alt auch hier wieder Resultate~(s.~\cite{Ho-shadow}), deren
Analoga f"ur Codes und Gitter k"urzlich von N.~Elkies~\cite{El-Z,El-shadow}
beschrieben worden sind. 
Insbesondere ergibt sich eine weitere Charakterisierung 
der Babymonster-\SVOA $\VB$ sowie ---
als Folgerung --- die Ergebnisse f"ur Codes und Gitter.

\renewcommand{\baselinestretch}{1.0}
\footnotesize
\addcontentsline{toc}{chapter}{Literaturverzeichnis}
\bibliography{literatur}

\end{document}